\documentclass{amsart}

\theoremstyle{plain}

\theoremstyle{definition}

\theoremstyle{remark}

\usepackage{amssymb}

\def\cal{\mathcal}
\def\AA{{\cal A}}
\def\CC{{\cal C}}
\def\DD{{\cal D}}
\def\EE{{\cal E}}
\def\HH{{\cal H}}
\def\MM{{\cal M}}
\def\NN{{\cal N}}
\def\JJ{{\cal J}}

\def\SS{{\cal S}}

\def\LL{{\cal L}}
\def\PP{{\cal P}}

\def\WW{{\cal W}}
\def\ZZ{{\cal Z}}
\def\Z{{\mathbb Z}}
\def\C{{\mathbb C}}
\def\R{{\mathbb R}}
\def\N{{\mathbb N}}
\def\Q{{\mathbb Q}}
\def\e{{\epsilon}}

\def\n{\noindent}
\def\part{{\partial}}
\def\delbar{{\overline \partial}}
\def\dudtau{{\part u\over \part \tau}}
\def\dudt{{\part u\over \part t}}

\begin{document}

\title[Floer mini-max theory]{Floer mini-max theory, the Cerf diagram, and the spectral invariants}

\author{Yong-Geun Oh}
\address{Department of Mathematics\\University of Wisconsin\\Madison, WI
53706, ~USA \\and \\
Korea Institute for Advanced Study, Seoul, Korea}
\email{oh@math.wisc.edu}

\thanks{Partially supported by the US NSF grant \#DMS 0203593 \& 0503954}

\subjclass{53D35, 53D40} \keywords{irrational symplectic
manifolds, Hamiltonian functions, action functional, Cerf
bifurcation diagram, sub-homotopies, tight Floer cycles, handle
sliding lemma, spectral invariants, spectrality axiom}

\begin{abstract}
The author previously defined the spectral invariants, denoted by
$\rho(H;a)$, of a Hamiltonian function $H$ as the mini-max value of
the action functional $\AA_H$ over the Novikov Floer cycles in the
Floer homology class dual to the quantum cohomology class $a$. The
spectrality axiom of the invariant $\rho(H;a)$ states that the
mini-max value is a critical value of the action functional $\AA_H$.
The main purpose of the present paper is to prove this axiom for
{\it nondegenerate} Hamiltonian functions in {\it irrational}
symplectic manifolds $(M,\omega)$. We also prove that the spectral
invariant function $\rho_a: H \mapsto \rho(H;a)$ can be pushed down
to a {\it continuous} function defined on the universal ({\it
\'etale}) covering space $\widetilde{Ham}(M,\omega)$ of the group
$Ham(M,\omega)$ of Hamiltonian diffeomorphisms on general
$(M,\omega)$. For a certain generic homotopy, which we call a {\it
Cerf homotopy} $\HH = \{H^s\}_{0 \leq s\leq 1}$ of Hamiltonians, the
function $\rho_a \circ \HH: s \mapsto \rho(H^s;a)$ is piecewise
smooth away from a countable subset of $[0,1]$ for each non-zero
quantum cohomology class $a$.

The proof of this nondegenerate spectrality
relies on several new ingredients in the chain level Floer
theory, which have their own independent interest: a structure
theorem on the Cerf bifurcation diagram of the critical values of
the action functionals associated to a generic one-parameter
family of Hamiltonian functions,  a general structure theorem and
the handle sliding lemma of Novikov Floer cycles over such a
family and a {\it family version} of new transversality statements
involving the Floer chain map, and many others. We call this chain
level Floer theory as a whole the {\it Floer mini-max theory}.
\end{abstract}

\maketitle

\bigskip
\centerline{\bf Contents} \bigskip

\n \S1. Introduction and the main results
\smallskip

\n \S2. Preliminary and review of the rational case
\smallskip

{\it 2.1. Novikov Floer cycles and spectral invariants}
\par
{\it 2.2. Review of the rational case}
\smallskip

\n \S3. The Cerf bifurcation diagram
\smallskip

{\it 3.1. Generic bifurcations of the critical set of the action
functional}
\par
{\it 3.2. The Cerf bifurcation diagram}
\smallskip

\n \S4. Sub-homotopies and transversality
\smallskip

{\it 4.1 Definition of the Floer boundary map, re-visited}
\par
{\it 4.2. Definition of the Floer chain map, re-visited}
\par
{\it 4.3. Transversality of sub-homotopies}
\smallskip

\n \S5. The composition law of Floer's chain maps, re-visited
\smallskip

\n \S6. Structure of Novikov Floer cycles in a Cerf family
\smallskip

\n \S7. Handle sliding lemma and sub-homotopies
\smallskip

\n \S8. Parametric stability of tightness of Novikov Floer cycles
\smallskip

\n \S9. Proof of the nondegenerate spectrality
\smallskip

\n \S10. Spectral invariants of Hamiltonian diffeomorphisms
\smallskip

\n \S11. Applications to Hofer's geometry

\n \quad Appendix
\smallskip

\section*{\bf \S 1. Introduction and the main results}

Let $(M,\omega)$ be a closed symplectic manifold and $Ham(M,\omega)$
be the group of smooth Hamiltonian diffeomorphisms as a subgroup of
symplectic diffeomorphisms $Symp(M,\omega)$ with the smooth topology
on it. In a previous paper \cite{Oh5}, to each one-periodic
time-dependent {\it normalized} Hamiltonian function $H: S^1 \times
M \to \R$, we associated a family of symplectic invariants
$\rho(H;a)$ parameterized by the quantum cohomology classes $0 \neq
a \in QH^*(M)$, which we call the {\it spectral invariants} of $H$.
These are the mini-max values of the action functional
$$
\AA_H: \widetilde \Omega_0(M) \to \R
$$
$$
\AA_H([\gamma ,w]) = - \int w^*\omega - \int_0^1 H(t, \gamma(t))\,
dt
$$
over the semi-infinite Floer cycles in the homology class `dual' to
the quantum cohomology class $a$. (See \cite{Oh5} for a precise
meaning of the notion of the `dual' used here.) Here $\Omega_0(M)$
is the set of contractible loops on $M$ and $\widetilde \Omega_0(M)$
is the standard $\Gamma$-covering space \cite{HS}: Two pairs
$(\gamma, w)$ and $(\gamma, w')$ with $w, \, w':D^2 \to M$
satisfying $\part w =
\part w' = \gamma$ are equivalent if they satisfy
$$
\omega(w' \# \overline w) = 0 \quad \text{and } \, c_1(w' \#
\overline w) = 0. \leqno(1.1)
$$
We denote by $[\gamma,w]$ the equivalence class of $(\gamma,w)$
and by $\widetilde \Omega_0(M)$ the set of equivalence classes
$[\gamma,w]$. We provide $\widetilde\Omega_0(M)$ with the quotient
topology induced from the natural $C^\infty$ topology of the set
of pairs $(z,w)$. We denote the (\'etale) covering group of
$\pi:\widetilde \Omega_0(M) \to \Omega_0(M)$ by $\Gamma$, which is
defined by
$$
\Gamma = \frac{\pi_2(M)}{\ker \omega|_{\pi_2(M)} \cap \ker
c_1|_{\pi_2(M)}}
$$
and denote
$$
\text{Spec}(H) = \{\AA_H([z,w]) \mid \dot z = X_H(z) \}
$$
which is nothing but the set of critical values of the action
functional $\AA_H$.

The spectral invariants $\rho(H;a)$ can be regarded as the
invariants of the associated Hamiltonian paths
$$
\lambda = \phi_H: t \mapsto \phi_H^t.
$$
We say that two Hamiltonians $H, \, H' \mapsto \phi$ are
equivalent if the corresponding Hamiltonian paths $\phi_H$ and
$\phi_{H'}$ are path homotopic in $Ham(M,\omega)$. We denote by
$[\phi,H]$ its path homotopy class, and by
$\widetilde{Ham}(M,\omega)$ the set of path homotopy classes.
$\widetilde{Ham}(M,\omega)$ is the universal covering space.

We assume that all Hamiltonian functions $H$ are normalized by the
condition
$$
\int H_t \, d\mu = 0 \quad \text{for all $t \in [0,1]$} \leqno(1.2)
$$
and denote by
$$
\HH_m(M)= C^\infty_m(S^1 \times M)
$$
the set of normalized one-periodic functions. Here `$m$' stands
for the term `mean zero'. This set has one-one correspondence with
the set
$$
\PP(Ham(M,\omega),id)
$$
of Hamiltonian paths in $Ham(M,\omega)$ based at the identity.

We recall that a symplectic manifold $(M,\omega)$ is called {\it
rational} if its period group
$$
\Gamma_\omega:= \omega(\Gamma) = \{ \omega(A) \mid A \in \pi_2(M)
\}
$$
is discrete. In \cite{Oh5}, for the {\it rational} symplectic
manifold $(M,\omega)$, we have proved
$$
\rho(H;a) = \rho(H';a) \leqno(1.3)
$$
when $H, \, H' \mapsto \phi$ and $[\phi,H] = [\phi,H']$,
irrespective of nondegeneracy of Hamiltonians. In particular, the
function $\rho$ induces a well-defined function
$$
\rho:\widetilde{Ham}(M,\omega) \times QH^*(M) \to \R \leqno(1.4)
$$
by setting $\rho(\widetilde \phi;a):=\rho(H;a)$. Our proof in
 \cite{Oh5} of (1.3) for the rational case relies on the following 5
facts:

\begin{itemize} \item[(1)] The set $\text{Spec}(H) \subset \R$, which is the
set of critical values of the action functional $\AA_H$ is a set of
measure zero (see \cite[Lemma 2.2]{Oh3}).

\item[(2)] For any two Hamiltonian functions $H, \, H'\mapsto \phi$
such that $[\phi, H] = [\phi, H']$, we have
$$
\text{Spec}(H) = \text{Spec}(H') \leqno(1.5)
$$
as a subset of $\R$ provided $H, \, H'$ satisfy the normalization
condition (1.2) (see \cite{Oh3} for the proof).

\item[(3)] The function $H \mapsto \rho(H;a)$ is continuous with
respect to the smooth topology on $C^\infty_m(S^1 \times M)$ (see
 \cite{Oh5} for its proof).

\item[(4)] {\bf (Spectrality Axiom)} For any $H$ and $a$, we have
$$
\rho(H;a) \in \text{Spec}(H). \leqno(1.6)
$$

\item[(5)] The only continuous functions on a connected space (e.g.,
the interval $[0,1]$) to $\R$, whose values lie in a measure zero
subset, are constant functions.
\end{itemize}

The author previously proved the facts (1)-(3) for the general cases
in a series of papers \cite{Oh2, Oh4, Oh5}, especially including the
irrational cases. However we were able to prove the spectrality
axiom (4) only for the rational case at the time of writing the
paper \cite{Oh5}. In this paper, we prove this spectrality axiom for
{\it nondegenerate} Hamiltonian functions in the general symplectic
manifolds, especially for the irrational symplectic manifolds.
\medskip

{\noindent\bf Theorem I (Nondegenerate spectrality).} {\it Let $(M,\omega)$
be an arbitrary closed symplectic manifold. For any one-periodic
nondegenerate Hamiltonian function $H: S^1 \times M \to \R$,
$\rho(H;a)$ is a critical value of $\AA_H$, i.e.,
$$
\rho(H;a) \in \text{\rm Spec}(H)
$$
for any given quantum cohomology class $0 \neq a \in QH^*(M)$.
}
\medskip

We {\it cannot} prove the spectrality axiom for general degenerate
Hamiltonian functions.  We suspect that this is indeed not the
case when degeneracy of the Hamiltonian function is severe. It is
an interesting problem to prove or to find a counter example of
the spectrality axiom.

In fact, we prove the following stronger theorem for the
nondegenerate case.
\medskip

{\noindent\bf Theorem II (Homological essentialness).} {\it
Suppose that $H$ is nondegenerate and $a$ be any non-zero quantum
cohomology class. Then $\rho(H;a)$ is a homologically essential
critical value of $\AA_H$, in that there is a $($and so any$)$
generic choice of $J$ such that there is a Novikov Floer cycle
$\alpha$ of $(H,J)$ satisfying $[\alpha] = a^\flat$ and
$$
\rho(H;a) = \lambda_H(\alpha).
$$
}
\medskip

We refer to section 2 for the definition of a Novikov Floer cycle
$\alpha$ of the pair $(H,J)$ and its level $\lambda_H(\alpha)$.

Theorem I, which we call the {\it nondegenerate spectrality axiom},
is an essential ingredient to push down the function $\rho_a: H
\mapsto \rho(H;a)$ to a {\it continuous} function on the universal
covering space $\widetilde{Ham}(M,\omega)$. We refer to
\cite{Oh3}-\cite{Oh6} for a more detailed study of spectral
invariants and their applications.

Because we do {\it not} know the validity of spectral axiom for
degenerate Hamiltonians, the scheme of the proof used to prove
(1.3) for the rational $(M,\omega)$ cannot be applied to
degenerate Hamiltonians. In this regard, the following result is
still a non-trivial theorem to prove.
\medskip

{\noindent\bf Theorem III (Homotopy invariance).} {\it For any pair $(H,K)$,
not necessarily nondegenerate, satisfying $H \sim K$, we have
$$
\rho(H;a) = \rho(K;a).
$$
}
\medskip

To prove Theorem III for degenerate Hamiltonians, we will use
Theorem I together with continuity property of the spectral
invariant function
$$
H \mapsto \rho(H;a)
$$
and some simple calculation of the Hamiltonian algebra in the spirit
of the proof of \cite[Lemma 5.1]{Oh2} or \cite[Theorem 5.1]{Oh5}.
This theorem immediately gives rise to the following theorem.
\medskip

{\noindent\bf Theorem IV.}  {\it Let $\widetilde{Ham}(M,\omega)$ be
the above universal $($\'etale$)$  covering space and equip it with
the quotient topology induced from $\PP(Ham(M,\omega),id)$. Then the
function
$$
\rho_a: \widetilde{Ham}(M,\omega) \to \R
$$
defined by $\rho_a(\widetilde \phi): = \rho(H;a)$ for a $H$
satisfying $[H] = \widetilde \phi$ is a continuous function in the
natural topology of $\widetilde{Ham}(M,\omega)$.
}
\medskip

In the proof \cite{Oh5} of the general spectrality axiom for the
rational case, we have used the fact, in an essential way, that the
period group $\omega(\Gamma)$ is discrete and so $\text{Spec}(H)$ is
a closed subset of $\R$. For the irrational case, the argument for
the rational case cannot be applied because $\text{Spec}(H) \subset
\R$ is {\it not a closed subset but a dense subset} (of measure
zero) of $\R$. In the classical mini-max theory (see \cite{BR} for
example) where the {\it global} gradient flow of the functional
exists, a proof of the convergence result of the mini-max sequence
heavily relies on the Palais-Smale type condition and the
deformation lemma. In our case the global flow does not exist and
the set of critical values is dense and hence {\it there is no way
to deform the space itself}. Therefore in the point of view of the
critical point theory, the action functional $\AA_H$ on an
irrational manifold $(M,\omega)$ belongs to the highly pathological
realm.

To overcome these difficulties and prove criticality of the mini-max
value $\rho(H;a)$, on an {\it irrational} symplectic manifold
$(M,\omega)$, we will {\it work with the relevant mini-max cycles}
instead of either trying to deform the whole space or trying to
prove convergence of the mini-max sequence of {\it individual
critical points}. For this purpose, in the {\it nondegenerate} case,
we use some intricate arguments involving the Novikov Floer cycles
and the Floer chain map in the context of the chain level Floer
theory. The proof in turn relies on a structure theorem of the Cerf
bifurcation diagram of the action functional, a careful
re-examination of the whole construction of the basic operators in
the Floer homology theory and the composition law of the Floer chain
map, and a structure theorem of general Novikov Floer
cycles in a generic one-parameter family of Hamiltonian functions.
This latter structure theorem, Theorem 6.7, is
closely related to the picture arising in the First Cancellation
Theorem in the classical Morse theory (see \cite{Mi} and section 5
for more explanations).

As a byproduct of our proof, we also prove the following piecewise
smoothness of the spectral invariants under a {\it Cerf homotopy}
of Hamiltonians (See Definition 3.8 for the definition of Cerf
homotopy and Theorem 9.5 for a more precise description of
non-differentiable points.)
\medskip

{\noindent \bf Theorem V. } {\it Let $\HH = \{H^s\}_{0 \leq s \leq 1}$ be
a smooth Cerf homotopy of Hamiltonians. Then the function
$$
s \in [0,1] \mapsto \rho(H^s;a)
$$
is continuous piecewise-smooth with a countable number of non-differentiable
points.}
\medskip

One may hope to use some limiting argument to study the degenerate
cases as in the rational case. (See section 2 for the limiting
arguments.) This attempt meets a serious difficulty in the
irrational case. The difficulty in proving existence of such a
critical point $[z,w]$ by a limiting argument, does not lie in the
part of periodic orbit $z$ but lies in the part $w$, because the
structure of the covering group $\Gamma$ or rather its image
$\Gamma_\omega \subset \R$, the period group of $(M,\omega)$, of
the evaluation homomorphism $[\omega]: \Gamma \to \R$ is much more
complex for the irrational case than for the rational case.

Besides the proof of the main theorem, in the course of our proof of
the spectrality axiom for an irrational $(M,\omega)$, we prove many
basic results in the Floer theory itself which seem to touch the
heart of the chain level theory in the way Smale's handle body
theory does in the classical Morse theory \cite{Mi}. We expect that
these will bring further consequences to symplectic topology in the
future.

It is instructive to compare the ways how we maneuver the Floer
cycles in \cite{Oh2} and in the present paper. One of the important
matters in our chain level theory is to transfer a Floer cycle of
one Hamiltonian to that of the other Hamiltonian along a given
homotopy $\HH = \{H(\eta)\}_{0 \leq \eta \leq 1}$ of Hamiltonian
functions in a way that we can control the amount of `sliding of
handles', i.e., so that we can control the levels of the cycles
during the transfer. For this purpose, we used a piecewise linear
homotopy in \cite{Oh2}, which we call the {\it adiabatic homotopy}.
On the other hand, in the present paper, we use {\it sub-homotopies}
of $(\HH,j)$. (See Definition 3.9 for the definition of
sub-homotopy). For a given homotopy
$$
(\HH,j): [0,1] \to C^\infty([0,1],\HH_m(M)) \times C^\infty([0,1],
j_\omega)
$$
we call $(\HH_{\eta\eta'},j_{\eta\eta'})$ a {\it sub-homotopy} of
$(\HH,j)$ for $0\leq \eta, \, \eta \leq 1$ where the latter is
defined by the homotopy
$$
\HH_{\eta\eta'}: s  \mapsto H((1-s)\eta + s\eta'),\quad
j_{\eta\eta'}: s \mapsto j((1-s)\eta + s\eta')
$$
for $s \in [0,1]$. Although this difference is marginal in the
rational case, it turns out to be an essential improvement for the
irrational case. In fact, our usage of sub-homotopies is one
essential ingredient that enables us to prove certain continuity
property of the levels of the transferred Floer cycles along the
given homotopy $\HH$. (See Proposition 8.4 and Lemma 8.5.) Such a
continuity property would have been impossible to prove if we had
used the adiabatic homotopy as in \cite{Oh2} {\it due to the fact
that $\text{Spec}(H)$ is a dense subset of $\R$}: Since
$\text{Spec}(H)$ is dense, however fine a partition of $[0,1]$ we
would choose, the levels of cycles could slide down to a lower
level in each step of transfer along the piecewise-linear
homotopy. As a result, we would not have been able to control the
amount of handle sliding in the end if we used the
piecewise-linear homotopy which only approximates the given
homotopy $\HH$. In hindsight, this is a natural thing to do in the
point of view of Hamiltonian fibrations, which we will elaborate
further elsewhere. However since we use sub-homotopies of a fixed
homotopy of the pair $(\HH,j)$, we have much less freedom to
achieve transversality than the general case unlike when we use
the adiabatic homotopies in \cite{Oh2} and so need to prove the
transversality required to study such sub-homotopies. This leads
us to a novel transversality statement (Theorem 4.6) which has not
appeared in the previous literature and which is one of the
essential ingredients in our proof.

A brief description of contents of the paper is in order. After a
review of the proof for the rational case from \cite{Oh5} in section
2, the rest of the paper deals with general symplectic manifolds,
which are not necessarily rational. In section 3, we prove a
structure theorem of the {\it Cerf bifurcation diagram} of the
action functionals for a generic one-parameter family of
Hamiltonians connecting two nondegenerate Hamiltonians $H_1$ and
$H_2$. In section 4, we study the transversality issue of {\it
sub-homotopies} of the given homotopy or a one-parameter family
$(\HH,j)$ consisting of Hamiltonians and almost complex structures,
and prove the main theorem, Theorem 4.6, in that regard. In section
5, we carefully describe the gluing construction needed for the
proof of the {\it composition law} of the Floer's {\it chain} maps
in a way that will be used in our proof. In section 6, we prove a
structure theorem of Novikov Floer cycles over a one-parameter
family of Hamiltonian functions. In section
7, we recall the handle sliding lemma from \cite{Oh2} and provide a
version thereof in terms of  the sub-homotopy, instead of the
adiabatic homotopy used in \cite{Oh2}. In section 8, we prove a
general parametric stability theorem of the tightness of Novikov
Floer cycles under a Cerf homotopy. As a by-product of the this proof,
we derive Theorem V too. After all these preparations, we
carry out the proof of Theorem II (and hence of Theorem I as a
corollary) for the nondegenerate case in section 9. In section 10,
we explain how we push down the spectral invariants to
$\widetilde{Ham}(M,\omega)$ as a continuous function in the natural
topology on it, and in section 11 we provide two immediate
applications to Hofer's geometry of the Hamiltonian diffeomorphism
group.

To make the main stream of arguments transparent without bogging
down with technicalities involved with the transversality problem,
as in \cite{Oh5, Oh6}, we will assume that $(M,\omega)$ is strongly
semi-positive, i.e., satisfies the condition that there is no
spherical homology class $A$ satisfying
$$
\omega(A) > 0 \quad \text{and } \, 2 -n \leq c_1(A) < 0.
$$
We will remove this assumption imposed in \cite{Oh5, Oh6} and in the
current paper all at once in a sequel to these papers. The
spectrality axiom was first announced in the author's preprint
``Mini-max theory, spectral invariants and geometry of the
Hamiltonian diffeomorphism group'' (ArXiv:math.SG/0206092) in
general. Since then, three papers, \cite{Oh5, Oh6} and the present
paper, have grown out of the paper. Except the spectrality axiom for
degenerate Hamiltonians in irrational symplectic manifolds and the
treatment of the case that is not strongly semi-positive, all the
results in this preprint are now proved  in the three papers with
complete details and some corrections.

Finally we would like to point out that the spectrality axiom, or
(1.3), is a crucial ingredient in Entov's work \cite{En} in his
applications of spectral invariants to the study of the
quasimorphisms and the commutator length of Hamiltonian
diffeomorphisms. A proof of the spectrality axiom for the
nondegenerate case is outlined in \cite[section 3]{En}. (See Part 4
of the page 76 of \cite{En}.) We emphasize that {\it the standard
argument used in the outline cannot be applied to the irrational
symplectic manifolds} because the set of critical values of the
action functional is dense and the argument of `pushing down the
cycles under the flow' cannot easily go through in the irrational
case, as we had mentioned before.

Recently Usher \cite{usher} gave an algebraic proof of the main
theorem, Theorem II in a general abstract context of Floer homology.

We thank the Korea Institute for Advanced Study for providing the
financial support and excellent research environment during the
writing of the present paper.  We also thank M. Usher for pointing
out an incorrect statement in Theorem 3.7 (3) in the previous
version of this paper and sending his preprint \cite{usher} : This
results in our modification of the proof of Theorem 8.3 and that of
Step 2 and 3 in section 9 from the previous proofs thereof. Our
usage of Proposition 8.8 is partly influenced by Usher's paper
\cite{usher} where a similar statement is an important ingredient in
his algebraic proof of Theorem II. We also thank him for many
helpful comments during our preparation of this version of the
paper.

\bigskip

\bigskip

\n{\bf Notations}
\medskip \begin{itemize}

\item[(1)] $\JJ_\omega = \text{the set of $\omega$-compatible almost
complex structures}$

\item[(2)] $J = \{J_t\}_{0 \leq t \leq 1}$ is a smooth one-periodic
family with $J_0 = J_1$, and $j_\omega$ is the set of such $J$'s.
We just denote $j_\omega = C^\infty(S^1,\JJ_\omega)$.
\smallskip

\item[(3)] $\PP(j_\omega) = C^\infty([0,1],j_\omega)$. We denote by $j$
a general element of $\PP(j_\omega)$.
\smallskip

\item[(4)] $H: S^1 \times M \to \R$ is a one-periodic family of
functions that satisfy the normalization condition $\int_M H_t \,
d\mu= 0$ where $d\mu$ is the Liouville measure. We denote by
$\HH_m(M) = C^\infty_m(S^1\times M)$ the set of such $H$'s.
\smallskip

\item[(5)] $\HH = \{H(\eta)\}_{0 \leq \eta \leq 1}$ is a one-parameter
family of $H$'s mentioned in (4). We denote the set of such
$\HH$'s by $\PP(\HH_m(M)) = C^\infty([0,1],\HH_m(M))$
\smallskip

\item[(6)] $\Gamma_\omega := \omega(\Gamma)$= the period group of
$(M,\omega)$.
\end{itemize}

\section*{\bf \S 2. Preliminary and review of the rational case}

\subsection*{\it 2.1. Novikov Floer cycles and spectral invariants}

\medskip

Suppose that $\phi \in Ham(M,\omega)$ is  nondegenerate in the
sense of Lefshetz fixed point theory: the derivative $T_p\phi:
T_pM \to T_pM$ has no eigenvalue one at any fixed point $p \in M$.
We will call a one-periodic Hamiltonian $H: S^1 \times M \to \R$
{\it nondegenerate} if $\phi_H^1 = \phi$ is a nondegenerate
diffeomorphism. Note that the nondegeneracy of $H$ depends only on
its time-one map $\phi_H^1$. We denote by $\text{Per}(H)$ the set
of contractible one-periodic orbits of $H$.

For each nondegenerate $H:S^1 \times M \to \R $, we consider the
free $\Q$ vector space over
$$
\text{Crit}\AA_H = \{[z,w]\in \widetilde\Omega_0(M) ~|~ z \in
\text{Per}(H)\}. \leqno(2.1)
$$
Following \cite{Fl3,HS}, we give the following definition.
\medskip

{\noindent\bf Definition 2.1.} We call the formal sum
$$
\beta = \sum _{[z, w] \in \text{Crit}\AA_H} a_{[z, w]} [z,w], \,
a_{[z,w]} \in \Q \leqno(2.2)
$$
a {\it Novikov Floer chain} (or simply a {\it Floer chain}) if
there are only finitely many non-zero terms in the expression
(2.2) above any given level of the action. We call $[z,w]$
 a {\it generator} of the chain $\beta$ and denote
$$
[z,w] \in \beta
$$
if $a_{[z,w]} \neq 0$. We also say that $[z,w]$ {\it contributes
to} $\beta$ in that case. We denote by $CF_k(H)$ the set of Floer
chains whose generators all have degree $k$, i.e., satisfy
$$
\mu_H([z,w]) = k
$$
where $\mu_H$ is the Conley-Zehnder index of $[z,w]$ \cite{CZ}.
\medskip

Note that $CF_*(H)$ is a graded $\Q$-vector space. This is
infinite dimensional as a $\Q$-vector space in general, unless
$\pi_2(M) = 0$.

Now we consider a Floer chain
$$
\beta = \sum a_{[z,w]} [z,w], \quad a_{[z,w]} \in \Q.
$$
The following notion plays an essential role for the mini-max
argument via the Floer homology theory in \cite{Oh2, Oh5}.
\medskip

{\noindent\bf Definition 2.2.}~  Let $\beta$ be a Floer chain of
a given degree $k$. We define the {\it level} of the cycle $\beta$
and denote by
$$
\lambda_H(\beta) =\max_{[z,w]} \{\AA_H([z,w]) ~|~a_{[z,w]}  \neq
0\, \text{ in }\, (2.2) \} \leqno(2.3)
$$
if $\beta \neq 0$, and just put $\lambda_H(0) = -\infty$ as usual.
We call any element $[z,w]$ with $\AA_H([z,w]) = \lambda_H(\beta)$
a {\it peak} of $\beta$.
\medskip

The level $\lambda_H$ induces a filtration of $CF_*(H)$ and so
induces a natural non-Archimedean topology. (See
\cite[Appendix]{Oh5}.) We regard each $CF_k(H)$ as a topological
vector space with respect to this topology.

For a given one-periodic family $J = \{J_t\}_{0 \leq t \leq 1}$ of
compatible almost complex structures, we consider the Floer
boundary map
$$
\part= \part_{(J,H)}: CF_*(H) \to CF_*(H).
$$
In section 4, we will briefly review construction of $\part$ in a
way that is useful for our formulation of transversality problem of
sub-homotopies. One can easily check that $\part$ or all the natural
operators arising in the Floer complex are continuous with respect
to the above mentioned topology. We refer to \cite[Appendix]{Oh5}
for a precise description of the topology and for the proof of this
continuity statements.
\medskip

{\noindent\bf Definition 2.3.}~ We say that a Floer chain $\beta \in CF(H)$
is a {\it Floer cycle} if $\part \beta = 0$ and a {\it Floer
boundary} if $\beta = \part \delta$ for a Floer chain $\delta$.
Two Floer chains $\beta, \, \beta'$ are said to be {\it
homologous} if $\beta' - \beta$ is a boundary.
\medskip

Let
$$
a = \sum a_A q^{-A}, \quad a_A \in H^*(M)
$$
be a non-zero quantum cohomology class. We denote by $\Gamma(a)
\subset \Gamma$ the set of $A$'s for which the coefficient $a_A$
is non-zero. By the definition of the Novikov ring, we can
enumerate $\Gamma(a)$ so that
$$
- \lambda_1 < - \lambda_2 < \cdots  <- \lambda_j < \cdots
$$
where $\lambda_j = \omega(A_j)$. We call the first term $a_1
q^{-A_1}$ the {\it leading order term} of the quantum cohomology
class $a$.

Next, for each given quantum cohomology class $0 \neq a \in
QH^*(M)$, we consider the Floer homology class $a^\flat$ dual to
$a \in QH^*(M)$. (See [O4] for its precise meaning.) Then we
associate the following mini-max value of the action functional
$\AA_H$
$$
\rho(H;a) = \inf_{\alpha}\{\lambda_H(\alpha) \mid \alpha \in
\ker\part_H \subset CF_n(H)\, \text{with }\, [\alpha] = a^\flat \}
\leqno(2.4)
$$
to each given pair $(H,a)$. We like to emphasize that {\it this
definition itself manifests neither finiteness of the mini-max value
nor its spectral property that $\rho(H;a)$ is a critical value of
$\AA_H$}. The finiteness was proved in \cite{Oh5} for a general
symplectic manifold $(M,\omega)$ whether it is rational or not.
However, we were able to prove the spectral property only for the
rational case at the time of writing \cite{Oh5}, which we now
review.

\medskip
\subsection*{\it 2.2. Review of the rational case}
\medskip

We first recall an important notion of {\it canonical thin
cylinder} between two nearby loops. We denote by $J_{ref}$ a fixed
compatible almost complex structure and by $\exp$ the exponential
map of the metric
$$
g: = \omega(\cdot, J_{ref}\cdot).
$$
Let $\iota(g)$ be the injectivity radius of the metric $g$. As
long as $d(x,y) < \iota(g)$ for the given two points of $M$, we
can  write
$$
y = \exp_x(\xi)
$$
for a unique vector $\xi  \in T_xM$. As usual, we write the unique
vector $\xi$ as
$$
\xi = (\exp_x)^{-1}(y). \leqno(2.5)
$$
Therefore if the $C^0$ distance $d_{C^0}(z,z')$ between the two
loops
$$
z, \, z': S^1 \to M
$$
is smaller than $\iota(g)$, we can define the canonical map
$$
u^{can}_{zz'}: [0,1] \times S^1 \to M
$$
by
$$
u^{can}_{zz'}(s,t) = \exp_{z(t)}(\xi_{zz'}(t)),\quad\text{or
}\quad \xi_{zz'}(t) = (\exp_{z(t)})^{-1}(z'(t)). \leqno(2.6)
$$
It is important to note that the image of $u^{can}_{zz'}$ is
contained in a small neighborhood of $z$ (or $z'$), and uniformly
converges to $z_\infty$  when $z$ and $z'$ converge  to a loop
$z_\infty$ in the $C^1$ topology. Therefore $u^{can}_{zz'}$ also
picks out s a canonical homotopy class, denoted by
$[u^{can}_{zz'}]$, among the set of homotopy classes of the maps
$u: [0,1] \times S^1 \to M$ satisfying the given boundary
condition
$$
u(0,t) = z(t), \quad u(1,t) = z'(t).
$$
The following lemma is an important ingredient in our proof, which
will be used to overcome irrationality of symplectic manifolds
later.
\medskip

{\noindent\bf Lemma 2.4.} {\it Let $z, \, z': S^1 \to M$ be two smooth loops
and $u^{can}$ be the above canonical cylinder. Then as
$d_{C^1}(z,z') \to 0$, then the map $u^{can}_{zz'}$ converges in
the $C^1$-topology, and its geometric area $Area(u^{can})$
converges to zero. In particular, we have the followings:
\begin{itemize}

\item[(1)]  For any bounding disc $w$ of $z$, the bounding disc
$$
w' : = w \# u^{can}_{zz'}
$$
of $w'$ is pre-compact in the $C^1$-topology of the maps from the
unit disc.

\item[(2)]
$$
\int_{u^{can}_{zz'}} \omega \to 0 \leqno(2.7)
$$
as $d_{C^1}(z,z') \to 0$ as $z' \to z$.
\end{itemize}
}

\begin{proof} (1) is an immediate consequence of the explicit form
(2.6) of $u_{zz'}^{can}$ and from the standard property of the
exponential map.

On the other hand, from the explicit expression (2.6) of the
canonical thin cylinder and from the property of the exponential
map, it follows that the geometric area
$\text{Area}(u_{i\infty}^{can})$ converges to zero as
$d_{C^1}(z,z') \to 0$ by an easy area estimate. Since $z,\, z'$
are assumed to be $C^1$, it follows $u_{zz'}^{can}$ is $C^1$ and
hence the inequality
$$
\text{Area}(u_{i\infty}^{can}) \geq
\Big|\int_{u_{i\infty}^{can}}\omega\Big|.
$$
This implies
$$
\lim_{j \to \infty} \int_{u_{i\infty}^{can}}\omega = 0,
$$
which finishes the proof.
\end{proof}

The following theorem was previously proved by the author in
\cite{Oh5}. We duplicate its proof here to highlight differences
between the rational and the irrational cases, and to motivate the
scheme of our proof in the irrational case.
\medskip

{\noindent\bf Theorem 2.5. \cite[Theorem 7.1]{Oh5}.} {\it Suppose
that $(M,\omega)$ is rational. Then for any smooth one-periodic
Hamiltonian function $H: S^1 \times M \to \R$, we have
$$
\rho(H;a) \in \text{\rm Spec}(H)
$$
for each given quantum cohomology class $0 \neq a \in QH^*(M)$.
}

\begin{proof}  We need to show that the mini-max value $\rho(H;a)$
is a critical value, or that there exists $[z,w] \in \widetilde
\Omega_0(M)$ such that
$$
\begin{aligned}
& \AA_H([z,w]) = \rho(H;a) \\
& d\AA_H([z,w]) = 0, \quad \text{i.e., } \quad \dot z = X_H(z).
\end{aligned}
$$
The finiteness of the value  $\rho(H;a)$ was proved in \cite{Oh5}.
If $H$ is nondegenerate, we just use the fixed Hamiltonian $H$. If
$H$ is degenerate, we approximate $H$ by a sequence of nondegenerate
Hamiltonians $H_i$ in the $C^2$ topology. Let $[z_i,w_i] \in
\text{Crit}\AA_{H_i}$ be a peak of the Floer cycle $\alpha_i \in
CF_*(H_i)$,  such that
$$
\lim_{j \to \infty} \AA_{H_i}([z_i,w_i]) = \rho(H;a). \leqno(2.8)
$$
Such a sequence can be chosen by the definition of $\rho(\cdot;
a)$ and its finiteness property.

Since $M$ is compact and $H_i \to H$ in the $C^2$ topology,  and
$\dot z_i = X_{H_i}(z_i)$ for all $i$, it follows from the
standard boot-strap argument that $z_i$ has a subsequence, which
we still denote by $z_i$, converging to some loop $z_\infty: S^1
\to M$ satisfying $\dot z = X_H(z)$. Now we show that the sequence
$[z_i,w_i]$  are pre-compact on $\widetilde \Omega_0(M)$.  Since
we fix the quantum cohomology class $0 \neq a \in QH^*(M)$ (or
more specifically since we fix its degree) and since the Floer
cycle is assumed to satisfy $[\alpha_i] = a^\flat$, we have
$$
\mu_{H_i}([z_i,w_i]) = \mu_{H_j}([z_j,w_j]).
$$

{\noindent\bf Lemma 2.6.} {\it When $(M,\omega)$ is rational, $\text{Crit}
\AA_K \subset \widetilde \Omega_0(M)$ is a closed subset of $\R$
for any smooth Hamiltonian $K$, and is locally compact in the
subspace topology  of  the covering space
$$
\pi: \widetilde \Omega_0(M) \to \Omega_0(M).
$$
}
\begin{proof} First note that when $(M,\omega)$ is rational, the
covering group $\Gamma$ of $\pi$ above is discrete. Together with
the fact that the set of solutions of $\dot z = X_K(z)$ is compact
(on compact M),  it follows that
$$
\text{Crit}(\AA_K) =  \{[z,w] \in \widetilde\Omega_0(M)
 \mid \dot z = X_K(z) \}
$$
is a closed subset which is also locally compact.
\end{proof}

Now consider the bounding discs of $z_\infty$ given by
$$
w_i' = w_i \# u_{i\infty}^{can}
$$
for all sufficiently large $i$, where $u_{i\infty}^{can} =
u^{can}_{z_iz_\infty}$ is the canonical thin cylinder between
$z_i$ and $z_\infty$. We note that as $i \to \infty$ the geometric
area of $u_{i\infty}^{can}$ converges to $0$.

We compute the action of the critical points $[z_\infty, w_i'] \in
\text{Crit}\AA_H$,
$$
\begin{aligned} \AA_H([z_\infty, w_i'])  & = - \int_{w_i'} \omega -
\int_0^1 H(t, z_\infty(t)) \, dt
\\
 & = - \int_{w_i} \omega - \int_{u_{i\infty}^{can}} - \int_0^1
H(t, z_\infty(t)) \, dt
\\
 & = \Big(- \int_{w_i} \omega  - \int_0^1 H_i(t, z_i(t)) \, dt\Big)
 \end{aligned}\leqno(2.9)$$
$$\begin{aligned}
 & \quad - \Big(\int_0^1 H(t, z_\infty(t))
- \int_0^1 H_i(t, z_i(t)) \Big) - \int_{u_{i\infty}^{can}} \omega \\
\quad\qquad& = \AA_{H_i}([z_i, w_i]) - \Big(\int_0^1 H(t,
z_\infty(t)) - \int_0^1 H_i(t, z_i(t)) \Big)\\ &\qquad -
\int_{u_{i\infty}^{can}} \omega.
\end{aligned}
\leqno(2.10)
$$
Since $z_i$ converges to $z_\infty$ uniformly and $H_i \to H$, we
have
$$
- \Big(\int_0^1 H(t, z_\infty(t)) - \int_0^1 H(t, z_i(t)) \Big)
\to 0. \leqno(2.11)
$$
Therefore combining (2.7), (2.8) and (2.11), we derive
$$
\lim_{i \to \infty} \AA_H([z_\infty, w_i'])  =  \rho(H;a).
$$
In particular $\AA_H([z_\infty, w_i'])$ is a Cauchy sequence,
which implies
$$
\Big|\int_{w_i'} \omega - \int_{w_j'}\omega \Big| =
\Big|\AA_H([z_\infty, w_i']) - \AA_H([z_\infty, w_j'])\Big| \to 0
$$
i.e.,
$$
\int_{w_i' \# \overline w_j'} \omega \to 0.
$$
Since $\Gamma$ is discrete and $\int_{w_i' \# \overline w_j'}
\omega\in \Gamma$, this indeed implies that
$$
\int_{w_i' \# \overline w_j'} \omega  = 0 \leqno(2.12)
$$
for all sufficiently large $i, \, j \in \Z_+$. Since the set
$\Big\{\int_{w_i'} \omega\Big\}_{i \in\Z_+}$ is bounded, we
conclude that the sequence $\int_{w_i'} \omega$ eventually
stabilize, by choosing a subsequence if necessary. Going back to
(2.9), we derive that the actions
$$
\AA_H([z_\infty, w_i'])
$$
themselves stabilize and so we have
$$
\AA_H([z_\infty,w_N']) = \lim_{i \to \infty}\AA_H([z_\infty,
w_i']) = \rho(H;a)
$$
for a fixed sufficiently large $N \in \Z_+$. This proves that
$\rho(H;a)$ is indeed the value of $\AA_H$ at the critical point
$[z_\infty,w_N']$. This finishes the proof.
\end{proof}

In fact, an examination of the above proof proves a stronger fact
that the mini-maxing sequence $[z_i,w_i]$ is precompact {\it for
the rational case}, which we now explain. We recall that if $H, \,
H'$ are nondegenerate and sufficiently $C^2$-close, there exists a
canonical one-one correspondence between the set of associated
Hamiltonian periodic orbits. We call an {\it associated pair} any
pair $(z,z')$ of Hamiltonian periodic orbits of $H, \, H'$ mapped
to each other under this correspondence. We will give the proof of
the following proposition in the Appendix.
\medskip

{\noindent\bf Proposition 2.7.} {\it Suppose that $H, \, H'$ are
nondegenerate and sufficiently $C^2$ close. Let $(z,
 z')$ be an associated pair of $H, \, H'$. Then we have
$$
\mu_H([z, w]) = \mu_{H'}([z',w\# u^{can}_{zz'}]). \leqno(2.13)
$$
}

We derive
$$
\begin{aligned} 2 c_1([w_i'\# \overline w_j']) & = 2c_1([w_i \#
u_{i\infty}^{can}\# \overline{w_j \# u_{j\infty}^{can}}]) \\
& = 2c_1([w_i\# u_{i\infty}^{can}\# \overline u_{j\infty}^{can}
\# \overline w_j])\\
& =  \mu_{H_i}([z_i,w_i]) - \mu_{H_i}([z_i, w_j \#
u_{j\infty}^{can}\# \overline u_{i\infty}^{can}]).
\end{aligned}\leqno(2.14)
$$
The third equality comes from the index formula
$$
\mu_H([z,w\# A]) = \mu_H([z,w]) - 2c_1(A)
$$
(see \cite[Appendix]{Oh5} for the details of its proof). On the
other hand, we derive
$$
\mu_{H_i}([z_i,w_j \# u_{j\infty}^{can}\# \overline
u_{i\infty}^{can}]) = \mu_{H_i}([z_i, w_j\# u_{z_jz_i'}^{can}]) =
\mu_{H_j}([z_j,w_j]) \leqno(2.15)
$$
when $i, \, j$ are sufficiently large. Here the first equality
follows since $u_{i\infty}^{can} \# \overline u_{i\infty}^{can}$
is homotopic to the canonical thin cylinder $u^{can}_{z_jz_i'}$,
and the second comes from (2.13). On the other hand, $[z_i,w_i]$
and $[z_j,w_j]$ satisfy
$$
\mu_{H_i}([z_i,w_i]) = \mu_{H_j}([z_j,w_j]) \leqno(2.16)
$$
because they are generators of Floer cycles $\alpha_i$ and
$\alpha_j$ both representing the same Floer homology class
$a^\flat$ and so having the same degree. Hence combining
(2.14)-(2.16), we obtain
$$
c_1([w_i'\# \overline{w}_j'])= 0 \leqno(2.17)
$$
for all sufficiently large $i, \, j$. Combining (2.12) and (2.17),
we have proved
$$
[z_\infty,w_i'] = [z_\infty, w_j'] \quad \text{in }\,
\widetilde{\Omega}_0(M).
$$
If we denote by $[z_\infty,w_\infty]$ this common element of
$\widetilde\Omega_0(M)$, we have proven that the sequence
$[z_i,w_i]$ converges to a critical point $[z_\infty,w_\infty]$ of
$\AA_H$ in the topology of the covering space $\pi:
\widetilde\Omega_0(M) \to \Omega_0(M)$. This finishes our
discussion about the rational case.

For the irrational case, the above argument breaks down since the
sequence $[z_\infty,w_i']$ used in the above proof will not
stabilize, and more seriously the action values
$\AA_H([z_\infty,w_i'])$ may accumulate at a value in $\R
\setminus \text{Spec}(H)$. Recall that in the irrational case,
$\text{Spec}(H)$ is a dense subset of $\R$. Therefore in the
irrational case, one needs to directly prove that the sequence has
a convergent subsequence in the natural topology of
$\widetilde\Omega_0(M)$. It turns out that the above limiting
arguments used for the rational case cannot be carried out due to
the possibility that the discs $w_i$ could behave wildly in the
limiting process.  As a result, proving such a convergence is not
possible in general even for the nondegenerate case for a given
mini-max sequence of {\it critical points} $[z_i,w_i]$ satisfying
(2.8). One needs to use a mini-max sequence of {\it cycles}
instead. This scheme is exactly what we have carried out in the
present paper. Because we use the Floer cycles and they are
defined only for nondegenerate Hamiltonians, we can prove the
spectrality axiom only for the nondegenerate case in this way.
Along the way, we develop many new ingredients in the chain level
Floer theory needed to carry out the scheme. We call our chain
level theory the {\it Floer mini-max theory}.

To go to the case of degenerate Hamiltonians, it is unavoidable to
use the approximation used arguments above as in the rational
case. It would be very interesting to see if this difficulty is
something intrinsic for this case. (See Remark 9.4 for some
related comments.)

\section*{\bf 3. The Cerf bifurcation diagram}

\subsection*{\it 3.1. Generic bifurcations of the critical set of the
action functional}

We first recall that for a generic one-parameter family $\HH =
\{H(\eta)\}_{0\leq \eta\leq 1}$, there are a {\it finite} number
of points
$$
\SS ing(\HH) = \{s_1, s_2, \ldots, s_{k_1}\} \subset [0,1]
$$
where there occurs either birth-death or death-birth type of
bifurcation of periodic orbits (see \cite{Lee} for a detailed proof
of this). Furthermore at each such $s_i$, there is exactly
one-periodic orbit $z_i$ of $\dot x = X_{H(s_i)}(x)$ for which a
continuous family of the pair $z^+(\eta), \, z^-(\eta)$ of periodic
orbits of $\dot x = X_{H(\eta)}(x)$ bifurcate from $z_i$ for $\eta$
with $|\eta - s_i| < \delta$, $\delta$ sufficiently small, that
satisfy
\smallskip

(1) $z^\pm(\eta) \to z_i$ as $\eta \to s_i$,\par

(2) the Conley-Zehnder indices satisfy
$$
\mu([z^+,w^+]) = \mu([z^-,w^-]) + 1 \leqno(3.1)
$$
where $w^+\sim w^- \# u$ for $u$ the thin cylinder between $z^+$
and $z^-$. This latter condition makes sense because $z^+$ and
$z^-$ are close when $\delta$ is sufficiently small, which depends
only on the homotopy $\HH$ independent of points $s_i \in \SS
ing(\HH)$. We denote this uniform $\delta$ as
$$
\delta_1= \delta_1(\HH). \leqno(3.2)
$$
In the course of studying a detailed structure of the Cerf-type
bifurcation diagram of the action functional, we will provide an
outline of a proof of the above statements for the reader's
convenience leaving more details to \cite{Lee}.

Let $\HH = \{H(\eta)\}_{0 \leq \eta \leq 1}$ be a homotopy (or a
one-parameter family) of smooth Hamiltonians. We denote
$$
\text{Crit}\AA_{\HH}: = \bigcup_{\eta \in [0,1]} \{ [z,w] \mid
\dot z = X_{H(\eta)}(z) \, \} = \bigcup_{\eta \in [0,1]}
\text{Crit}\AA_{H(\eta)} \leqno(3.3)
$$
and consider it as a subset of $[0,1] \times
\widetilde{\Omega}_0(M)$. For the simplicity of notations, we
sometimes denote
$$
\HH_m(M): = C^\infty_m(S^1 \times M). \leqno(3.4)
$$
We then define
$$
\PP \HH_m(M): = C^\infty([0,1],\HH_m(M)), \leqno(3.5)
$$
and
$$
\PP(\HH_m(M);H_0,H_1) : = \{\HH \in \PP \HH_m(M) \mid H(0) =
H_0,\, H(1) = H_1 \} \leqno(3.6)
$$
for any given Hamiltonians $H_1, \, H_2$. We also denote
$$
\HH_m^{nd}(M): = \{ H \in \HH_m(M) \mid H \quad\text{is
nondegenerate}\}.
$$

The following lemma is easy to prove from the definition and from
the standard facts on the first order ordinary differential
equation.
\medskip

{\noindent\bf Lemma 3.1.} {\it Let $H_\alpha,\, H_\beta$ be smooth and $\HH
\in \PP(\HH_m(M);H_\alpha,H_\beta)$. Then we have \begin{itemize}

\item[(1)] $\text{Crit}\AA_{\HH}$ is invariant under the deck
transformation of $\Gamma$ on $\widetilde \Omega_0(M)
\times[0,1]$.

\item[(2)] Under the covering map $\pi:\widetilde \Omega_0(M) \to
\Omega_0(M)$, $\text{Crit} \AA_{\HH}/\Gamma$ coincides with
$$
\text{Per}(\HH): = \cup_{\eta \in [0,1]}\text{Per}(H(\eta))
\subset [0,1] \times \Omega_0(M) \leqno(3.7)
$$
and in particular is compact.
\end{itemize}
}

Now we prove the following general structure theorem on
$\text{Crit}\AA_\HH$ for a generic homotopy $\HH$ such that the
end points $H(0)$ and $H(1)$ are nondegenerate.
\medskip

{\noindent\bf Proposition 3.2.} {\it Let $H_\alpha,\, H_\beta$ be two
nondegenerate Hamiltonian functions and $\HH \in
\PP(\HH_m(M);H_\alpha,H_\beta)$. Then there exists a dense subset
$$
\PP^{reg}(\HH_m(M);H_\alpha,H_\beta) \subset \PP(\HH_m(M);H_0,H_1)
$$
such that the subset
$$
\text{Crit}\AA_{\HH} \subset [0,1] \times \widetilde \Omega_0(M)
$$
becomes a smooth one-manifold with its boundary $\part
(\text{Crit}\AA_{\HH})$ contained in $\{0,1\}$ $\times$ $\widetilde
\Omega_0(M)$.
}

\begin{proof} We first note that the action by $\Gamma$ on
$\text{Crit}\AA_\HH$ is free. Therefore it is enough to prove that
$\text{Crit}\AA_{\HH}/\Gamma$ is a smooth one-manifold for some
generic choices of $\HH$. We also know
$$
\text{Crit}\AA_{\HH}/\Gamma = \text{Per}(\HH). \leqno(3.8)
$$
This in turn implies that it is enough to prove that
$\text{Per}(\HH)$ becomes a (compact) smooth one-manifold for some
generic choices of $\HH$. We now prove this statement.

We first remark that in the Fredholm analysis we are going to carry
out below, one needs to take a suitable Banach completion of the
various function spaces that appear. However this is a standard
procedure by now, and so we will not mention this technicality but
just work with $C^\infty$ function spaces. A good reference for this
matter and also for detailed calculations involving the action
functional is the paper \cite{W} by Weinstein.

We define the map
$$
\Phi: \Omega_0(M) \times \PP(\HH_m(M);H_\alpha, H_\beta) \times
[0,1] \to T\Omega_0(M) \times [0,1]
$$
by
$$
\Phi(z,\HH, \eta) \mapsto \Big(\dot z - X_{H(\eta)}(z), \eta
\Big). \leqno(3.9)
$$
Considering this as the composition of the section
$$
(z, \HH, \eta) \to (\dot z - X_{H(\eta)}(z), \HH, \eta)
$$
of the parametric tangent bundle
$$
T\Omega_0(M) \times \PP(\HH_m(M);H_\alpha,H_\beta) \times [0,1]
\to \Omega_0(M)\times \PP(\HH_m(M);H_\alpha,H_\beta) \times [0,1]
$$
and the projection map
$$
\Omega_0(M)\times \PP(\HH_m(M);H_\alpha,H_\beta) \times [0,1] \to
\Omega_0(M) \times [0,1],
$$
it is straightforward to check that the derivative of $\Phi$ is
surjective at all the zero points $(z,\HH, \eta)$ of $\Phi$, i.e.,
those satisfying $\dot z = X_{H(\eta)}(z)$. Therefore the {\it
universal} set of periodic orbits, denoted by
$$
\PP er:= \Phi^{-1}(o_{T\Omega_0(M)} \times
\PP(\HH_m(M);H_\alpha,H_\beta) \times [0,1]), \leqno(3.10)
$$
is a smooth submanifold of $\Omega_0(M) \times
\PP(\HH_m(M);H_\alpha,H_\beta) \times [0,1]$ by the implicit
function theorem.

Furthermore it is well-known that the linearization map
$$
\xi \mapsto \frac{D\xi}{dt} - DX_K(z)(\xi) \leqno(3.11)
$$
along the periodic orbit $z$ of a Hamiltonian $K \in \HH_m(M)$ is a
Fredholm operator of index zero after making a suitable Banach
completion of $\Omega_0(M)$ (see \cite{W} for example). This then is
translated into the statement that the projection map
$$
\pi_2: \PP er \subset \Omega_0(M) \times
\PP(\HH_m(M);H_\alpha,H_\beta) \times [0,1] \to
\PP(\HH_m(M);H_\alpha,H_\beta) \leqno(3.12)
$$
is a Fredholm map of index 1. Now by the Sard-Smale theorem, the
set of regular values of $\pi_2$, which we denote by
$$
\PP^{reg}(\HH_m(M);H_\alpha,H_\beta)
$$
is residual, and in particular dense. This finishes the proof.
\end{proof}

Next a simple version of the two-jet transversality implies the
following, whose proof we omit and refer to \cite{Lee} for more
details.
\medskip

{\noindent\bf Proposition 3.3.} {\it Let $\HH \in
\PP^{reg}(\HH_m(M);H_1,H_2)$. Then the set
$$
\SS ing(\HH)\!: =\! \{ \eta \in [0,1] \mid \!\text{the linearization map
{\rm(3.11)} is not surjective for $H(\eta)$} \}
$$
is finite. And there is another smaller dense subset of $$\HH \in
\PP^{reg}(\HH_m(M);H_1,H_2)$$ for which at each point $\eta \in
\SS ing(\HH)$, either birth-death or death-birth type of
bifurcation occurs as described in the beginning of this section.
}
\medskip \subsection*{\it 3.2. The Cerf bifurcation diagram}
\medskip

Next we introduce a notion of the {\it Cerf bifurcation diagram}
of the action functionals and study its structure for a generic
choice of the homotopy $\HH$.

\medskip{\noindent \bf Definition 3.4.}
Consider the set
$$
\begin{aligned} \Sigma(\HH)  = & \{(\eta,a) \mid   \eta \in [0,1],\, a =
\AA_{H(\eta)}([z,w]), \, [z,w] \in \text{Crit}\AA_{H(\eta)} \}\\
& \subset [0,1] \times \R. \end{aligned} \leqno (3.13)
$$
We call $\Sigma(\HH)$ the Cerf bifurcation diagram of the homotopy
$\HH$.
\medskip

There is the natural evaluation homomorphism
$$
g \in \Gamma \mapsto \omega(g) \in \Gamma_\omega  \subset \R.
\leqno (3.14)
$$
Via the homomorphism (3.14), $\Gamma$ naturally acts on $[0,1]
\times \R$ by
$$
\Gamma \times [0,1] \times \R \to [0,1] \times \R ;\quad g \cdot
(\eta,a) \mapsto  (\eta, a - \omega(g)) \leqno (3.15)
$$
which preserves $\Sigma(\HH)$. Now we consider the map $\widetilde
\Phi$
$$
\begin{aligned} \widetilde \Phi: \widetilde \Omega_0(M) & \times
\PP(\HH_m(M);H_\alpha,H_\beta) \times [0,1] \\
& \to \widetilde T\Omega_0(M) \times
\PP(\HH_m(M);H_\alpha,H_\beta) \times [0,1] \times \R
\end{aligned}
$$
defined by
$$
\widetilde\Phi([z,w], \HH, \eta) = (\dot z - X_{H(\eta)}(z), \HH,
\eta, \AA_{H(\eta)}([z,w])). \leqno (3.16)
$$
This map is equivariant under the obvious actions of $\Gamma$.

It follows from Proposition 3.2 that $\widetilde \Phi$ is
transverse to the submanifold
$$
o_{T\Omega_0(M)} \times [0,1] \times \R \subset T\Omega_0(M)
\times [0,1] \times \R.
$$
In particular, we know that
$$
\ZZ : = (\widetilde \Phi)^{-1}(o_{T\Omega_0(M)} \times [0,1]
\times \R) \leqno (3.17)
$$
is a smooth submanifold of
$$
\widetilde \Omega_0(M) \times \PP(\HH_m(M);H_\alpha,
H_\beta)\times [0,1].
$$
Obviously we have the natural (\'etale) covering map
$$
\Gamma \to \ZZ \to \PP er = \bigcup_{\HH \in
\PP(\HH_m(M);H_\alpha,H_\beta)} \text{Per}(\HH). \leqno (3.18)
$$
Now we consider the projection
$$
\Pi: \ZZ \subset \widetilde \Upsilon \to
\PP(\HH_m(M);H_\alpha,H_\beta).
$$
It is easy to check that this is a Fredholm map of index 1. By the
Sard-Smale theorem, there is another dense subset set of $\HH$'s
such that
$$
\Pi^{-1}(\HH) = : Z(\HH)
$$
becomes a smooth one dimensional manifold with boundary
$$
\part Z(\HH) = \Psi^{-1}(o_{T\Omega_0(M)} \times \{0,1\} \times
\R),
$$
{\it as long as  it is non-empty}. In addition, $\Gamma$ freely
acts on $Z(\HH)$ and hence comes a natural $\Gamma$ principal
bundle
$$
\pi: Z(\HH) \to Z(\HH)/\Gamma \leqno (3.19)
$$
where the quotient can be canonically identified with
$\text{Per}(\HH)$. By the standard a priori estimates on the
Hamilton equation, $Per(\HH)$ is shown to be a compact one
manifold with boundary.

Now consider the map
$$
ev: Z(\HH) \subset \widetilde \Omega_0(M) \times [0,1] \to [0,1]
\times \R
$$
defined by
$$
([z,w], \eta) \to \Big(\eta, \AA_{H(\eta)}([z,w]) \Big) \leqno
(3.20)
$$
for which the diagram
$$
\begin{matrix} Z(\HH) &  \longrightarrow  & [0,1] \times \R \\
\downarrow & \quad & \downarrow \\
[0,1] & \equiv & [0,1]
\end{matrix}
$$
commutes. Furthermore $ev$ is equivariant under the
fiber-preserving action of $\Gamma$, and the Cerf diagram
$\Sigma(\HH)$ is nothing but the image of $ev$. Hence there
induces a natural action of $\Gamma$ acting fiberwise on
$\Sigma(\HH)$ under the projection $\Sigma(\HH) \to [0,1]$.

Now we describe structure of the image of the map (3.20). We first
prove the following lemma, which states that the action cannot be
the same for two different critical points $[z_1,w_1]$ and
$[z_2,w_2]$ of the form $z_1 =z_2$, if we further require
$\mu_K([z_1,w_1]) = \mu_K([z_2,w_2])$.

\medskip{\noindent \bf  Lemma 3.5.} {\it Let $K$ be a nondegenerate Hamiltonian
function and consider the elements from $\text{Crit}(\AA_K)$. Then
for any element $z \in \text{Per}(K)$, $\AA_K([z,w_1]) =
\AA_K([z,w_2])$ and $\mu_K([z,w_1]) = \mu_K([z,w_2])$ if and only
if $[z,w_1] = [z,w_2]$. In particular, if we have
$$
\AA_K([z_1,w_1]) = \AA_K([z_2,w_2]), \quad \mu_K([z_1,w_1]) =
\mu_K([z_2,w_2]),
$$
then $z_1 \neq z_2$ unless $[z_1,w_1] = [z_2,w_2]$. }
\begin{proof} Suppose
$$
\mu_K([z,w_1]) = \mu_K([z,w_2]) \leqno (3.21)
$$
and
$$
\AA_K([z,w_1]) = \AA_K([z,w_2]) \leqno (3.22)
$$
From the assumption (3.22), we derive
$$
\int_{w_2\# \overline w_1}\omega = 0 \leqno (3.23)
$$
since we have
$$
\AA_K([z,w_1])- \AA_K([z,w_2]) = \int_{w_1\# \overline w_2}\omega.
$$
And from (3.21) and the index formula
$$
\mu_K([z,w_1]) = \mu_K([z,w_2]) - 2 c_1(w_1\# \overline w_2)
$$
we derive
$$
c_1(w_1\# \overline w_2) = 0. \leqno (3.24)
$$
Therefore (3.23) and (3.24) imply $[z,w_1] = [z,w_2]$ by the
definition of $\Gamma$-equivalence classes in (1.1). The converse
is obvious.  \end{proof}

\medskip{\noindent \bf  Definition 3.6.} Let $H_1,\, H_2$ be two nondegenerate
Hamiltonians and $\HH$ be a homotopy between them. Let
$\Sigma(\HH) \subset [0,1] \times \R$ be the associated Cerf
diagram.

\begin{itemize}
\item[(1)] We say that a point $(\eta,a) \in \Sigma(\HH)$ is a
{\it cusp} if $\eta \in \SS ing(\HH)$ and $\AA_{H(\eta)}([z,w]) =
a$, and a {\it generic cusp} if the second derivative
$d^2\AA_{H(\eta)}([z,w])$ at $[z,w]$ has exactly one-dimensional
kernel.

\item[(2)] We say that a point $(\eta,a) \in\Sigma(\HH)$ is a {\it
crossing} if there are two different $[z_1,w_1], \, [z_2,w_2] \in
\text{Crit}\AA_{H(\eta)}$ with
$$
\begin{aligned} a = \AA_{H(\eta)}([z_1,w_1]) & = \AA_{H(\eta)}([z_2,w_2])\\
\mu_{H(\eta)}([z_1,w_1]) & = \mu_{H(\eta)}([z_2,w_2])
\end{aligned}
$$
and a {\it nondegenerate crossing} if it has the property that the
corresponding branches intersect transversely. We denote by
$$
\CC  ross^{nd}(\HH) \subset [0,1]
$$
the set of nondegenerate crossings. If there are exactly one such
nondegenerate crossing at $\eta$, modulo the action of $\Gamma$,
whose associated pair of critical points $[z_i,w_i], \, i=1, \, 2$
in addition, we call the crossing a {\it generic crossing}.
\end{itemize}
\medskip

We note that due to the action of $\Gamma$ on the Cerf diagram, the set
$\CC ross^{nd}(\HH)$ is a countable infinite subset of $[0,1]$.
We like to emphasize that {\it $\CC ross^{nd}(\HH)$ does not
include the points $(\eta,a)$ with $a = \AA_{H(\eta)}([z_1,w_1]) =
\AA_{H(\eta)}([z_2,w_2])$ with different Conley-Zehnder indices}.
There exists a natural fiberwise action of the group $\Gamma$ on
$\Sigma(\HH)$ under the projection
$$
\pi_1: \Sigma(\HH) \subset [0,1] \times \R \to [0,1].
$$
With these definitions, we prove the following structure theorem
of the Cerf bifurcation diagram for a generic homotopy $\HH$.

\medskip{\noindent \bf  Theorem 3.7.} {\it Let $H_1, \, H_2$ be two nondegenerate
Hamiltonians. Then there exists a dense subset
$$
\PP^{Cerf}(\HH_m(M);H_0,H_1) \subset \PP^{reg}(\HH_m(M);H_0,H_1)
\subset \PP(\HH_m(M);H_0,H_1)
$$
such that for any element $\HH \in\PP^{Cerf}(\HH_m(M);H_0,H_1)$
its associated Cerf diagram $\Sigma(\HH)$ satisfies the following
list of the properties: \begin{itemize}

\item[(1)] $\Sigma(\HH)$ is the projection of the one-manifold
$Z(\HH)$, with boundary contained in $\pi_1^{-1}(\widetilde
\Omega_0(M) \times \{0, 1\}\times \R)$.

\item[(2)] All the crossings are nondegenerate and unique modulo
the action of $\Gamma$.

\item[(3)] There is exactly one cusp point $(\eta,a)$ unique,
modulo the action of $\Gamma$, corresponding to each point $\eta
\in \SS ing(\HH)$.

\item[(4)] $\SS ing(\HH) \cap \CC  ross^{nd}(\HH) = \emptyset$.
\end{itemize}

}

\begin{proof} We denote the diagonal subset of $\Omega_0(M) \times
\Omega_0(M)$ by  $\Delta$, and let
$$
\widetilde \Delta: = \pi^{-1}(\Delta) \subset \widetilde
\Omega_0(M) \times \widetilde \Omega_0(M)
$$
be its lifting to $\widetilde \Omega_0(M) \times \widetilde
\Omega_0(M)$ under the projection map
$$
\widetilde \Omega_0(M) \times \widetilde \Omega_0(M) \to
\Omega_0(M) \times \Omega_0(M).
$$
We consider the spaces
$$
\widetilde \LL_2(M) := (\widetilde \Omega_0(M) \times \widetilde
\Omega_0(M)) \setminus \widetilde \Delta
$$
and
$$
\widetilde \LL_2(M) \times \PP(\HH_m(M);H_1,H_2) \times [0,1] =:
\widetilde \Upsilon.
$$
There is a natural product action by $\Gamma \times \Gamma$ on the
product $\widetilde \Omega_0(M) \times \widetilde \Omega_0(M)$,
which preserves the subset $\widetilde \LL_2(M)$.

We define the map
$$
\widetilde \Psi: \widetilde \Upsilon \to T(\widetilde \LL_2(M))
\times [0,1] \times \R
$$
by
$$
\begin{aligned} \widetilde \Psi([z_1, w_1], [z_2, w_2], \HH, \eta)= &
\Big(\dot z_1 - X_{H(\eta)}(z_1), \dot z_2 -
X_{H(\eta)}(z_2), \eta,\\
& \quad \AA_{H(\eta)}([z_1,w_1]) - \AA_{H(\eta)} ([z_2,w_2])\Big).
\end{aligned}
\leqno (3.25)
$$
Note that there is a canonical identification
$$
T(\widetilde \LL_2(M)) \times [0,1] \times \R \cong T(\LL_2(M))
\times [0,1] \times \R
$$
as before and so we assume the image of $\widetilde \Psi$ lies
$T(\LL_2(M)) \times [0,1] \times \R$. The image of the map
$\widetilde \Psi$ is invariant under the diagonal action $\Gamma$
on $\widetilde\Upsilon$ and so pushes down to the map
$$
\Psi: \Upsilon \to T(\LL_2(M)) \times [0,1] \times \R \leqno (3.26)
$$
where $\Upsilon : = \widetilde \Upsilon/ \Gamma^{diag}$. In
general $\Upsilon$ is not Hausdorff.  Note that the quotient group
$$
\Gamma^{quot}: = (\Gamma\times \Gamma)/\Gamma^{diag}
$$
naturally acts on the domain and on the range of the map $\Psi$,
with respect to which $\Psi$ is equivariant: the action of $(A_1,
A_2) \in \Gamma \times \Gamma$ on $T(\LL_2(M)) \times [0,1] \times
\R$ is given by
$$
(\xi, \eta, a) \mapsto (\xi, \eta, a - \omega(A_1 \# \overline
A_2)). \leqno (3.27)
$$
It follows, again by a standard calculation of the linearization
of $\widetilde\Psi$, $\Psi$ is transverse to
$$
o_{T\LL_2(M)} \times [0,1] \times \{0\} \subset T\LL_2(M) \times
[0,1] \times \R.
$$
In particular, we know that
$$
\WW : = \Psi^{-1}(o_{T\LL_2(M)} \times [0,1] \times   \{0\})
\leqno (3.28)
$$
is a smooth submanifold of  $\Upsilon$.

Now we consider the projection
$$
\Pi: \WW \subset  \Upsilon \to \PP(\HH_m(M);H_1,H_2).
$$
One can check that this is a Fredholm map of index 0. By the
Sard-Smale theorem, there is a dense subset of $\HH$'s such that
{\it as long as it is non-empty},
$$
\Pi^{-1}(\HH) = : W(\HH)
$$
is a smooth zero dimensional manifold.

Furthermore by the dimension counting argument, one proves that
the phenomena of multiple (with multiplicity more than 2) crossing
or of having a nonempty intersection $\SS ing(\HH) \cap \CC
ross^{nd}(\HH) \neq \emptyset$ generically have negative
codimension and so can be avoided over a residual subset of
$\PP(\HH_m(M);H_1,H_2)$, which we denote by
$\PP^{gc}(\HH_m(M);H_1,H_2)$.

Now consider the map
$$
p: W(\HH) \subset \widetilde \LL_2(M) \times [0,1] \to [0,1]
\times \R
$$
defined by
$$
([z_1,w_1], [z_2,w_2], \eta) \to \Big(\eta, \AA_{H(\eta)}([z_1,w_1])
\Big) \leqno (3.29)
$$
for which the diagram
$$
\begin{matrix} W(\HH) &  \longrightarrow  & [0,1] \times \R \\
\downarrow & \quad & \downarrow \\
[0,1] & \equiv & [0,1]
\end{matrix}
$$
commutes, and is equivariant under the action of $\Gamma^{diag}$
on $W(\HH)$, and the action (3.15). The image of $p$ is precisely
the set of crossings in the Cerf diagram $\Sigma(\HH)$. Now the
theorem follows by making a choice of any homotopy $\HH$ from
$\PP^{Cerf}(\HH_m(M);H_1,H_2)$.  \end{proof}

For the completeness's sake, we would like to consider the {\it
slopes} of the branches at a crossing. Note that nondegenerate
crossing implies that their slopes are different.

Suppose that $[z_\eta,w_\eta]$ locally parameterizes a branch at
$\eta = \eta_0$. The slope of the branch at $\eta$ is given by
$$
\frac{d}{d\eta}\left(\AA_{H(\eta)}([z_\eta,w_\eta])\right).
$$
Using the fact that $z_\eta$ is a Hamiltonian periodic orbit of
$H(\eta)$, we compute
$$
\begin{aligned}
\frac{d}{d\eta}\left(\AA_{H(\eta)}([z_\eta,w_\eta])\right) & =
d\AA_{H(\eta)}([z_\eta,w_\eta]) - \int_0^1
\frac{\partial H}{\partial \eta}(t,z_\eta(t))\, dt \\
& = - \int_0^1 \frac{\partial H}{\partial \eta}(t,z_\eta(t))\, dt.
\end{aligned}
\leqno (3.30)
$$
We remark that the slopes of the branches of periodic orbits are the
same as those of their liftings to $\widetilde\Omega_0(M)$.

\medskip{\noindent \bf Definition 3.8.} Suppose that $H_\alpha, \, H_\beta$
are two nondegenerate Hamiltonians. We call a homotopy $\HH$
between them a {\it Cerf homotopy} or a {\it Cerf family} if $\HH$
satisfies the properties described in Theorem 3.7 and in
Proposition 3.8. We call a point
$$
\eta \in [0,1] \setminus \SS ing(\HH)
$$
{\it Cerf regular}.
\medskip

By definition, we have
$$
\PP^{Cerf}(\HH_m(M);H_1,H_2)  := \PP^{reg}(\HH_m(M);H_1,H_2)  \cap
\PP^{gc}(\HH_m(M);H_1,H_2).
$$
From now on, we will always assume that $\HH$ is a Cerf homotopy,
unless otherwise stated.

For the later purpose, we introduce a concept of {\it
sub-homotopy}, which will play a crucial role in our proof of the
spectrality.

\medskip{\noindent \bf Definition 3.9.}
For each given pair $\eta_1 < \, \eta_2 \in [0,1]$, we consider
the homotopy $\HH_{\eta_1\eta_2}$ of $\HH$ between $H(\eta_1)$ and
$H(\eta_2)$ defined by the reparameterization
$$
\HH_{\eta_1\eta_2}: s \mapsto H^s = H((1-s)\eta_1+s \eta_2).
\leqno (3.31)
$$
We call any such a homotopy  $\HH_{\eta_1\eta_2}$ a {\it
sub-homotopy of $\HH$}. When $\eta_1 > \eta_2$, we define the
corresponding sub-homotopy by
$$
\HH_{\eta_1\eta_2} : = (\HH^{-1})_{\eta_2\eta_1}
$$
where $\HH^{-1}$ is the time reversal homotopy defined by
$$
\HH^{-1}: \eta \mapsto H(1-\eta).
$$
\medskip

\section*{\bf \S 4. Sub-homotopies and transversality}

\medskip

\subsection*{\it 4.1. Definition of the Floer boundary map, re-visited}

Suppose $H$ is a nondegenerate one-periodic Hamiltonian function
and $J$ be a one-periodic one-parameter family of compatible
almost complex structure. We denote by
$$
j_\omega= C^\infty([0,1],\JJ_\omega)
$$
the set of such families of almost complex structures. We first
recall the construction of the Floer boundary map and the
transversality conditions needed to define the Floer homology
$HF_*(H,J)$ of the pair. We will also add some novel elements in the
exposition of the construction, which are needed for our formulation
of the transversality problem of sub-homotopies.

The following definition is useful for the later discussion.
\medskip

{\noindent\bf Definition 4.1.} Let $z, \, z' \in \text{Per}(H)$. We
denote by $\pi_2(z,z')$ the set of homotopy classes of smooth maps
$$
u: [0,1] \times S^1:= T  \to M
$$
relative to the boundary
$$
u(0,t) = z(t), \quad u(1,t) = z'(t).
$$
We denote by $[u] \in \pi_2(z,z')$ its homotopy class and by $C$ a
general element in $\pi_2(z,z')$.
\medskip

We define by $\pi_2(z)$ to be the set of relative homotopy classes
of the maps
$$
w: D^2 \to M; \quad w|_{\part D^2} = z.
$$
We note that there is a natural action of $\pi_2(M)$ on $\pi_2(z)$
and $\pi_2(z,z')$ by the obvious operation of a `gluing a sphere'.
Furthermore there is a natural map of $C \in \pi_2(z,z')$
$$
(\cdot) \# C: \pi_2(z) \to \pi_2(z')
$$
induced by the gluing map
$$
w \mapsto w \# u.
$$
More specifically we will define the map $w \# u: D^2 \to M$ in
the polar coordinates $(r,\theta)$ of $D^2$ by the formula
$$
w \# u:(r,\theta) = \begin{cases} w(2r,\theta) & \quad  \text{for }\, 0
\leq r \leq \frac{1}{2} \\
w(2r-1,\theta)  & \quad \text{for } \, \frac{1}{2} \leq r \leq 1
\end{cases}
\leqno(4.1)
$$
once and for all. There is also the natural gluing map
$$
\pi_2(z_0,z_1) \times \pi_2(z_1,z_2) \to \pi_2(z_0,z_2)
$$
$$
(u_1, u_2) \mapsto u_1\# u_2.
$$
We also explicitly represent the map $u_1\# u_2: T \to M$ in the
standard way once and for all similarly to (4.1).
\medskip

{\noindent\bf Definition 4.2.} We define the {\it relative
Conley-Zehnder index} of $C \in \pi_2(z,z')$ by
$$
\mu_H(z,z';C) = \mu_H([z,w]) - \mu_H([z',w\# C]) \leqno(4.2)
$$
for a (and so any) representative $u:[0,1] \times S^1 \times M$ of
the class $C$. We will also write $\mu_H(C)$, when there is no
danger of confusion on the boundary condition.
\medskip

It is easy to see that the right hand side of (4.2) does not depend
on the choice of bounding disc $w$ of $z$, and so  the function
$$
\mu_H: \pi_2(z,z') \to \Z
$$
is well-defined.
\medskip

{\noindent\it Remark 4.3.} In fact, the function $\mu_H: \pi_2(z,z')
\to \Z$ can be defined without assuming $z_0, \, z_1$ being
contractible, as long as $z_0$ and $z_1$ lie in the same component
of $\Omega(M)$: For any given map $u: T \to M$, choose a {\it
marked} symplectic trivialization
$$
\Phi: u^*TM \to T \times \R^{2n}
$$
that satisfies
$$
\Phi\circ \Phi^{-1}|_{[0,1] \times \{1\}}=id.
$$
We know that $z_0(t) = \phi_H^t(p_0)$ and $z_1(t) = \phi_H^t(p_1)$
for $p_0, \, p_1 \in \text{Fix}(\phi_H^1)$. Then we have two maps
$$
\alpha_{\Phi,i}: [0,1] \to Sp(2n), \quad i = 0, \, 1
$$
such that
$$
\Phi\circ d\phi_H^t(p_i)\circ \Phi^{-1}(i,t,v) = (i,t,
\alpha_{\Phi,i}(t) v)
$$
for $v \in \R^{2n}$ and $t \in [0,1]$. By the nondegeneracy of
$H$, the maps $\alpha_{\Phi,i}$ define elements in $\SS P^*(1)$.
(See the Appendix A.1 for the definition of $\SS P^*(1)$.) Then we
define
$$
\mu_H(z,z';C): = \mu_{CZ}(\alpha_{\Phi,0}) -
\mu_{CZ}(\alpha_{\Phi,1}).
$$
It is easy to check that this definition does not depend on the
choice of marked symplectic trivializations.
\medskip

We now denote by
$$
\MM(H,J;z,z';C)
$$
the set finite energy solutions of
$$
\frac{\part u}{\part \tau} + J \Big (\frac{\part u}{\part t} -
X_H(u)\Big ) = 0 \leqno(4.3)
$$
with the asymptotic condition and the homotopy condition
$$
u(-\infty) = z, \quad u(\infty) = z'; \quad [u] = C. \leqno(4.4)
$$
Here we remark that although $u$ is a priori defined on $\R \times
S^1$, it can be compactified into a continuous map $\overline u:
[0,1] \times S^1 \to M$ with the corresponding boundary condition
due to the exponential decay property of solutions $u$ of (4.2),
recalling we assume $H$ is nondegenerate. We will call $\overline
u$ the {\it compactified map} of $u$. By some abuse of notation,
we will also denote by $[u]$ the class $[\overline u]\in
\pi_2(z,z')$ of the compactified map $\overline u$.

We now recall that the Floer boundary map
$$
\part_{(H,J)}; CF_{k+1}(H) \to CF_k(H)
$$
is defined under the following conditions. (See \cite{Fl3},
\cite{HS}.)
\medskip

{\noindent\bf Definition 4.4 (The boundary map).} Let $H$ be nondegenerate.
Suppose that $J$ satisfies the following conditions:

\begin{itemize}

\item[(1)] For any pair $(z_0,z_1) \subset \text{Per}(H)$ satisfying
$$
\mu_H(z_0,z_1;C) = \mu_H([z_0,w_0]) - \mu_H([z_1,w_0\# C])=0,
$$
$\MM\left(H,J;z_0,z_1;C\right) = \emptyset$ unless $z_0 = z_1$ and
$C = 0$. When $z_0=z_1$ and $C=0$, the only solutions are the
stationary solution, i.e., $u(\tau) \equiv z_0 = z_1$ for all
$\tau \in \R$.

\item[(2)] For any pair $(z_0,z_1) \subset \text{Per}(H)$ and a
homotopy class $C \in \pi_2(z_0,z_1)$ satisfying
$$
\mu_H(z_0,z_1;C) =1,
$$
$\MM(H,J;z_0,z_1;C)/\R$ is transverse and compact and so a finite
set. We denote
$$
n(H,J;z_0,z_1;C) = \#(\MM(H,J;z_0,z_1;C)/\R)
$$
the algebraic count of the elements of the space
$\MM(H,J;z_0,z_1;C)/\R$. We set $n(H,J;z_0,z_1:C) = 0$ otherwise.

\item[(3)] For any pair $(z_0,z_2) \subset \text{Per}(H)$ and $C \in
\pi_2(z_0,z_2)$ satisfying
$$
\mu_H(z_0,z_2;C) = 2,
$$
$\MM(H,J;z_0,z_2;C)/\R$ can be compactified into a smooth
one-manifold with boundary comprising the collection of the broken
trajectories
$$
[u_1] \#_\infty [u_2]
$$
where $u_1 \in \MM(H,J;z_0,y:C_1)$ and $u_2 \in
\MM(H,J;y,z_2:C^2)$ for all possible $y \in \text{Per}(H)$ and
$C_1 \in \pi_2(z_0,y), \, C_2 \in \pi_2(y,z_2)$ satisfying
$$
C_1 \# C_2 = C; \quad [u_1] \in \MM(H,J;z_0,y;C_1)/\R, \quad [u_2]
\in \MM(H,J;y,z_2;C_2)/\R
$$
and
$$
\mu_H(z_0,y;C_1) = \mu_H(y,z_2;C_2) = 1.
$$
Here we denote by $[u]$ the equivalence class  represented by $u$.
\end{itemize}
We call any such $J$ {\it $H$-regular} and call any such pair
$(H,J)$ {\it Floer regular}.
\medskip

The upshot is that for a Floer regular pair $(H,J)$ the Floer
boundary map
$$
\part=\part_{(H,J)}: CF_*(H) \to CF_*(H)
$$
is defined and satisfies $\part\part = 0$ and so the Floer homology
$HF(H,J): = \ker \part/\text{im }\part$ is defined. For any given
nondegenerate $H$, the set of $H$-regular $J$'s is dense in
$j_\omega$ {\it under the assumption of semi-positivity}. (See
\cite{Fl3}, \cite{HS} for the proof.) We denote by
$$
j_\omega^{reg}(H) \subset j_\omega = C^\infty([0,1],\JJ_\omega)
$$
the set of $H$-regular $J$'s.

\medskip
\subsection*{\it 4.2. Definition of the Floer chain map, revisited}

\smallskip

 Now we study the Floer regularity of the
triple $(\HH,j;\rho)$. We need to study an intermediate problem.
Consider the pair $(\HH_\R,j_\R)$ of maps
$$
\begin{aligned} \HH_\R & : \R \to C^\infty_m([0,1]\times M):= \HH_m(M) \\
& j_\R : \R \to j_\omega
\end{aligned}
$$
that are {\it asymptotically constant} i.e., such that there
exists $R> 0$ such that
$$
H(\tau)\equiv H(\pm\infty), \quad j(\tau) \equiv j(\pm\infty)
$$
for $\tau$ with $|\tau| \geq R$. To each such a pair, we associate
non-autonomous analog to (4.3) still with the condition (4.4).
\medskip

{\noindent\bf Definition 4.5 (The chain map).}
We say that {\it $(\HH_\R, j_\R)$ is Floer regular} if the
following holds: \begin{itemize}

\item[(1)] For any pair $z_0 \in \text{Per}(H_0)$ and $z_1 \subset
\text{Per}(H_1)$ satisfying
$$
\mu_{\HH_\R}(z_0,z_1;C) =0,
$$
$\MM(\HH_{\R},j_{\R};z_0,z_1;C)$ is transverse and compact, and so
a finite set. We denote
$$
n(\HH_\R,j_\R;z_0,z_1;C) = \#(\MM(\HH_\R,j_\R;z_0,z_1;C))
$$
the algebraic count of the elements in
$\MM(\HH_\R,j_\R;z_0,z_1;C)$. We set $n(\HH_\R,j_\R;z_0,z_1:C) =
0$ otherwise.

\item[(2)] For any pair $z_0 \in \text{Per}(H_0)$ and $z_1 \in
\text{Per}(H_1)$ satisfying
$$
\mu_{\HH_\R}(z_0,z_2;C) = 1,
$$
$\MM(H,J;z_0,z_2;C)$ is transverse and can be compactified into a
smooth one-manifold with boundary comprising the collection of the
broken trajectories
$$
u_1 \#_\infty u_2
$$
where
$$
\begin{aligned} (u_1,u_2) \in  & \MM(\HH_\R,j_\R; z_0,y:C_1) \times
\MM(H(\infty),J(\infty);y,z_2:C^2); \\
& \mu_{\HH_\R}(z_0,y;C_1) = 0, \, \mu_H(y,z_2;C_2) = 1
\end{aligned}
$$
or
$$
\begin{aligned} (u_1,u_2) \in & \MM(H(-\infty), J(-\infty);z_0,y:C^1)
\times
\MM(\HH_\R,j_\R;y,z_2:C_1) ; \\
& \mu_{\HH_\R}(z_0,y;C_1) = 1, \, \mu_H(y,z_2;C_2) = 0
\end{aligned}
$$
and $C_1 \# C_2 = C$ for all possible such $y \in \text{Per}(H)$
and $C_1 \in \pi_2(z_0,y)$, $C_2 \in \pi_2(y,z_2)$.
\end{itemize}
We say that $(\HH_\R,j_\R)$ are {\it Floer regular} if it
satisfies these conditions.
\medskip

Again for any fixed $\HH_\R$, the residual property of $j_\R$'s for
which $(\HH_\R,j_\R)$ are regular is well-known {\it for the
semi-positive case}. (See \cite{Fl3}, \cite{HS}.)

Now suppose that $\HH$ is a homotopy connecting two nondegenerate
Hamiltonians $H_0$ and $H_1$. We denote
$$
\PP(j_\omega): = C^\infty([0,1], j_\omega)
$$
the set of smooth one-parameter family $j= \{J(s)\}_{0 \leq s\leq
1}$ with $J(s) \in j_\omega$. We define a function $\rho :\R \to
[0,1]$ of the type
$$
\begin{aligned}
\rho(\tau) & = \begin{cases} 0 \, \quad \text {for $\tau \leq -R$}\\
                    1 \, \quad \text {for $\tau \geq R$}
                    \end{cases}
\end{aligned}
$$
for some $R > 0$. We call $\rho$ a {\it monotone} cut-off function
if it satisfies $\rho'(\tau) \geq 0$ for all $\tau$'s in addition.

Each such pair $(\HH,j)$ and a cut-off function $\rho$ define a
pair
$$
\HH_\R = \HH^\rho, \quad j_\R = j^\rho
$$
where $\HH^\rho$  is the reparameterized homotopy $\HH^\rho =
\{H^\rho\}_{\tau \in \R}$ defined by
$$
\tau \mapsto H^{\rho}(\tau,t,x) = H(\rho(\tau),t,x).
$$
We call $\HH^\rho$ the {\it $\rho$-elongation} of $\HH$ or the
{\it $\rho$-elongated homotopy} of $\HH$. The same definition
applies to  $j$. Therefore such a triple $(\HH,j;\rho)$ associates
the {\it non-autonomous} equation
$$
\begin{cases} \frac{\part u}{\part \tau} + J^{\rho(\tau)}\Big(\frac{\part
u}{\part t}
- X_{H^{\rho(\tau)}}(u)\Big) = 0\\
\lim_{\tau \to -\infty}u(\tau) = z^-,  \lim_{\tau \to
\infty}u(\tau) = z^+
\end{cases}
\leqno(4.5)
$$
with the boundary condition
$$
u(-\infty) = z_0, \quad u(\infty) = z_1. \leqno(4.6)
$$
We denote by
$$
\MM((\HH,j;\rho);z_0,z_1;C)
$$
the set of finite energy solutions of (4.5)-(4.6) satisfying the
topological condition $[u] = C$ in $\pi_2(z_0,z_1)$. We say that
$(\HH,j;\rho)$ is {\it Floer regular} if the $\rho$-elongation
$(\HH^\rho,j^\rho)$ is Floer regular in the sense of Definition
4.5.

\medskip
\subsection*{\it 4.3. Transversality of sub-homotopies}

\smallskip

With these preparations, we now launch our main study of the
transversality question on the sub-homotopies of a given pair
$(\HH,j)$.

Let $\HH$ be a Cerf homotopy and $j \in \PP(j_\omega)$ be given.
As for the homotopy $\HH$ of Hamiltonian functions, we define the
{\it sub-homotopy} $j_{\eta\eta'}$ by
$$
s\in [0,1] \mapsto J((1-s)\eta + s \eta') \leqno(4.7)
$$
for each given $0 \leq \eta \leq \eta' \leq 1$. When $\eta >
\eta'$, we define
$$
j_{\eta\eta'}: = (j^{-1})_{\eta'\eta}.
$$
The main purpose of the present section is to study the question
whether one can define the Floer chain map
$$
h_{(\eta\eta';\rho)}:CF_*(H(\eta)) \to CF_*(H(\eta'))
$$
for {\it sufficiently many} points of $\eta' \in [0,1]$, when
there are given a generic homotopy $(\HH,j)$ and a Cerf regular
point $\eta \in [0,1]$.

The following theorem is the main theorem in that regard.
\medskip

{\noindent\bf Theorem 4.6.} {\it Let $\HH \in
\PP^{Cerf}(\HH_m(M);H_\alpha,H_\beta)$ be a Cerf homotopy
connecting two nondegenerate Hamiltonians $H_\alpha, \, H_\beta$.
Let $J_i, \, i = \alpha, \, \beta$  $(H_i,J_i)$ be Floer regular
and fix a cut-off function $\rho: \R \to [0,1]$. Then we have the
following:

\begin{itemize} \item[(1)] There exists a dense subset
$$
\PP^{tran}(j_\omega;\HH) \subset \PP(j_\omega)
$$
such that for any element $j$ from $\PP^{tran}(j_\omega;\HH)$
there exists a residual subset of $[0,1]$ containing $\{0,1\}$,
denoted by
$$
I(\HH,j) \subset [0,1],
$$
at each point $\eta$ of which the pair $(H(\eta), J(\eta))$ is
Floer regular.

\item[(2)] For any $\eta \in I(\HH,j)$, there exists a residual subset
$$
\PP^{sub}(j_\omega, \HH;\eta) \subset \PP^{tran}(j_\omega;\HH)
$$
such that for any $j \in \PP^{sub}(j_\omega;\HH;\eta)$ there
exists a residual subset
$$
I(\HH,j;\eta) \subset I(\HH,j) \subset [0,1]
$$
such that for any $\eta' \in I(\HH,j;\eta)$ the triple
$$
(\HH_{\eta\eta'},j_{\eta\eta'};\rho)
$$
is Floer regular and hence the Floer chain map
$$
h_{(\HH_{\eta\eta'},j_{\eta\eta'};\rho)}:=
h_{(\HH_{\eta\eta'}^\rho,j_{\eta\eta'}^\rho)}: CF_*(H_0) \to
CF_*(H_1)
$$
is defined and satisfies
$$
h_{(\HH_{\eta\eta'},j_{\eta\eta'};\rho)} \circ
\part_{(H_\alpha,J_\alpha)} = \part_{(H_\beta,J_\beta)} \circ
h_{(\HH_{\eta\eta'},j_{\eta\eta'};\rho)}.
$$
\end{itemize}
}

The rest of the section will be occupied by the proof of this
theorem. We fix a Cerf homotopy $\HH$ and a cut-off function
$\rho$.

We consider the case (1) first. For each fixed $z_0 \in
\text{Per}(H_0)$ and $z_1 \in \text{Per}(H_1)$, and a class $C \in
\pi_2(z_0,z_1)$, we consider the space
$$
C^\infty(z_0,z_1; C): = \{ u:\R \times S^1 \to M \mid u \,
\text{satisfies (4.6) and } \, E_{(H,J)}(u) < \infty \}
$$
and its $W^{1,p}$-completion with respect to a suitably weighted
Sobolev norm on $C^\infty(z_0,z_1; C)$. (See \cite{Fl1}.) We denote
the corresponding weighted $W^{1,p}$-space as
$$
W^{1,p}(z_0,z_1; C).
$$
Then we consider the assignment
$$
(u,J) \mapsto \frac{\part u}{\part \tau} +
J^{\rho(\tau)}\Big(\frac{\part u}{\part t} -
X_{H^{\rho(\tau)}}(u)\Big) \leqno(4.8)
$$
as a section of a vector bundle $\EE(z_0,z_1;C)$ over
$W^{1,p}(z_0,z_1;C) \times j_\omega$ whose fiber is given by
$$
L^p(z_0,z_1;C): = L^p\Big(\Lambda^{0,1}(u^*TM)\Big).
$$
We denote this section by
$$
\overline\part_{H}: W^{1,p}(z_0,z_1;C) \times j_\omega \to
\EE(z_0,z_1;C), \leqno(4.9)
$$
and denote by $\underline 0$ the zero section of the vector bundle
$\EE(z_0,z_1;C)$. Then it is well-known \cite{Fl1}, \cite{Fl3} that
the covariant linearization of $\overline\part_{H}$ is surjective
and so the zero set
$$
\MM(H;z_0,z_1;C) : = (\delbar_H)^{-1}(\underline 0) \subset
W^{1,p}(z_0,z_1;C) \times j_\omega \leqno(4.10)
$$
is a smooth submanifold whose image is indeed contained in $
C^\infty(z_0,z_1; C) \times j_\omega$. Furthermore the projection
$$
\Pi_{(H;z_0,z_1;C)}: \MM(H;z_0,z_1;C) \to j_\omega \leqno(4.11)
$$
is a Fredholm map of index $\mu_H(C)$. By the Sard-Smale theorem,
the set $j_\omega^{reg}(H) \subset j_\omega$ is a residual subset
and so dense.

To prove (1), we consider $j: [0,1] \to j_\omega$ with $j(0) =
J_0$ and $j(1) = J_1$. Then applying the Sard-Smale theorem, it is
enough to consider the set of smooth paths $j: [0,1] \to j_\omega$
that are transverse to $\Pi_{(H;z_0,z_1;C)}$ for all triple
$(z_0,z_1;C)$. Since $H$ is assumed to be nondegenerate, there are
only finitely many pairs $(z_0,z_1)$ and so only countably many
possible triples $(z_0,z_1;C)$. We denote by
$$
\PP^{tran}(j_\omega;\HH)
$$
the set of such $j= \{J(\eta)\}_{0 \leq \eta \leq 1}$. Then
$$
\MM(\HH,j;z_0,z_1;C):= \bigcup_{\eta \in
[0,1]}\MM(H(\eta),J(\eta);z_0,z_1;C) \subset [0,1] \times
W^{1,p}(z_0,z_1;C)
$$
is a smooth manifold of dimension $\mu_H(C) + 1$. There is a
canonical projection
$$
\pi_j:\MM(\HH,j;z_0,z_1;C) \to [0,1].
$$
By the classical Sard theorem, the set of regular values of
$\pi_j$, denoted by $I(\HH,j) \subset [0,1]$, is residual and so
dense. This finishes the proof of (1).

For the proof of (2), we first note that there is a natural map
$$
\SS ub_\eta:\PP(j_\omega) \times [0,1] \to \PP(j_\omega)
$$
defined by taking the sub-homotopy
$$
\SS ub_\eta(j, \eta') = j_{\eta\eta'}; \quad j_{\eta\eta'}(s): =
j((1-s)\eta + s\eta'). \leqno(4.12)
$$
We need to study Floer regularity of the triple
$$
(\HH_{\eta\eta'},j_{\eta\eta'};\rho).
$$
Fix an element $j \in \PP^{tran}(j_\omega;\HH)$. Note that by
Definition 4.5, the triple $(\HH_{\eta\eta},j_{\eta\eta};\rho)$ is
Floer regular if $j \in \PP^{trans}(j_\omega;\HH)$ and if the
corresponding moduli spaces satisfying
$$
\MM((\HH_{\eta\eta},j_{\eta\eta};\rho);z,z';C) = \emptyset
$$
{\it except when $z=z'$ and $C = \underline 0$.} In particular,
the are only a {\it finite} number of moduli spaces to consider,
each of which consists of a single stationary element.

Since the structure of the set $\text{Per}(H(\eta'))$ will vary as
$\eta'$ changes,  we need to find a way of encoding the asymptotic
conditions of (4.5) corresponding to
$(\HH_{\eta\eta'},j_{\eta\eta'})$ as $\eta'$ varies. For example
it may experience bifurcations as $\eta'$ changes. {\it Here
enters again the property of Cerf homotopy $\HH$ in a crucial
way}.

We partition $[0,1]$ into
$$
0= s_0 < s_1 < s_2 < \cdots < s_{k_1} = 1
$$
where $\SS ing(\HH) = \{ s_1, \ldots, s_k\}$, and denote $I_n =
(s_{n-1}, s_n)$. Since a bifurcation of periodic orbits occurs only at the points $s_n$,
we can smoothly parameterize $\text{Per}(H(\eta))$ on each $I_n$.
For each given pair $(n,n')$ and a pair of points $(\eta,\eta') \in
I_{n} \times I_{n'}$, we denote by $(z,z')$ the pair of branches
$z=z(\eta)$ for $\eta \in I_n$ and $z'=z'(\eta')$ for $\eta \in
I_{n'}$. For each given such a pair, we denote by
$$
\MM_{n}^{n'}((\HH_{\eta\eta},j_{\eta\eta};\rho);z,z';C) \leqno(4.13)
$$
the corresponding moduli space of solutions of (4.5). Therefore we
need to consider only a finite number of possibilities. To study
the transversality of (4.13), we need to set up the Fredholm
theory for each given $(n,n';z,z')$ and a given $C \in
\pi_2(z,z')$.

Let $\eta \in [0,1]$ be given and $\eta \in I_n$, and fix another
interval $I_{n'}$. We consider the assignment
$$
\delbar_{(\HH;\eta)}: (u,j,\eta') \mapsto
\delbar_{(\HH_{\eta\eta'},j_{\eta\eta'};\rho)}(u) \leqno(4.14)
$$
for $\eta' \in I_{n'}$. The linearization is surjective as before
and so the universal moduli space
$$
(\delbar_{(\HH;\eta)})^{-1}(\underline
0)=:\MM_{n}^{n'}(\HH;\rho;z,z';C)
$$
is a smooth submanifold of
$$
C^\infty(\R \times S^1;z,z';C) \times \PP(j_\omega) \times I_{n'}.
$$
The projection map
$$
\Pi_\HH:\MM_{n}^{n'}(\HH;\rho;z,z';C) \to \PP(j_\omega) \leqno (4.14)
$$
is a Fredholm map of index $\mu_{\HH_{\eta\eta'}}(C) + 1$ in
general. We need to study the cases $\mu_{\HH{\eta\eta'}}(C) = 0$
and $\mu_{\HH_{\eta\eta'}}(C)=1$. We will provide the details only
for the case $\mu_{\HH_{\eta\eta'}}(C) = 0$ and leave the other
case to the readers.

Again by the Sard-Smale theorem, the regular values of $\Pi_\HH$
is residual for each choice of $(n',z';C)$ with $C \in
\pi_2(z,z')$. We define by
$$
\PP^{sub}(j_\omega;\HH;\eta)
$$
the intersection of the sets of regular values over all possible
$(n',z';C)$ and $\PP^{tran}(j_\omega;\HH)$. Certainly we have
$$
\PP^{sub}(j_\omega;\HH;\eta)\subset \PP^{tran}(j_\omega;\HH)
$$
and is a residual subset of $\PP(j_\omega)$, since there are only
countably many choices of $(n',z';C)$.

By the definition of $\PP^{sub}(j_\omega;\HH;\eta)$, the preimage
$$
\Pi_\HH^{-1}(j) = \bigcup_{\eta'\in I_{n'}}
\MM_{n}^{n'}((\HH_{\eta\eta'},j_{\eta\eta'};\rho);z,z';C) \subset
C^\infty(z,z';C) \times I_{n'} \leqno (4.15)
$$
is a smooth manifold of dimension 1 for each given $n'$, which
forms a fibration over $I_{n'}$. We denote by
$$
\pi_n^{n'}: \Pi_\HH^{-1}(j) \to I_{n'} \leqno (4.16)
$$
the natural projection map. Denote by $I_{n'}^{reg}[n]$ the set of
regular values of $\pi_n^{n'}$. Then we have only to define
$I(\HH,j;\eta)$ to the union
$$
I(\HH,j;\eta): = \cup_{n'=1}^{k_1}I_{n'}^{reg}[n].
$$
This proves the proof of (2) and so the theorem.

Similar discussion applies to the case $\mu_{\HH_{\eta\eta'}}(C) =
1$ and omitted.

\section*{\bf \S 5. The composition law of Floer's chain maps, revisited}

In this section, we will re-examine the well-known composition raw
$$
h_{\alpha\gamma} = h_{\beta\gamma} \circ h_{\alpha\beta} \leqno (5.1)
$$
of the Floer's canonical isomorphism \cite{Fl3}
$$
h_{\alpha\beta}: HF_*(H_\alpha) \to HF_*(H_\beta). \leqno (5.2)
$$

We first carefully review this construction {\it in the chain
level}. Although the isomorphism (5.2) {\it in homology} depends
only on the end Hamiltonians $H_\alpha$ and $H_\beta$, the
corresponding chain map depends on the homotopy $\HH =
\{H(\eta)\}_{0 \leq \eta \leq 1}$ between $H_\alpha$ and
$H_\beta$, and also on the homotopy $j = \{J(\eta)\}_{0\leq \eta
\leq 1}$. Let us fix nondegenerate Hamiltonians $H_\alpha, \,
H_\beta$ and a homotopy $\HH$ between them. We also fix a homotopy
$j = \{J(\eta)\}_{0 \leq \eta \leq 1}$ of compatible almost
complex structures and a cut-off function $\rho:\R \to [0,1]$. The
Floer chain map is defined by considering the non-autonomous
equation
$$
\begin{cases} \frac{\part u}{\part \tau} + J^{\rho(\tau)}\Big(\frac{\part
u}{\part t}
- X_{H^{\rho(\tau)}}(u)\Big) = 0\\
\lim_{\tau \to -\infty}u(\tau) = z^-,  \lim_{\tau \to
\infty}u(\tau) = z^+
\end{cases}
\leqno (5.3)
$$
with the condition that $w^+$ and $w^-\# u$ are homotopic to each
other relative to the boundary. We denote this condition by
$$
w^+ \sim w^- \# u. \leqno (5.4)
$$
One consequence of (5.4) is
$$
[z^+,w^+] = [z^+, w^-\# u] \quad \text{in } \quad \Gamma
$$
but the latter is a much weaker condition than the former. The
asymptotic condition with (5.4) is equivalent to (4.6)
corresponding to the class $C=[u]$.

Considering any such solution $u$ as a path in the covering space
$\widetilde \Omega_0(M)$, we will also write the asymptotic
condition as
$$
\lim_{\tau \to -\infty}u(\tau) = [z^-,w^-], \quad \lim_{\tau \to
\infty}u(\tau) = [z^+,w^+]. \leqno (5.5)
$$
$h_{(\HH,j;\rho)}$ has degree 0 and satisfies
$$
\part_{(J^1,H^1)} \circ h_{(\HH,j;\rho)} = h_{(\HH,j;\rho)} \circ
\part_{(J^0,H^0)}.
\leqno (5.6)
$$
In general two such maps $h_{(\HH,j;\rho_1)}$ and
$h_{(\HH,j;\rho_2)}$ are chain homotopic to each other in the
sense of (5.7) below. In the end of this section, we will
carefully study the dependance of this chain map on the cut-off
functions.
\medskip

{\noindent\it Remark 5.1.} One may directly consider the
$\R$-family of Hamiltonians that are asymptotically constant as
some literature do, which may remove the additional choice $\rho$
in the construction and be more natural in the point of view of
Hamiltonian fibrations. However in this paper we prefer to use our
approach because the usual homotopy is defined over $[0,1]$.
Composing the homotopy with the cut-off function $\rho: \R \to
[0,1]$ automatically makes the corresponding elongated homotopy
asymptotically constant, which is a necessary requirement for the
many parts of analysis of the perturbed Cauchy-Riemann equation:
e.g., the gluing theorem and the index theorem of the solutions of
the perturbed Cauchy-Riemann equation. In fact, using the cut-off
function, we can consider the elongation of the sub-homotopy
$\HH_{\eta\eta'}$ of a given homotopy $(\HH,j)$ {\it in one step}
by considering the functions $\rho$ of the type satisfying
$$
\begin{aligned}
\rho(\tau) & = \begin{cases} \eta \quad \text{for } \, \tau \leq -R \\
\eta' \quad \text{for } \, \tau \geq R
\end{cases} \\
\eta & \leq \rho(\tau) \leq \eta'.
\end{aligned}
$$
We will elaborate this point elsewhere.

When we are given a homotopy $(\overline j, \overline \HH)$ of
homotopies with $\overline j = \{j_\kappa\}$, $\overline\HH =
\{\HH_\kappa\}$, we also define the elongations
$\HH^{\overline\rho}$ of $\HH_\kappa$ by a homotopy of cut-off
functions $\overline \rho=\{\rho_\kappa\}$: we have
$$
\HH^{\overline\rho} = \{ \HH_\kappa^{\rho_\kappa} \}_{0 \leq
\kappa \leq 1}.
$$
Consideration of the parameterized version of (5.3) for $ 0 \leq
\kappa \leq 1$ defines the chain homotopy map
$$
H_{\overline\HH} :CF_*(H_\alpha) \to CF_*(H_\beta)
$$
which has degree $+1$ and satisfies
$$
h_{(j_1, \HH_1;\rho_1)} - h_{(j_0,\HH_0:\rho_0)} =
\part_{(J^1,H^1)} \circ H_{\overline\HH} + H_{\overline\HH} \circ
\part_{(J^0,H^0)}. \leqno (5.7)
$$
Again the map $\HH_{\overline\HH}$ depends on the choice of a
homotopy $\overline j$ and $\overline\rho = \{\rho_\kappa\}_{0
\leq \kappa \leq 1}$ connecting the two functions $\rho_0, \,
\rho_1$. Therefore we will denote
$$
H_{\overline \HH} = H_{(\overline \HH,\overline j; \overline
\rho)}.
$$
(5.7) in particular proves that two such chain maps (5.2) for
different homotopies $(j_0,\HH_0;\rho_0)$ and $(j_1, \HH_1;\rho_1)$
connecting the same end points are chain homotopic \cite{Fl3} and so
proves that the isomorphism (5.2) in homology is independent of the
homotopies $(\overline\HH,\overline j)$ or of $\overline \rho$.

Now we re-examine the equation (5.3).  The key analytic fact in
the proof of (5.6) or (5.7) is an a priori upper bound of the
energy
$$
E_{(\HH,j;\rho)}(u) := \frac{1}{2}\int_{-\infty}^\infty\int_0^1
\Big(\Big| \frac{\part u}{\part \tau} \Big|_{J^{\rho(\tau)}}^2+
\Big|\frac{\part u}{\part t} -
X_{H^{\rho(\tau)}}(u)\Big|_{J^{\rho(\tau)}}^2\Big)\, dt\,d\tau
$$
for the solutions $u$ of (5.3) with (5.4). In this respect, we
recall the following standard identity.
\medskip

{\noindent\bf Lemma 5.2.} {\it Let $(\HH,j)$ be any pair as above,
not necessarily generic. Suppose that $u$ satisfies $(5.3)$, has
finite energy and satisfies
$$
\lim_{j \to \infty} u(\tau^-_j) = [z^-,w^-], \quad \lim_{j \to
\infty}u(\tau^+_j) = [z^+,w^+]
$$
for some sequences $\tau^\pm_j$ with $\tau^-_j \to -\infty$ and
$\tau_j^+ \to \infty$. Then we have
$$
\begin{aligned} \AA_F([z^+,w^+]) & - \AA_H([z^-,w^-])\\
& = - \int \Big|\dudtau \Big|_{J^{\rho_1(\tau)}}^2 -
\int_{-\infty}^\infty \rho'(\tau)\int_0^1 \Big( \frac{\part
H^s}{\part s}\Big|_{s = \rho(\tau)} (t, u(\tau,t))\Big) \,
dt\,d\tau
\end{aligned}
\leqno (5.8)
$$
}
\medskip

{\noindent\bf Corollary 5.3.} {\it Let $(\HH,j)$ and $u$ be as in Lemma 5.2.
\begin{itemize} \item[(1)] Suppose that $\rho$ is monotone in addition. Then we
have
$$
\begin{aligned} \AA_F([z^+,w^+]) - \AA_H([z^-,w^-]) & \leq - \int
\Big|\dudtau \Big|_{J^{\rho_1(\tau)}}^2 + \int_0^1
-\min_{x, \, s} \Big( \frac{\part H^s_t}{\part s}\Big)\, dt
\end{aligned}\leqno (5.9 )
$$
$$
\begin{aligned}
& \leq \int_0^1 -\min_{x, \, s} \Big( \frac{\part H_t^s}{\part s}
\Big)\, dt.
\end{aligned}\leqno (5.10)
$$
And $(5.10)$ can be rewritten as the upper bound for the energy
$$
\begin{aligned} \int \Big|\dudtau \Big|_{J^{\rho_1(\tau)}}^2 & \leq
\AA_H([z^+,w^+]) - \AA_F([z^-,w^-]) \\
& \qquad + \int_0^1 -\min_{x, \, s} \Big( \frac{\part H^s_t}{\part
s}\Big)\, dt.
\end{aligned}
\leqno (5.11)
$$
\item[(2)] For a general $\rho$, we instead have
$$
\begin{aligned} \AA_F([z^+,w^+]) - \AA_H([z^-,w^-]) & \leq - \int
\Big|\dudtau \Big|_{J^{\rho_1(\tau)}}^2 + \int_0^1
\max_{x, \, s} \Big|\frac{\part H^s_t}{\part s}\Big|\, dt
\end{aligned} \leqno (5.12)
$$
$$
\begin{aligned}
& \leq \int_0^1 \max_{x, \, s} \Big|\frac{\part H_t^s}{\part s}
\Big|\, dt.
\end{aligned}\leqno (5.13)
$$
And $(5.13)$ can be rewritten as the upper bound for the energy
$$
\begin{aligned} \int \Big|\dudtau \Big|_{J^{\rho_1(\tau)}}^2 & \leq
\AA_H([z^+,w^+]) - \AA_F([z^-,w^-]) \\
& \qquad + \int_0^1 \max_{x, \, s} \Big|\frac{\part H^s_t}{\part
s}\Big|\, dt.
\end{aligned}
\leqno (5.14)
$$
\end{itemize}
}

\begin{proof} The proof is an immediate consequence of (5.8) and
omitted.
\end{proof}

Here we would like to emphasize that the upper estimates
(5.10)-(5.11) or (5.13)-(5.14) do not depend on the choice of $j$
or of $\rho$, but depend only on the homotopy $\HH$ itself.

Motivated by the upper estimate (5.10), we introduce the following
definition
\medskip

{\noindent\bf Definition 5.4.} Let $\HH = \{H(s)\}_{0 \leq s \leq 1}$
be a homotopy of Hamiltonians. We define the {\it negative part of
the variation} and the {\it positive part of the variation} of
$\HH$ by
$$
\begin{aligned} E^-(\HH) & := \int_0^1 -\min_{x, \, s} \Big( \frac{\part
H_t^s}{\part s} \Big)\, dt. \\
E^+(\HH) & := \int_0^1 \max_{x, \, s} \Big( \frac{\part
H_t^s}{\part s} \Big)\, dt.
\end{aligned}
$$
And we define the {\it total variation} $E(\HH)$ of $\HH$ by
$$
E(\HH) = E^-(\HH) + E^+(\HH).
$$
\medskip

If we denote by $\HH^{-1}$ the time reversal of $\HH$, i.e., the
homotopy given by
$$
\HH^{-1}: s \in [0,1] \mapsto H^{1-s}
$$
then we have the identity
$$
E^\pm(\HH^{-1}) = E^\mp(\HH) \quad \text{and } E(\HH^{-1}) =
E(\HH).
$$
With these definitions, applied to a pair $(\HH,j)$ such that
their ends $H(0)$ and $H(1)$ are nondegenerate, the a priori
energy estimate (5.11) can be written as
$$
\int \Big|\dudtau \Big|_{J^{\rho(\tau)}}^2 \leq  -
\AA_F([z^+,w^+]) +\AA_H([z^-,w^-]) + E^-(\HH) \leqno (5.15)
$$
for a monotone $\rho$, and (5.14) as
$$
\int \Big|\dudtau \Big|_{J^{\rho(\tau)}}^2 \leq  -
\AA_F([z^+,w^+]) + \AA_H([z^-,w^-]) + E(\HH) \leqno (5.16)
$$
for a general $\rho$.

Note that when $\HH$ is the linear homotopy
$$
\HH^{lin}: s \mapsto (1-s)H_1 + sH_2
$$
between $H_1$ and $H_2$, $E^\pm(\HH)$ and $E(\HH)$ just become
$E^\pm(H_2 -H_1)$, and $\|H_2-H_1\|$, respectively: We recall the
definitions
$$
\begin{aligned} E^-(H) & = \int_0^1 - \min_xH_t\, dt, \quad E^+(H) =
\int_0^1 \max_xH_t\, dt \\
\|H\|& = E^+(H) + E^-(H) = \int_0^1(\max_x H_t - \min_x H_t)\,dt.
\end{aligned}
$$
Therefore, taking the infimum of $E(\HH)$ over all $\HH$ with
fixed end points $H(0) = H_0$ and $H(1) = H_1$, we have the
inequality
$$
\inf_{\HH} \Big\{E(\HH) \mid H(0) = H_0, \, H(1) = H_1 \Big\} \leq
\|H_1 - H_0\|
$$
which is a strict inequality in general. It seems to be an
interesting problem to investigate the geometric meaning of the
quantity in the left hand side. This will be a subject of the
future study.

Next, we consider the triple
$$
(H_\alpha, \, H_\beta, \, H_\gamma)
$$
of Hamiltonians and homotopies $\HH_1, \, \HH_2$ connecting from
$H_\alpha$ to $H_\beta$ and $H_\beta$ to $H_\gamma$ respectively.
We define their concatenation $\HH_1 \# \HH_2 = \{H_3(s)\}_{1 \leq
s \leq 1}$ by
$$
H_3(s) = \begin{cases} H_1(2s) &\quad 0 \leq s \leq \frac{1}{2} \\
H_2(2s-1) & \quad \frac{1}{2} \leq s \leq 1.
\end{cases}
$$
From the definition of $E^\pm$ and $E$ for the homotopy $\HH$
above, we immediately have the following lemma
\medskip

{\noindent\bf Lemma 5.5.} {\it All $E^\pm$ and $E$ are additive under the
concatenation of homotopies. In other words, for any triple
$(H_\alpha,H_\beta,H_\gamma)$ and homotopies $\HH_1, \, \HH_2$ as
above, we have
$$
E^\pm(\HH_1 \# \HH_2) = E^\pm(\HH_1) + E^\pm(\HH_2).
$$
The same additivity holds for $E$.
}

Next we note that due to the choice of the cut-off function
$\rho$, the continuity equation (5.3) is {\it autonomous}  for the
region $|\tau|
> R$ i.e., is invariant under the translation by $\tau$. When we
are given a triple $(H_\alpha, \, H_\beta, \, H_\gamma)$, this
fact enables us to glue solutions of two such equations
corresponding to the pairs $(H_\alpha,H_\beta)$ and
$(H_\beta,H_\gamma)$ respectively.

Now a more precise explanation is in order. For a given pair of
cut-off functions
$$
\rho = (\rho_1, \rho_2)
$$
and a positive number $R > 0$, we define an elongated homotopy of
$\HH_1 \# \HH_2$
$$
\HH_1 \#_{(\rho;R)} \HH_2 =\{ H_{(\rho;R)}(\tau) \}_{-\infty <
\tau < \infty}
$$
by
$$
H_{(\rho;R)}(\tau,t,x)
= \begin{cases} H_1(\rho_1(\tau + 2R) ,t,x) & \quad \tau \leq 0 \\
H_2(\rho_2(\tau - 2R), t, x) & \quad \tau \geq 0.
\end{cases}
\leqno (5.17)
$$
Note that
$$
H_{(\rho;R)} \equiv \begin{cases} H_\alpha  \quad & \text{for } \, \tau
\leq - (R_1+2R)\\
H_\beta \quad  & \text{for } \, -R \leq  \tau \leq R\\
H_\gamma \quad & \text{for } \, \tau \geq R_2 + 2R
\end{cases}
\leqno (5.18)
$$
for some sufficiently large $R_1, \, R_2 > 0$ depending on the
cut-off functions $\rho_1, \, \rho_2$ and the homotopies $\HH_1,
\, \HH_2$ respectively. {\it In particular this elongated homotopy
is always smooth, even when the usual glued homotopy $\HH_1\#
\HH_2$ may not be so.} We define the elongated homotopy
$j_1\#_{(\rho;R)} j_2$ of $j_1\# j_2$ in a similar way.

For an elongated homotopy $(j_1\#_{(\rho;R)} j_2, \HH_1 \#_{(\rho,
R)} \HH_2)$, we consider the associated perturbed Cauchy-Riemann
equation
$$
\begin{cases} \frac{\part u}{\part \tau} +
J_3^{\rho(\tau)}\Big(\frac{\part u}{\part t}
- X_{H_3^{\rho(\tau)}}(u)\Big) = 0\\
\lim_{\tau \to -\infty}u(\tau) = z^-,  \, \lim_{\tau \to
\infty}u(\tau) = z^+
\end{cases}
\leqno (5.19)
$$
with the condition (5.4). The following lemma will be used later.
\medskip

{\noindent\bf Lemma 5.6.} {\it Let $H_\alpha, \, H_\beta, \, H_\gamma$ be
given, and let $\HH_i$ for $i = 1,\, 2$ be homotopies between them
respectively. Fix a generic homotopies $j_i$ for $i=1, \, 2$. For
any given pair $\rho = (\rho_1,\rho_2)$ of monotone cut-off
functions and a positive number $R > 0$, we consider the elongated
homotopy
$$
j_1 \#_{(\rho;R)} j_2 = \{J_{(\rho;R)}(\tau)\}_{\tau \in \R},
\quad \HH_1 \#_{(\rho;R)} \HH_2 = \{H_{(\rho;R)}(\tau)\}_{\tau \in
\R}
$$
defined as above, and the associated equation (5.19). Then for any
finite energy solution of (5.19), we have the inequality
$$
\AA_{H_\gamma}(u(\infty)) - \AA_{H_\alpha}(u(-\infty)) \leq - \int
\Big|\dudtau\Big|^2_{J_{(\rho;R)}(\tau)} + (E^-(\HH_1) +
E^-(\HH_2)). \leqno (5.20)
$$
This can be rewritten also as an a priori energy bound
$$
\begin{aligned} \int \Big|\dudtau\Big|_{J_{(\rho;R)}(\tau)}^2 & \leq
 - \AA_{H_\gamma}(u(\infty)) + \AA_{H_\alpha}(u(-\infty)) \\
& \qquad + (E^-(\HH_1) + E^-(\HH_2)).
\end{aligned}
\leqno (5.21)
$$
}
\begin{proof} This is a special case of Lemma 5.2 because
$$
E^-(\HH_1 \# \HH_2) = E^-(\HH_1) + E^-(\HH_2)
$$
from Lemma 5.5, and  $\HH_1 \#_{(\rho;R)} \HH_2 =
\{H_{(\rho;R)}(\tau)\}$ itself is an elongation of the homotopy
$\HH_1 \# \HH_2$ corresponding to the cut-off function
$$
\rho(\tau) = \begin{cases} \frac{1}{2}\rho_1(\tau + 2R) & \quad \tau \leq
0
\\
\frac{1}{2} + \frac{1}{2}\rho_2(\tau - 2R) & \quad \tau \geq 0
\end{cases}
\leqno (5.22)
$$
which remains monotone when $\rho_1, \, \rho_2$ are monotone.
\end{proof}

Now let $u_1$ and $u_2$ be given solutions of (5.3)-(5.4) associated
to $\rho_1$ and $\rho_2$ respectively. If we define the pre-gluing
map $u_1 \#_R u_2$ by the formula
$$
u_1\#_R u_2(\tau,t) =
\begin{cases} u_1(\tau + 2R,t)  & \quad\text{for }\, \tau \leq -R \\
u_2(\tau - 2R, t) & \quad\text{for }\, \tau \geq R
\end{cases}
$$
and a suitable interpolation between them by a partition of unity
on the region $ -R \leq \tau \leq R$, the assignment defines a
diffeomorphism
$$
(u_1, u_2, R) \to u_1 \#_R u_2 \leqno (5.23)
$$
from
$$
\MM\Big(j_1,\HH_1;[z_1,w_1],[z_2,w_2]\Big) \times
\MM\Big(j_2,\HH_2;[z_2,w_2],[z_3,w_3]\Big) \times (R_0, \infty)
$$
onto its image, provided $R_0$ is sufficiently large. Denote by
$\overline \part_{(\HH,j;\rho)}$ the corresponding perturbed
Cauchy-Riemann operator
$$
u \mapsto \frac{\part u}{\part \tau} +
J_3^{\rho(\tau)}\Big(\frac{\part u}{\part t} -
X_{H_3^{\rho(\tau)}}(u)\Big)
$$
acting on the maps $u$ satisfying the asymptotic condition imposed
in (5.19) and fixed homotopy condition $[u] = C \in
\pi_2(z^-,z^+)$. By perturbing $u_1\#_R u_2$ by the amount that is
smaller than the error for $u_1 \# _R u_2$ to be a genuine
solution, i.e., less than a weighted $L^p$-norm, for $p > 2$,
$$
\|\overline \part_{(\HH,j;\rho)}(u_1\#_{(\rho;R)} u_2)\|_p
$$
in a suitable weighted $W^{1,p}$ space of $u$'s (see \cite{Fl1,
Fl2}), one can construct a unique genuine solution near $u_1 \#_R
u_2$. By an abuse of notation, we will denote this genuine solution
also by $u_1 \#_R u_2$. Then the corresponding map defines an
embedding
$$
\begin{aligned} \MM\Big(j_1,\HH_1;[z_1,w_1],[z_2,w_2]\Big) & \times
\MM\Big(j_2,\HH_2;[z_2,w_2],[z_3,w_3]\Big)
\times (R_0, \infty) \to \\
& \to  \MM\Big(j_1\#_{(\rho;R)} j_2,\HH_1\#_{(\rho;R)}
\HH_2;[z_1,w_1],[z_3,w_3]\Big).
\end{aligned}
$$
Especially when we have
$$
\mu_{H_\beta}([z_2,w_2]) - \mu_{H_\alpha}([z_1,w_1]) =
\mu_{H_\gamma}([z_3,w_3]) - \mu_{H_\beta}([z_2,w_2]) = 0
$$
 both $\MM(j_1,\HH_1;[z_1,w_1],[z_2,w_2])$ and
$\MM(j_2,\HH_2;[z_2,w_2],[z_3,w_3])$ are compact, and so consist
of a finite number of points. Furthermore the image of the above
mentioned embedding exhausts the `end' of the
$$
\MM\Big(j_1\#_{(\rho;R)} j_2,\HH_1\#_{(\rho;R)}
\HH_2;[z_1,w_1],[z_3,w_3]\Big)
$$
and the boundary of its compactification consists of the broken
trajectories
$$
u_1\#_{(\rho; \infty)} u_2 = u_1 \#_\infty u_2.
$$
This then proves the gluing identity
$$
h_{\HH_1\#_{(\rho;R)} \HH_2} = h_{(\HH_1;\rho_1)} \circ
h_{(\HH_2;\rho_2)} \leqno (5.24)
$$
(modulo any filtration order we want) for a sufficiently large $R >
0$, suppressing the $j$-dependence. Here we remind the readers that
the homotopy $\HH_1\#_{(\rho;R)} \HH_2$ itself is an elongated
homotopy of the glued homotopy $\HH_1 \# \HH_2$.

Now we study the chain map $h_{(\HH,j;\rho)}$ more closely,  and
compare the maps for two different choices of the cut-off function
$\rho$ in relation to the level changes of the transferred Floer
cycles. Suppose that $\rho_1$ and $\rho_2$ are two cut-off
functions. We choose a one-parameter family
$$
\kappa \mapsto [0,1] \to (1-\kappa)\rho_1 + \kappa \rho_2 :=
\rho(\kappa)
$$
and consider the homotopy of homotopies of Hamiltonians
$$
(\kappa,\tau) \mapsto H(\rho(\kappa)(\tau),\cdot, \cdot).
$$
By considering the family $\overline j$ and $\overline \HH$
associated to this homotopy, the chain homotopy map $H_{(\overline
j, \overline \HH; \overline \rho)}$ is constructed by counting
solutions of (5.3) at non-regular, but {\it parametrically
regular} parameters $\kappa$. We will elaborate this remark in the
proof of the following main result in this section.
\medskip

{\noindent\bf Proposition 5.7.} {\it Let $(\HH,j)$ be a given
homotopy between $(H^0,J^0)$ and \linebreak $(H^1,J^1)$. Let
$\rho_1, \, \rho_2$ be two cut-off functions and $\overline\rho$ be
the homotopy $(5.24)$ between them. Let $(\overline j, \overline
\HH)$ be the associated homotopy of homotopies over $\kappa \in
[0,1]$. Then we have
$$
h_{(\HH,j;\rho_1)} - h_{(\HH,j;\rho_2)} =
\part_{(J^1,H^1)}\circ H_{(\overline \HH, \overline j;\overline\rho)}
+ H_{(\overline \HH, \overline j;\overline\rho)}\circ
\part_{(J^0,H^0)}.
$$
Furthermore the inequality
$$
\lambda_{H^1}(H_{(\overline \HH, \overline j;\overline\rho)}
(\alpha)) \leq \lambda_{H^0}(\alpha) + E(\HH) \leqno (5.25)
$$
holds for any Floer cycle $\alpha$ of $(J^0,H^0)$.
}
\begin{proof} It remains to prove (5.25). The map $H_{(\overline j,
\overline \HH;\overline\rho)}$ is determined by the number of the
pairs
$$
(\kappa, u)
$$
(See \cite{Fl3}) where $\kappa \in (0,1)$ and $u$ is a solution of
$$
\begin{cases} \frac{\part u}{\part \tau} +
J^{\rho_\kappa(\tau)}\Big(\frac{\part u}{\part t}
- X_{H^{\rho_\kappa(\tau)}}(u)\Big) = 0\\
\lim_{\tau \to -\infty}u(\tau) = z^-,  \lim_{\tau \to
\infty}u(\tau) = z^+
\end{cases}
$$
satisfying
$$
w^+ \sim w^- \# u;\quad \mu_{\HH_\kappa}([u]) = -1.
$$
Therefore we can estimate the level change
$$
\AA_{(H_\beta,J_\beta)}([z^+,w^+]) -
\AA_{(H_\alpha,J_\alpha)}([z^-,w^-])
$$
applying (5.8) for the function $\rho = \rho_k = (1-\kappa)\rho_1
+ \kappa\rho_2$. Once we have this, (5.25) follows in the same way
(5.16) was derived. This finishes the proof.
\end{proof}

\section*{\bf \S 6. Structure of  Novikov Floer cycles in a Cerf family}

In this section, we prove a general structure theorem, Theorem
6.9, of Novikov Floer cycles for any Cerf homotopy
$\HH=\{H(\eta)\}$ of Hamiltonian functions.

Let $H$ be  any one-periodic Hamiltonian and consider the
perturbed Cauchy- Riemann equation
$$
\frac{\part u}{\part \tau} + J \Big (\frac{\part u}{\part t} -
X_H(u)\Big ) = 0 \leqno (6.1)
$$
for a generic $J$. We call a solution $u$ {\it stationary} if it
is $\tau$-independent. We define
$$
A_{(J,H)} := \inf\Big\{ \int \Big| {\part u\over
\part \tau}\Big |^2_J  \,  \Big|  \, u \, \text{ satisfies (6.1)
and is not stationary } \Big\}.
$$
The positivity of $A_{(J,H)}$ is an easy consequence of the Gromov
compactness type theorem (see \cite{Oh6} for details of such a
proof).
\medskip

We also introduce the following invariant of a compact family $K
\subset \JJ_\omega$ of compatible almost complex structures. Let
$$
K;[0,1]^n \to  \JJ_\omega
$$
be a $n$-parameter family in the
$C^1$-topology, and define $A(\omega;K)$ be the constant
$$
\begin{aligned} A(\omega;K) = \inf_{\kappa \in [0,1]^n}
\Big\{A(\omega,J(\kappa)) \Big\}.
\end{aligned}
\leqno (6.2)
$$
This is always positive (see \cite{Oh6} for the proof in a similar
context), and enjoys the following lower semi-continuity property.
\medskip

{\noindent\bf Proposition 6.1.} {\it $A(\omega;K)$ is lower semi-continuous
in $K$. In other words, for any given $K$ and $0 < \e <
A(\omega;K)$, there exists some $\delta = \delta(K,\e) > 0$ such
that for any $K'$ with $\|K' - K\|_{C^1} \leq \delta$ we have
$$
A(\omega;K') \geq A(\omega;K) - \e.
$$
}
\begin{proof} The arguments used in the proof of this proposition
is similar to the one used in \cite[section 4]{Oh6}.

Suppose to the contrary that there exists some $0 < \e <
A(\omega;K)$ for which there are sequences $\delta_k \to 0$, $K_k$
with $\|K_k - K\|_{C^1} \leq \delta_k$ and
$$
A(\omega;K_k) < A(\omega;K) - \e
$$
for all $k$. By the definition of $A(\omega;K_k)$, this implies
that there exist non-constant $J_{k,t_k}$-holomorphic spheres
$w_k$ for $t_k \in [0,1]$ such that we have
$$
0< E_{J_{k,t_k}}(w_k) < A(\omega;K) - \e. \leqno (6.3)
$$
By choosing a subsequence, we may assume that $t_k \to t_\infty
\in [0,1]$ and $J_{k,t_k} \to J_{t_\infty}$ in the $C^1$-topology.
By the energy bound (6.3) for $w_k$ and since $K_k \to K$ in the
$C^1$ topology, we can produce a $J_\infty$-holomorphic map
$$
w_\infty = \sum_\ell w_{\infty, \ell}
$$
whose total energy satisfies
$$
E_{J_{t_\infty}}(w_\infty) \leq A(\omega;K) - \e. \leqno (6.4)
$$
Then by definition of $A(\omega;K)$, (6.4) implies all the
components of $w_\infty$ must be constant. By choosing a
subsequence if necessary, it follows that $w_k$ converges to a
constant map, say $p \in M$, in the $C^1$-topology, and also
$$
\lim_{k \to \infty} E_{J_{k,t_k}}(w_k) = 0.
$$
In particular, the image of $w_k$ is contained in a (contractible)
Darboux neighborhood $U$ of $p$ and so we have
$$
\int w_k^*\omega = 0.
$$
On the other hand, by the compatibility of $J_{k,t}$ to $\omega$,
we also have
$$
E_{J_{k,t_k}}(w_k) = \int w_k^*\omega,
$$
which in turn implies $E_{J_{k,t_k}}(w_k) = 0$ for all
sufficiently large $k$. But this contradicts the assumption that
$w_k$ are nonconstant. This finishes the proof.
\end{proof}

By a similar argument, we prove the following proposition.
\medskip

{\noindent\bf Proposition 6.2.} {\it Let $(J,H)$ be a pair with $H$
nondegenerate. Then for any given $0 < \e < \min\{A_{(J,H)},
A(\omega;J)\}$, there exists some $\delta_2 = \delta_2(J,H,\e) >
0$ such that for any $(J', H')$ with $\|(J',H') -
(J,H)\|_{C^\infty} \leq \delta$ we have
$$
A_{(J',H')} \geq \min\{A_{(J,H)}, A(\omega;J)\} - \e.
$$
}
\begin{proof} Suppose the contrary that there exists some $0\!<\!e\! <\!
\min\{A_{(J,H)}, A(\omega;J)\}$ for which there are sequences
$\delta_k \to 0$, $(J_k,H_k)$ with $\|(J_k,H_k)- (J,H)
\|_{C^\infty} \leq \delta_k$ and
$$
A_{(J_k,H_k)} < \min\{A_{(J,H)}, A(\omega;J)\} - \e
$$
for all $k$. By the definition of $A_{(J,H)}$, this implies that
there exist solutions $u_k$ of (6.1) for $(J_k,H_k)$ such that we
have
$$
0< E_{(J_k,H_k)}(u_k) < \min\{A_{(J,H)}, A(\omega;J)\} - \e. \leqno(6.5)
$$
By the energy bound (6.5) for $u_k$ and since $(J_k,H_k) \to
(J,H)$ in the $C^\infty$ topology, we can produce a
cusp-trajectory
$$
u_\infty = \sum_\ell u_{\infty, \ell}
$$
of (6.1) for $(J,H)$ with
$$
E_{(J,H)}(u_{\infty}) \leq \min\{A_{(J,H)}, A(\omega;J)\}- \e.
$$
Therefore by the definitions of $A_{(J,H)}$ and $A(\omega;J)$, it
follows that all the summand $u_{\infty,\ell}$ in $u_\infty$ must
be trivial in that all the principal components are stationary and
all the bubble components are constant. In particular, we have
$E_{(J,H)}(u_\infty) = 0$ and $u_k$ uniformly converges to a
periodic orbit $z_\infty$ of $H$ in the (fine) $C^\infty$
topology. These then imply
$$
\lim_{k \to \infty} E_{(J_k,H_k)}(u_k) = 0. \leqno (6.6)
$$
On the other hand, since we assume that $H$ is nondegenerate and
since  we can make $H_k$ arbitrarily $C^\infty$-close to $H$, by
taking sufficiently large $k$, $H_k$ are nondegenerate and there
is a canonical one-one correspondence between $\text{Per}(H)$ and
$\text{Per}(H_k)$ for each $k$. In particular, there exists some
$k_0 \in \Z_+$ and $c_1
> 0$ such that
$$
d_{C^0}(z,z') \geq c_1
$$
for any periodic orbits $z \neq z'$ of $H_k$ with $k \geq k_0$.
Since $J_k \to J$ in the $C^\infty$ topology, this and (6.6) imply
that all $u_k$ must satisfy
$$
u_k(-\infty) = u_k(\infty): \leqno (6.7)
$$
For otherwise it would imply that there exists $c_2 > 0$
independent of $k$ such that $E_{(J_k,H_k)}(u_k) > c_2$ (See Lemma
A.3 for the proof of this statement), which would contradict
(6.6).

We now recall the following lemma from \cite[Lemma 4.5]{Oh6}
\medskip

{\noindent\bf Lemma 6.3.} {\it Suppose that $u: \R \times S^1 \to M$ is any
finite energy solution of
$$
\begin{cases}
\dudtau + J\Big(\dudt - X_H(u)\Big) = 0 \\
E_J(u) = \int |{\part u \over \part \tau }|_{J_t}^2 < \infty.
\end{cases}
$$
that satisfies
$$
u(-\infty,t) = u(\infty, t).
$$
Then $\int_{\R \times S^1} u^*\omega$ converges, and we have
$$
E_J(u) = \int_{\R \times S^1} u^*\omega.
$$
}

Going back to our proof of Proposition 6.2, (6.7) and this lemma
imply
$$
E_{(J_k,H_k)}(u_k) = \int u_k^*\omega. \leqno (6.8)
$$
On the other hand, since $u_k$ satisfies (6.7) and converges to a
periodic orbit $z_\infty$ in the (fine) $C^\infty$ topology, it
defines a cycle homologous to the one dimensional cycle $z_\infty$
and hence $\int u_k^*\omega = 0$. In turn (6.8) implies we must
have
$$
E_{(J_k,H_k)}(u_k) = 0
$$
which contradicts the assumption that $u_k$ are non-stationary.
This finishes the proof.
\end{proof}

For a choice of $(\HH,j)$ with $\HH
\in\PP^{Cerf}(\HH_m(M);H_1,H_2)$, we define
$$
\begin{aligned} \NN t(\HH,j)& = \{\eta \in [0,1] \backslash \SS ing(\HH)
\mid \text{the pair $(J(\eta),H(\eta))$} \\
& \qquad \text{ has a solution of (4.3) with the Fredholm index 0}
\}.
\end{aligned}
\leqno (6.9)
$$
The following proposition was proved in section 4 (Theorem 4.6).
We rephrase Theorem 4.6 (1) to manifest what the transversality
means in a more concrete context. It provides a structure theorem
of $\NN t(\HH,j)$ for a generic choice of $j$ for a given homotopy
$\HH$. We will consider only the case where $\HH$ is a Cerf
homotopy in the statement.
\medskip

{\noindent\bf Proposition \& Definition 6.4.} {\it Let $\HH$ be a Cerf
homotopy and $\NN t(\HH,j)$ be as above. Then for a generic choice
of the homotopy $j$ of almost complex structures, the followings
hold: \begin{itemize}

\item[(1)] At each $\eta \in \NN t(\HH,j)$, the equation {\rm (4.3)-(4.4)}
for $(J,H) = (J(\eta),H(\eta))$ has exactly one pair $z^+, \, z^-
\in \text{Per}(H(\eta))$, and one non-trivial solution $u$
connecting $z^\pm$ which has Fredholm index 0.

\item[(2)] Write $\NN t(\HH,j)$ as an increasing union
$$
\NN t(\HH,j) = \bigcup_{N = 1}^\infty \NN t_N(\HH,j) \leqno (6.10)
$$
where we define
$$
\begin{aligned} \NN t_N(\HH,j) = \{ \eta \in \NN t(\HH,j) \mid\,  &
\text{the
unique solution } \, u \\
& \quad \text{ satisfies }\, E_{(J(\eta),H(\eta))} \leq N \}.
\end{aligned}
$$
Then $\NN t_N(\HH,j)$ is a compact zero dimensional manifold for
each given $N\in \Z_+$ and in particular a finite subset of
$[0,1]$. In particular, there are only {\it countably many} such
points $\eta \in [0,1] \backslash \SS ing(\HH)$.
\end{itemize}
We call a corresponding homotopy $(\HH,j)$ a {\em Floer homotopy}
and call any point $\eta$ lying in $[0,1] \setminus (\SS ing(\HH)
\cup \NN t(\HH,j))$ a {\em Floer point} of the Floer homotopy
$(\HH,j)$.}
\medskip

For a Floer homotopy $(\HH,j)$, the Floer homology
$HF_*(H(\eta),J(\eta))$ is well-defined at any Floer point $\eta$.
Furthermore the subset $\SS ing(\HH) \cup \NN t(\HH,j)$ is at most
countable and so the subset $[0,1] \setminus (\SS ing(\HH) \cup
\NN t(\HH,j))$ of Floer points in $[0,1]$ is residual and dense.

Next we compare the levels of Novikov Floer {\it cycles} and of
their boundaries. Suppose that $\eta, \, \eta'$ are sufficiently
close so that the associated pairs are defined. Let $(z,z')$ be an
associated pair of $(H(\eta),H(\eta'))$. We denote by
$$
\MM_{\eta\eta'}^\rho(z,z';[u^{can}_{zz'}]): =
\MM\Big((\HH_{\eta\eta'},j_{\eta\eta'};\rho);z,z':[u^{can}_{zz'}]\Big)
\leqno (6.11)
$$
the set of solutions $u$ of (5.3)-(5.4) homotopic to the canonical
cylinder $u^{can}_{zz'}$ relative to the boundary, and satisfying
$$
u(-\infty) = z \in \text{Per}(H(\eta)), \quad u(\infty) = z' \in
\text{Per}(H(\eta'))
$$
A straightforward calculation using (5.9) gives rise to the
following proposition whose proof we leave to the readers.
\medskip

{\noindent\bf Proposition 6.5.} {\it Let $F$ and $H$ be
nondegenerate Hamiltonians and $\HH$ be a homotopy between them. Fix
any, monotone or not, cut-off function $\rho$. Let $u_1, \, u_2$ be
solutions satisfying $(5.3)$ with fixed asymptotic conditions
$$
u_i(-\infty) = [z^-,w^-], \quad u_i(\infty) = [z^+,w^+]
$$
for some $w^-$ and be homologous to each other in that the
compactified torus $u_1 \# \overline u_2$ satisfies
$$
u_1 \# \overline u_2 \sim 0.
$$
Then we have
$$
|E_{(\HH,j;\rho)}(u_1) -E_{(\HH,j;\rho)}(u_2)| \leq E(\HH). \leqno(6.12)
$$
In particular, when $\HH \equiv H, \, j \equiv J$, we have
$$
E_{(J,H)}(u_1) = E_{(J,H)}(u_2).
$$
}

Let $\eta \in [0,1]$ be a point for which $H(\eta)$ is
nondegenerate. We then introduce the following function
$$
\begin{aligned} \ell_{(\HH,j)}^\rho(\eta, \eta'): = \max_{z, \, z', u}
\Big\{\int \Big| {\part u \over \part
\tau}\Big|^2_{J^{\rho(\tau)}}\,  \Big| \quad & u \in
\MM_{\eta\eta'}^\rho(z,z';[u^{can}_{zz'}]), \\
& (z,z') \, \text{is an associated pair} \Big\} \\
\end{aligned}
\leqno (6.13)
$$
defined for $\eta'$ such that  $|\eta - \eta'|$ is sufficiently
small and $H(\eta')$ is nondegenerate.

And for a given $\eta_\infty \in \SS ing(\HH)$, we consider the
bifurcation pair $(z^+(\eta)$, $z^-(\eta))$ for the pair
$(J(\eta),H(\eta))$ as $\eta \to \eta_\infty \in \SS ing(\HH)$. We
define the function
$$
\ell_{(\HH,j; \eta_\infty)}(\eta): = \max_{u}\Big\{\int \Big|
{\part u \over \part \tau}\Big|^2_{J(\eta)}\, \Big| \, u \in
\MM(z^+(\eta),z^-(\eta);[u^{can}]) \Big\} \leqno (6.14)
$$
defined for $\eta$ with $|\eta - \eta_\infty|$ sufficiently small,
where
$$
\MM(z^+(\eta),z^-(\eta);[u^{can}]) :=
\MM\Big(H(\eta),J(\eta);z^+(\eta),z^-(\eta);[u^{can}]\Big).
$$
\medskip

{\noindent\bf Proposition 6.6.} {\it Let $\eta, \, \eta'$ be nondegenerate
points such that $|\eta -\eta'|$ is sufficiently small and
$(z,z')$ be an associated pair such that
$\MM_{\eta\eta'}^{\rho_1}(z,z';[u^{can}_{zz'}])$ is nonempty for a
cut-off function $\rho_1$. Then there exists a function
$C_{(\HH,j)}(r)$ depending only on $(\HH,j)$ such that \begin{itemize}

\item[(1)] $C_{(\HH,j)}(r) \to 0$ as $r \to 0$ and

\item[(2)] $\ell_{(\HH,j)}^\rho( \eta, \eta')
\leq C_{(\HH,j)}(|\eta -\eta'|).$
\end{itemize}

Similar estimate also holds for the function
$\ell_{(\HH,j;\eta_\infty)}$ uniformly over $\eta_\infty \in \SS
ing(\HH)$.
}
\begin{proof} Let $u \in \MM_{\eta\eta'}(z,z';[u^{can}_{zz'}])$.
Using the fact that $u$ is a solution of (5.3), we derive
$$
\int \Big|\dudtau \Big|_{J^{\rho_1(\tau)}}^2 \leq -
\AA_{H(\eta')}([z', w\# u^{can}]) + \AA_{H(\eta)}([z,w]) +
E(\HH_{\eta\eta'}) \leqno (6.15)
$$
from (5.16). On the other hand, from the definition of the action
functional, we have
$$
\begin{aligned} \AA_{H(\eta')}([z', w\# u^{can}])& = -\int_{w \# u^{can}}
\omega - \int_0^1 H(\eta')(t, z'(t)) \, dt \\
\AA_{H(\eta)}([z,w]) & =  -\int_{w} \omega - \int_0^1 H(\eta)(t,
z(t)) \, dt
\end{aligned}
$$
and so we get
$$
\begin{aligned} \AA_{H(\eta')}([z', w\# u^{can}]) & -
\AA_{H(\eta)}([z,w]) \\
& = -\int_{u^{can}} \omega -\int_0^1 \Big(H(\eta')(t,z'(t)) -
H(\eta)(t,z(t))\Big) \, dt.
\end{aligned}
\leqno (6.16)
$$
Substituting (6.16) into (6.15), we get
$$
\int \Big|\dudtau \Big|_{J^{\rho_1(\tau)}}^2 \leq \int_{u^{can}}
\omega + \int_0^1 \Big(H(\eta')(t,z'(t)) - H(\eta)(t,z(t))\Big) \,
dt + E(\HH_{\eta\eta'}). \leqno (6.17)
$$
This can be estimated above by
$$
\text{Area}_g(u^{can}_{zz'}) + \Big|\int_0^1 \Big(H(\eta)(t,z'(t))
- H(\eta)(t, z(t))\Big) \, dt \Big| + E(\HH_{\eta\eta'}). \leqno(6.18)
$$
On the other hand, one can easily estimate
$$
\text{dist}_{C^1}(z,z') \leq \widetilde C_1(\|H(\eta') -
H(\eta)\|_{C^2}) \leqno (6.19)
$$
for the associated pair $(z,z')$ with $z \in \text{Per}(H(\eta))$
and $z' \in \text{Per}(H(\eta'))$ by a function $\widetilde
C_1(r)$ such that $\widetilde C_1(r) \to 0$ as $r \to 0$.

 From the explicit expression (2.6) of the canonical thin cylinder,
(6.19) immediately implies
$$
\text{Area}_g(u^{can}_{zz'}) \leq \widetilde
C_1(\|H(\eta') - H(\eta)\|_{C^2}) \leq C_1(|\eta - \eta'|)
$$
by another function $\widetilde C_1$ satisfying $C_1(r) \to 0$ as
$r \to 0$. On the other hand, (6.19) obviously implies
$$
\Big|\int_0^1 \Big(H(\eta)(t,z'(t)) - H(\eta)(t, z(t))\Big) \, dt
\Big| \leq C_2(|\eta-\eta'|)
$$
for a function $C_2 = C_2(r)$ with the same property.

Finally since $\HH$ is smooth and $[0,1]\times M$ is compact, we
 have
$$
E(\HH_{\eta\eta'}) \leq C_3(|\eta - \eta'|)
$$
for some positive function $C_3$ satisfying $C_3(r) \to 0$ as
$r\to 0$. Now, defining $C_{(\HH,j)}$ by
$$
C_{(\HH,j)}(r) = C_1(r) + C_2(r) + C_3(r),
$$
we have finished the proof. The statement about
$\ell_{(\HH,j;\eta_\infty)}$ follows immediately from (6.17) and
(6.19) applied to the case of $\eta' = \eta$ and $z' = z^+(\eta),
\, z = z^-(\eta)$. Note that in the latter case, the term
$E(\HH_{\eta\eta'})$ drops out.
\end{proof}

Now we are ready to state the main theorem of this section, which
concerns a general structure of Novikov Floer cycles over a Floer
homotopy $(\HH,j)$. This theorem in particular says that for a
Floer homotopy $(\HH,j)$, if there is a Floer trajectory issued
{\it at a peak of the given Floer cycle}, all such trajectories
must be `long' uniformly over the interval $[0,1]$ in that its
energy has uniform positive lower bound. We like to emphasize that
such a property strongly depends on the {\bf cycle} property of
the chains.

This theorem is closely related to the well-known picture arising in
the First Cancellation Theorem in the Morse theory \cite{Mi}. More
precisely, consider a one-parameter family of smooth functions $f_s$
such that $f_s, \, 0 \leq s < 1$ are all Morse but $f_1$ contains a
unique degenerate critical points $p$. Let $(p^+(s),p^-(s))$ for $s
\in [0,1]$ with $p^+(1) = p^-(1)$ be the continuous family of
cancellation pair of critical points of Morse functions $f_s$ of
Morse indices of $(k+1, k)$ that appear in the First Cancellation
Theorem \cite{Mi}, the cancellation theorem implies that $p^+(s)$
cannot contribute to the (Morse) homology $H_{k+1}(M)$. Once one
gets rid of the critical point $p^+(s)$, one can easily see that
there is a constant $A > 0$, independent of $s \in [1-\e, 1)$ such
that all the gradient trajectories issued at any critical point $p$
have length greater than $A$. In particular, the level of the
boundary $\part(p)$ is less than the level of $p$ at least by $A$.
\medskip

{\noindent\bf Theorem 6.7.} {\it Let $\eta \in [0,1] \setminus (\SS ing(\HH)
\cup \NN t(\HH,j))$ be any Floer point. Then there exist constants
$ \delta_1(\HH) > 0$ and $A^1_{(\HH,j)}
> 0$ independent of $\eta$ and $[z,w]$, depending only on
$(\HH,j)$, such that the followings hold :

\begin{itemize}

\item[(1)] if $[z,w] \in \text{Crit}\AA_{H(\eta)}$, $[z',w']
\in\part_{(H(\eta),J(\eta)}([z,w])$, and
$$
\AA_{H(\eta)}([z,w]) - \AA_{H(\eta)}([z',w']) < A^1_{(\HH,j)} \leqno(6.20)
$$
then we have $d(\eta, \SS ing(\HH))< \delta_1(\HH)$ and
$$
[z,w] = [z^+(\eta),w], \quad [z',w'] = [z^-(\eta), w\# u^{can}].
\leqno (6.21)
$$
\item[(2)] there exists another constant $\delta < \delta_1(\HH)$ such
that if $d(\eta,\SS ing(\HH)) < \delta$, we can write
$\part_{(J(\eta_i),H(\eta_i))}([z^+(\eta_i),w_i^+])$ as
$$
\part_{(J(\eta_i),H(\eta_i))}([z^+(\eta_i),w_i^+])
= [z^-(\eta_i), w_i^+\# u_i] + \beta_i \leqno (6.22)
$$
for a chain $\beta$ satisfying
$$
\lambda_{H(\eta)}(\beta) \leq \AA_{H(\eta)}([z^-(\eta),w\#
u^{can}] - A^1_{(\HH,j)}.
$$
\item[(3)] Furthermore the constant $\min\{A^1_{(\HH,j)},
A(\omega;j)\}$ can be chosen to be lower semi-continuous in $j$.
\end{itemize}
}
\begin{proof}  We first note that
$$
\AA_{H(\eta)}([z,w]) >
\lambda_{H(\eta)}(\part_{(J(\eta),H(\eta))}([z,w]))
$$
in general as the Floer boundary map always decreases the level.
Furthermore if $d(\eta,\SS ing(\HH)) \geq \delta_1$, where
$d(\eta,\SS ing(\HH))$ is the distance of $\eta$ to $\SS
ing(\HH)$, Lemma A.3 (or rather its proof) implies that there
exists $A=A(\HH,j,\delta_1) > 0$ such that
$$
\AA_{H(\eta)}([z,w]) - \AA_{H(\eta)}([z',w']) \geq A
$$
for any $[z,w] \in \text{Crit}\AA_{H(\eta)}$ and $[z',w'] \in
\part([z,w])$.

Therefore from now on, we assume
$$
d(\eta,\SS ing(\HH)) < \delta_1
$$
and so the associated bifurcation branches $(z^+(\eta),
z^-(\eta))$ issued at one of the points in $\SS ing(\HH)$ are
defined. We define
$$
\begin{aligned} A^1_{(\HH,j)} : =  \inf_{\eta, u} \Big\{
E_{(H(\eta),J(\eta))}(u) \, & \Big|\,
\eta \in [0,1] \setminus (\SS ing(\HH) \cup \NN t(\HH,j)), \\
& \qquad 0 < d(\eta,\SS ing(\HH)) < \delta_1(\HH), \\
& \qquad u
\not \in \MM(\HH,j;\rho;z^+(\eta),z^-(\eta);[u^{can}]) \Big\}.
\end{aligned}
$$
By definition, $A^1_{(\HH,j)}$ satisfies (6.20). The rest of the
proof will show that  $A^1_{(\HH,j)} > 0$ and $A^1_{(\HH,j)}$ is
lower semi-continuous in $(\HH,j)$.

Suppose that there exists a sequence $\eta_i \in [0,1] \backslash
\SS ing(\HH) \cup \NN t(\HH,j)$ with $\eta_i \to \eta_\infty \in
(0,1)$ such that there exist critical points $[z_i,w_i] \in
\text{Crit}\AA_{H(\eta_i)}$ and $[z_i',w_i'] \in
\part_{(J(\eta_i),H(\eta_i))}([z_i,w_i])$ satisfying
$$
\AA_{H(\eta_i)}([z_i,w_i]) - \AA_{H(\eta_i)}([z_i',w_i']) \to 0.
\leqno (6.23)
$$
It will be enough to prove that (6.21) must hold under the given
assumption, if $i$ is sufficiently large.

The rest of the proof will be divided into 5 steps.
\medskip
\n{\it Step 1: $\eta_\infty \in \SS ing(\HH)$.}
\smallskip

Suppose to the contrary that $\eta_\infty \in [0,1] \setminus \SS
ing(\HH)$. Then the lower semi-continuity of $A_{(J,H)}$ stated in
Proposition 6.2 implies
$$
A_{(H(\eta),J(\eta))} \geq
\frac{1}{2}A_{(H(\eta_\infty),J(\eta_\infty))} > 0
$$
for all $\eta$
sufficiently close to $\eta_\infty$. In turn this implies
$$
\AA_{H(\eta_i)}([z_i,w_i]) - \AA_{H(\eta_i)}([z_i',w_i']) =
E_{J(\eta_i)}(u) \geq  \frac{1}{2}A_{(J(\eta_\infty),
H(\eta_\infty))}
$$
for any pair $[z_i,w_i] \neq [z_i',w_i']$ for all sufficiently
large. This contradicts (6.23).

\medskip
\n{\it Step 2.}
\smallskip

By the definition of the Cerf homotopy $\HH$, there is exactly one
birth-death or death-birth point at each $\eta_\infty \in \SS
ing(\HH)$. Then we can repeat the proof of Proposition 6.2 with
minor modifications, and prove that except the critical points of
the form $[z^+(\eta), w]$, where $z^+(\eta)$ is the upper one of
the bifurcation pair $(z^+(\eta),z^-(\eta))$ issued at
$\eta_\infty$ as stated in Proposition 3.3, the levels of all
other critical points are dropped down by a positive amount, say
$B > 0$, under the action by the boundary map
$\part_{(J(\eta),H(\eta))}$, i.e.,
$$
\lambda_{H(\eta_i)}(\part_{(J(\eta_i),H(\eta_i))}([z_i,w_i])) \leq
\AA_{H(\eta_i)}([z_i,w_i]) - B.
$$
And $B$ can be chosen to be independent of any given Floer point
$\eta$. This proves that $[z_i,w_i]$ must be of the form
$[z^+(\eta_i), w_i]$.

\medskip
\n{\it Step 3: Finish-up of the proof of (1)}
\smallskip

We now analyze the term $\part_{H(\eta_i)}([z^+(\eta_i), w_i])$.
First, we claim that we have
$$
\AA_{H(\eta_i)}([z^+(\eta_i),w_i]) - \AA_{H(\eta_i)}([z,w]) \geq B
\leqno (6.24)
$$
for all the generators $[z,w] \in \part_{H(\eta_i)}([z^+(\eta_i),
w_i])$, {\it except } possibly $$[z^-,w^+_i\# u^{can}_{z^+z^-}],$$
if we choose $B$ smaller if necessary. We now make this statement
more precise.

Let $(\HH,j)$ be a Floer homotopy in the sense of Definition 6.4,
and $\eta \in [0,1]$ be a Floer point. We denote by $\alpha$ a
Floer cycle associated to $(J(\eta),H(\eta))$. Write
$$
\SS ing(\HH) =\{ s_1, \ldots, s_{k_1}\}.
$$
When $d(\eta,\SS ing(\HH))$ is sufficiently small, let
$(z^+(\eta), z^-(\eta))$ be the bifurcation pair issued from a
point $s_i \in \SS ing(\HH)$. We will prove the following lemma in
the Appendix, which will finish the proof of (1).
\medskip

{\noindent\bf Lemma 6.8.} {\it Suppose $(\HH,j)$ and $\eta$ are as above.
Let $[z^+(\eta), w]$ be the generator of $\alpha$ having the
highest level among the generators of the form $[z^+(\eta),w]$.
Then there exists a constant $B = B(\HH,j) > 0$ such that there
exists $\delta = \delta(\HH,j,\e)$ with
$$
0 < \delta < \min_{1\leq i\leq k_1-1}\{|s_{i+1} - s_i|\}
$$
such that for any $\eta$ with $dist(\eta,\SS ing(\HH)) < \delta$,
any element $u$ in
$$
\MM\Big(J(\eta),H(\eta); [z^+(\eta),w^+], [z, w]\Big)
$$
satisfies
$$
E_{(J(\eta), H(\eta))}(u) \geq B
$$
unless $[z,w] = [z^-,w_i\# u^{can}]$.
}

\medskip
\n{\it Step 4: Proof of $(6.22)$.}
\smallskip
To prove (2), we first study the moduli space
$$
\MM\Big(J(\eta_i),H(\eta_i); [z^+(\eta_i),w_i], [z^-(\eta_i),
w_i\# u^{can}]\Big)
$$
of solutions $u$ of (6.1) in the class prescribed by the condition
$$
u(-\infty) = z^+, \, u(\infty)= z^-, \quad w\# u^{can} \sim w\# u.
$$
We first note that the estimate for $\ell_{(\HH,j;\eta_\infty)}$
in Proposition 6.6 implies that all the elements in this moduli
space has small energy and so is `localized' near the unique
degenerate periodic orbit $z_\infty$ of $H(\eta_\infty)$. In other
words, the image of all the elements $u$ in
$\MM(J(\eta_i),H(\eta_i); [z^+(\eta_i),w^+], [z^-(\eta_i), w^+\#
u^{can}])$ is contained in a small neighborhood of that of
$z_\infty$, which can be chosen as small as we want by choosing a
sufficiently large $i$.

Then by considering the standard bifurcation picture near a generic
degenerate periodic orbit with respect to a suitable family
$J=\{J_t\}_{0 \leq t \leq 1}$ of almost complex structures (see
 \cite{Fl2, Lee} for some relevant explanations), there exists such a
family for which we have {\it precisely one} Floer trajectory $u_i$
from $z^+(\eta_i)$ to $z^-(\eta_i)$ that is homotopic to the
canonical thin cylinder. We refer to \cite[Theorem 9.9]{Lee} for the
precise statement and its proof. We like to remark that proving such
a statement directly involves highly technical analytical estimates
because it involves an analysis of the Floer moduli space near a
{\it degenerate} Hamiltonian $H(\eta_\infty)$. The proof was
outlined in \cite{Fl1} by Floer himself and later completed by
Yi-Jen Lee \cite{Lee}.

Once we have this analytical theorem at our disposal, applying the
standard cobordism argument over a homotopy from the above
mentioned $J$ to our $J(\eta_i)$, we have proven that the matrix
coefficient becomes
$$
\langle
\part_{(J(\eta_i),H(\eta_i))}([z^+(\eta_i),w_i]),[z^-(\eta_i),
w_i\# u^{can}] \rangle = 1.
$$
({\it Here we remind the readers that in this paper we are
assuming that $(M,\omega)$ is strongly semi-positive and so all
the matrix coefficients of the basic operators in the Floer
homology have integer coefficients. For the general $(M,\omega)$,
this matrix coefficient may become rational numbers but will not
still be zero.}) Therefore we have
$$
\part_{(J(\eta_i),H(\eta_i))}([z^+(\eta_i),w_i])
= [z^-(\eta_i), w_i\# u^{can}] + \beta_i \leqno (6.25)
$$
where $\beta_i$ does not have $[z^-(\eta_i), w_i\# u^{can}]$ as
one of its generators.

\medskip
\n{\it Step 5: Finish-up of the proof of $(2)$}
\smallskip

Now we estimate the level of $\beta_i$. But Lemma 6.8 implies
$$
\AA_{H(\eta_i)}([z^+(\eta_i),w_i]) - \AA_{H(\eta_i)}([z,w]) \geq B
\leqno (6.26)
$$
for any generator $[z,w] \in \beta_i$ where $B=B(\HH,j)$. Combined
with Step 4, this finishes the proof of (2).

\medskip
\n{\it Step 6: The lower semi-continuity of $A_{(\HH,j)}^1$.}
\smallskip

The proof of the lower semi-continuity of $\min\{A^1_{(\HH,j)},
A(\omega;j)\}$ in $j$ can be proceeded as the proof of proof of
Proposition 6.2 and so omitted. This finally finishes the proof of
Theorem 6.7.
\end{proof}

Theorem 6.7 gives rise to the following proposition.
\medskip

{\noindent\bf Proposition 6.9.} {\it Let $a\neq 0$ be a given quantum
cohomology class and denote by $\alpha$ a Floer cycle with
$[\alpha] = a^\flat$. Then there exists $0 < \delta_3 \leq
\delta_1$ with $\delta_3 = \delta_3(\HH,j)$ such that for any
Floer point $\eta$ satisfying $d(\eta,\SS ing(\HH)) < \delta_3$
for any tight Floer cycle $\alpha$ of $(H(\eta),J(\eta))$, any of
its peaks cannot have the form
$$
[z^+(\eta), w]
$$
for any bounding disc $w$.
}
\begin{proof} Suppose to the contrary that there exists a sequence
of the Floer points $\eta_k \to \eta_\infty \in \SS ing(\HH)$ such
that there exists a sequence of tight Floer cycles $\alpha_k$
whose peaks have the form $[z^+(\eta_k),w_k]$. By the definition
of Cerf homotopies, there exists at most one such peak of the form
$[z^+(\eta_k),w_k]$ provided $k$ is sufficiently large. We first
derive
$$
\lambda_{H(\eta_k)}(\part_{(H(\eta_k),J(\eta_k))}
([z^+(\eta_k),w_k])) = \AA_{H(\eta_k)}([z^-(\eta_k),w_k\#
u^{can}]) \leqno (6.27)
$$
from (6.25). Since $[z^+(\eta_k),w_k]$ is a peak of $\alpha_k$, we
have
$$
\lambda_{H(\eta_k)}(\alpha_k) = \AA_{H(\eta_k)}([z^+(\eta_k),w_k])
$$
and Proposition 6.5 implies
$$
\begin{aligned} \AA_{H(\eta_k)}([z^-(\eta_k),w_k\# u^{can}]) & \geq
\lambda_{H(\eta_k)}(\alpha_k) - \ell(\HH,j;\eta_0)(|\eta_k
- \eta_0|) \\
& = \rho(H(\eta_k);a) - \ell(\HH,j;\eta_0)(|\eta_k - \eta_0|).
\end{aligned}\leqno (6.28)
$$
From the proof of Theorem 6.7 or more specifically from (6.24), we
have
$$
\lambda_{H(\eta_k)}\Big(\part_{(H(\eta_k),J(\eta_k))}([z,w])\Big)
\leq \rho(H(\eta_k);a) - B \leqno (6.29)
$$
if $z \neq z^+(\eta_k)$ where $B$ is the same constant used in the
proof of Theorem 6.7. On the other hand if $z = z^+(\eta)$, then
obviously we have
$$
\begin{aligned} \lambda_{H(\eta_k)}\Big(\part_{(H(\eta_k),J(\eta_k))}&
([z^+(\eta_k),w_k])\Big) \\   =& \AA_{H(\eta_k)}([z^-(\eta_k),
w_k\# u^{can}]) > \AA_{H(\eta_k)}([z^-(\eta_k), w\# u^{can}]) \\
 =&\lambda_{H(\eta_k)}\Big(\part_{(H(\eta_k),J(\eta_k))}([z^+(\eta_k),w])\Big)
\end{aligned}\leqno (6.30)
$$
for any other $[z^+(\eta_k),w] \neq [z^+(\eta_k),w_k]$. Therefore
if we write
$$
\part_{(H(\eta_k),J(\eta_k))}(\alpha_k) \!=\!
\part_{(H(\eta_k),J(\eta_k))}([z^+(\eta_k),w_k]) +
\part_{(H(\eta_k),J(\eta_k))} (\alpha_i - [z^+(\eta_k),w_k])
$$
(6.29) and (6.30) imply
$$
\lambda_{H(\eta_k)}\Big(\part_{(H(\eta_k),J(\eta_k))} (\alpha_i -
[z^+(\eta_k),w_k])\Big)\! <\!
\lambda_{H(\eta_k)}\Big(\part_{(H(\eta_k),J(\eta_k))}([z^+(\eta_k),w_k])
\Big). \leqno (6.31)
$$
In particular, $\part_{(H(\eta_k),J(\eta_k)}(\alpha_k)$ cannot
vanish which contradicts that $\alpha_k$ is a cycle. Hence the
proof.
\end{proof}

\section*{\bf \S 7. Handle sliding lemma and sub-homotopies}

In this section, we recall another important element in the chain
level theory, the {\it handle sliding lemma} introduced in
\cite{Oh2}. We state the most natural version of the handle sliding
lemma which uses the {\it sub-homotopies} of the given homotopy
$\HH$ instead of the {\it piecewise-linear approximation} of $\HH$
which was used in \cite{Oh2}.

We start with the following lemma from \cite{Oh2}. Partly for the
reader's convenience and also because we need to add some important
points to the lemma, we repeat its proof here.
\medskip

{\noindent\bf Lemma 7.1. \cite[Lemma 5.1]{Oh2}} {\it Let $(\HH,j)$
be one-parameter family such that $\HH \in
\PP^{reg}(\HH_m(M);H_1,H_2))$. For each $\eta \in [0,1] \backslash
\SS ing(\HH)$, we define
$$
\begin{aligned} A^0_{(H(\eta),J(\eta))}  = & \inf_{u} \Big\{ \int
\Big|{\part u \over
\part \tau} \Big|^2_{J(\eta)}\Big|\quad  u \, \text{satisfies (6.1)}, \\
& \qquad \text{ is not stationary and }\, \text{Index } u = 0
\Big\}
\end{aligned}
$$
and
$$
A^{reg,0}_{(\HH,j)} = \inf_{s \in [0,1] \backslash \SS ing(\HH)}
A^0_{(H^s, J^s)}. \leqno (7.1)
$$
Then $A^{reg,0}_{(\HH,j)}$ is strictly positive.
}
\begin{proof} Suppose the contrary that $A^{reg,0}_{(\HH,j)} = 0$,
i.e., that there exists a sequence $\eta_k \in [0,1] \backslash
\SS ing(\HH)$ with $\eta_k \to \eta_\infty \in (0,1)$ and $u_k$
solutions of (6.1) for $(H(\eta_k),J(\eta_k))$ such that
$$
\int \Big|{\part u_k \over \part \tau}\Big|^2_{J(\eta_k)} \to 0,
\quad \text{Index } u_k = 0. \leqno (7.2)
$$
Then we must have, by choosing a subsequence if necessary,
$$
\eta_\infty \in \SS ing(\HH)
$$
and a degenerate periodic orbit $z_\infty$ of $\dot x =
X_{H(\eta_\infty)}(x)$ such that $u_k \to z_\infty$ uniformly and
so
$$
u_k(\infty), \, u_k(-\infty) \to z_\infty.
$$
Since $u_k(\pm\infty)$ are solutions of $\dot x =
X_{H(\eta_k)}(x)$, they must be the pair described in (1) right
above (3.1) in section 3.1 and hence
$$
\text{Index }(u_k) = \mu([z^+(\eta_k),w^+_k]) - \mu([z^-(\eta_k),
w^-_k] = 1.
$$
But this contradicts the index condition in the definition of
$A^0_{(H^s,J^s)}$ which finishes the proof.
\end{proof}

Next we define
$$
A^{sing}_{(\HH,j)} = \min_k \Big\{ A_{(H^{s_k}, J^{s_k})} ~|~ s_k
\in \SS ing(\HH)\Big\}. \leqno (7.3)
$$
This is again positive by a Gromov type compactness theorem. Now
we have the following crucial definition of a family version of
the constant $A_{(H,J)}$ suitable for our purpose. We define
$$
A^0_{(\HH,j)} = \min \Big\{ A^{reg,0}_{(\HH,j)},\,
A^{sing}_{(\HH,j)}, A(\omega;j)\Big \} \leqno (7.4)
$$
which we know is strictly positive.
\medskip

{\noindent\it Remark 7.2.} We would like to point out that all the
invariants $A^{sing}_{(\HH,j)}$, $A^{reg,0}_{(\HH,j)}$ and
$A^0_{(\HH,j)}$ are defined in terms of the Floer boundary
equation (4.3), not in terms of the continuity equation (4.5).
Furthermore, it follows from the same kind of proof as the proofs
of Proposition 6.1 and 6.2 that for a fixed $\HH$,
$$
\min\{A^0_{(\HH,j)}, A(\omega;j)\}
$$
is lower semi-continuous in $j$.

We state the following simple lemma, or rather an observation from
this remark, which follows immediately from the definitions of
$A^0_{(\HH,j)}$ and of the sub-homotopy in Definition 3.8. This
turns out to play an important role in our proof of the main theorem
later, and is one of the reasons why we have to use the
sub-homotopies of the given homotopy.
\medskip

{\noindent\bf Lemma 7.3.} {\it Let $(\HH,j)$ be a given homotopy. Then we
have the inequality
$$
A^0_{(\HH_{\eta_1\eta_2},j_{\eta_1\eta_2})} \geq A^0_{(\HH,j)}
\leqno (7.5)
$$
for the sub-homotopy $\HH_{\eta_1\eta_2}$ between any two Floer
points $0 \leq \eta_1 \leq \eta_2 \leq 1$. If $\eta_1 > \eta_2$
instead, then we have
$$
A^0_{(\HH_{\eta_1\eta_2}^{-1},j_{\eta_1\eta_2}^{-1})} \geq
A^0_{(\HH^{-1},j^{-1})}.
$$
}
\begin{proof} The proof is an immediate consequence of the
definitions of $A^0_{(\HH,j)}$ and of the sub-homotopy in general.
\end{proof}

We now recall the following {\it handle sliding lemma} from
\cite{Oh2}. We, however, add an important improvement from that of
\cite{Oh2}: here we used the sub-homotopy (4.8) of $\HH$ instead of
the linear homotopy that was used in \cite{Oh2}. It turns out that
this usage of sub-homotopies is the most natural and the optimal
choice, in that the constant associated to the given $(\HH,j)$ in
Definition 5.4 can be used for all its sub-homotopies. The proof
here is taken from \cite{Oh2}.
\medskip

{\noindent\bf Proposition 7.4 (The handle sliding lemma).} {\it Let
$j = \{J^\eta\} \in \PP^{Cerf}(j_\omega;\HH)$ be a $($two
parameter$)$ family of almost complex structures and
$\{H(\eta)\}_{0\leq \eta \leq 1}$ be a generic family of
Hamiltonians. Fix a cut-off function $\rho$. Let $A^0_{(\HH,j)}$ be
the constant defined in $(7.4)$ and let $\eta_1, \eta_2 \in [0,1]$.
\begin{itemize} \item[(1)] Then there exists $\delta = \delta(\HH,j)
> 0$ such that if $|\eta_1 - \eta_2| < \delta$, any finite energy
solution $u$ with
$$
\text{Index }u = 0
$$
of $(5.3)$ must either satisfy
$$
\int \Big| {\part u \over \part \tau}\Big|^2_{J^{\rho(\tau)}} \leq
\varepsilon(\delta) \leqno (7.6)
$$
or
$$
\int \Big| {\part u \over \part \tau}\Big|^2_{J^{\rho(\tau)}} \geq
A^0_{(\HH,j)} - \varepsilon(\delta) \leqno (7.7)
$$
where for $\varepsilon(\delta) \to 0$ as $0 < \delta \to 0$,
provided $\delta\leq \delta_0$.

\item[(2)] In addition, in the case of $(7.6)$, $u$ is homotopic to the
canonical cylinder $u^{can}_{z^-z^+}$ between $z^-$ and $z^+$ the
asymptotic periodic orbits of $u$. In particular, we have
$$
\int \Big| {\part u \over \part \tau}\Big|^2_{J^{\rho(\tau)}} \to
0
$$
as $\delta \to 0$.

\item[(3)]  Furthermore the same constant $A^0_{(\HH,j)}$ for $(7.7)$ can
be used for all the sub-homotopies $\HH_{\eta_1\eta_2}$ for any
two Floer points $\eta_1, \, \eta_2$.
\end{itemize}
}

\begin{proof} We prove this by contradiction. Suppose the contrary
that there exists some $\varepsilon > 0$, $\eta_1$ and $\eta_i$
with $\eta_i \to \eta_1$ as $i \to \infty$, and solutions $u_i$
that satisfy
$$
\text{Index }u_i = 0,
$$
$$
\frac{\part u_i}{\part \tau} + J^{\rho(\tau)} \Big(\frac{\part
u_i}{\part t} - X_{H^{\rho(\tau)}}(u_i)\Big) = 0 \leqno (7.8)
$$
and
$$
\varepsilon < \int \Big|{\part u_i \over \part
\tau}\Big|^2_{J^{\rho(\tau)}} < A_{(\HH,j)}^0 - \varepsilon. \leqno(7.9)
$$
In particular, the right half of (7.9) implies the uniform bound
on the energy of $u_i$. As $i \to \infty$, the equation (7.8)
converges to (6.1) with $(H,J) = (H(\eta_1),J(\eta_1))$. By
Gromov's type compactness theorem, we have a cusp trajectory
$$
u_\infty = \sum_k u_{\infty,k}
$$
which is the limit of a subsequence where each $u_{\infty,k}$ is a
solution of (5.1) for $H = H(\eta_1)$, possibly with a finite
number of bubbles attached. We also have
$$
\lim_i E_{(H^{\rho(\tau)}, J^{\rho(\tau)})}(u_i) = \sum_k
E_{(H^{\rho(\tau)}, J^{\rho(\tau)})}(u_{\infty,k})
$$
On the other hand the left half of (7.9) implies that at least one
of $u_{\infty,k}$ is not trivial, i.e., not stationary.

Now we consider three cases separately: the first is the one where
$\eta_1 \in \SS ing(\HH)$ and the second where $\eta_1 \in \NN
t(\HH,j)$ and the rest where $$\eta_1 \in [0,1] \backslash (\SS
ing(\HH)\cup \NN t(\HH,j)).$$ When $\eta_1 \in \SS ing$, we must
have
$$
\lim_i E_{(H^{\rho(\tau)}, J^{\rho(\tau)})}(u_i) \geq
\min\{A_{(H(\eta_1),J(\eta_1))}^{sing}, A(\omega;j)\} \geq
A_{(\HH,j)}^0
$$
which gives rise to a contradiction to (7.9) when $i$ is
sufficiently large. On the other hand, if $\eta_1 \in \NN
t(\HH,j)$, the cusp curve must contain a component $u_\infty$ that
has Index 0 and is non-constant. Again the right hand side
inequality of (7.9) prevents this from happening. Finally when
$\eta_1 \in [0,1] \backslash (\SS ing(\HH) \cup \NN t(\HH,j))$,
the index condition $\text{Index } u_i = 0$ and the transversality
condition imply that all the components $u_{\infty,k}$ must be
constant which again contradicts to the left hand side inequality
of (7.9) if $i$ is sufficiently large. This finishes the proof of
the handle sliding lemma.
\end{proof}

\section*{\bf \S 8. Parametric stability of tightness of Novikov Floer cycles}

This is the key section which will involve all the results we
proved in section 3-7, especially Theorem 3.7, Theorem 4.6,
Proposition 5.7, Theorem 6.7 and Proposition 7.4.

Let
$$
a = \sum a_A q^{-A}, \quad a_A \in H^*(M) \leqno (8.1)
$$
be a non-zero quantum cohomology class. We denote by $\Gamma(a)
\subset \Gamma$ the set of $A$'s for which the coefficient $a_A$
is non-zero. By the definition of the Novikov ring, we can
enumerate $\Gamma(a)$ so that
$$
- \lambda_1 < - \lambda_2 < \cdots  <- \lambda_j < \cdots
$$
where $\lambda_j = \omega(A_j)$. We call the first term $a_1
q^{-A_1}$ the {\it leading order term} of the quantum cohomology
class $a$ and denote by $Ld(a)$.

We recall from \cite{Oh5} that for a given quantum cohomology class
$a\neq 0$, we define the mini-max value of $\AA_H$
$$
\rho(H;a) = \inf_{\alpha}\{\lambda_H(\alpha) \mid \alpha \in \ker
\part, \, \text{ with }\, [\alpha] = a^\flat \}
$$
for a nondegenerate Hamiltonian $H$ for which the Floer homology
$HF_*(H,J)$ is defined for a generic choice of $J$. As we pointed
out in \cite{Oh5}, the number $\rho(H;a)$ is independent of the
choice of $J$.

The following notion of {\it tight Floer cycles} introduced in
\cite{Oh5} is important in the proof of Theorem II.
\medskip

{\noindent\bf Definition 8.1.} Let $(H,J)$ be a Floer regular pair
so that the Floer complex $(CF_*(H), \part_{(H,J)})$ is defined.
Let $\alpha$ be a Floer cycle of $H$ and $a \in QH^*(M)$ be the
corresponding quantum cohomology class with $[\alpha] = a^\flat$.
We call the Floer cycle $\alpha$ of $(H,J)$ {\it tight} if it
realizes the mini-max value, i.e.,
$$
\lambda_H(\alpha) = \rho(H;a).
$$
We call a critical value $\lambda$ of $\AA_H$ a {\it homologically
essential critical value} of $\AA_H$, if there exists $J$ such
that $(H,J)$ is Floer regular, and $\lambda = \lambda_H(\alpha)$
for a tight Floer cycle of $(H,J)$.
\medskip

Although the homologically essentialness of a critical point
$[z,w]$ of $\AA_H$ depends on the choice of $J$, the following
proposition proves that the homologically essentialness of a
critical value is independent of the choice of $H$-regular $J$'s.
\medskip

{\noindent\bf Proposition 8.2.} {\it Let $J, \, J'$ be two $H$-regular
one-periodic family of almost complex structures. Suppose that a
Floer cycle $\alpha$ of $(H,J)$ is tight and assume that $j$ is a
homotopy connecting $J, \, J'$ that is $\HH$-regular for the
constant homotopy $\HH\equiv H$. Then the transferred cycle
$$
h_{(\HH,j;\rho)}(\alpha)
$$
defines a tight Floer cycle for the pair $(H,J')$.
}
\begin{proof} Since $[h_{(H,j;\rho)}(\alpha)] = a^\flat$, we have
$$
\lambda_H(h_{(H,j;\rho)}(\alpha)) \geq \rho(H;a) =
\lambda_H(\alpha). \leqno (8.2)
$$
Next let $[z',w']$ be a peak of $h_{(H,j;\rho)}(\alpha)$. By the
definition of $h_{(H,j;\rho)}(\alpha)$, there is $[z,w] \in
\alpha$ such that $\MM((H,j;\rho);[z,w],[z',w']) \neq \emptyset$.
Then (5.8) applied to the constant homotopy $\HH \equiv H$ implies
$$
\AA_H([z',w']) \leq \AA_H([z,w])
$$
which implies
$$
\lambda_H(h_{(H,j;\rho)}(\alpha)) \leq \lambda_H(\alpha). \leqno (8.3)
$$
Combining (8.2) and (8.3), we have proved
$$
\lambda_H(h_{(H,j;\rho)}(\alpha)) = \rho(H;a)
$$
which finishes the proof.
\end{proof}

Now we fix a homotopy $\HH \in \PP^{Cerf}(\HH_m(M);\e f, H)$ and
$j \in \PP^{tran}(j_\omega;\HH)$, satisfying the properties
described in Theorem 4.6,
$$
\HH=\{H(\eta)\}_{0 \leq \eta \leq 1}, \quad j = \{J(\eta)\}_{0\leq
\eta \leq 1}
$$
such that
$$
H(0) = H_\alpha, \quad H(1) = H_\beta
$$
are nondegenerate. In particular, by the choice of $j$, the Floer
homology $HF_*(H(\eta),J(\eta))$ is defined for any $\eta \in
I(\HH,j)$.

Let $\eta_0$ be a Floer point at which there exists a tight Floer
cycle $\alpha_0 \in CF_*(H(\eta_0))$. We fix a homotopy
$$
j^{\eta_0} \in \PP^{sub}(j_\omega;\HH;\eta_0); \quad j^{\eta_0} =
\{ J^{\eta_0}_t \}_{0 \leq t \leq 1} \leqno (8.4)
$$
that satisfies the properties of Theorem 4.6 for the pair
$(\HH,j^{\eta_0})$ and also satisfies
$$
J^{\eta_0}(\eta_0) = J(\eta_0). \leqno (8.5)
$$
The main result of this section is to show that the tight Floer
cycle $\alpha_0$ of $(H(\eta_0),J(\eta_0))$, {\it after perturbing
$j$ slightly to $j^{\eta_0}$ as above}, is {\it parametrically
stable} in that there exists $\delta
> 0$ such that for any Floer point $\eta$ with $|\eta - \eta_0| <
\delta$, $\rho(H(\eta);a)$ is the level of a tight Floer cylce of
the pair $(H(\eta),J^{\eta_0}(\eta))$. Note that once we know the
latter fact, Proposition 8.2 implies that $\rho(H;a)$ is indeed
the level of a tight Floer cycle of $(H(\eta),J(\eta))$ for the
original homotopy $(\HH,j)$.

We fix a cut-off function $\rho_0: \R \to [0,1]$ and fix one
$j^{\eta_0}$ as in $(8.4)$ and satisfying $(8.5)$. Then consider the
corresponding chain map
$$
h_{(\eta_0\eta;\rho_0)}: =
h_{(\HH_{\eta_0\eta},j^{\eta_0}_{\eta_0\eta};\rho_0)} :
CF_*(H(\eta_0)) \to CF_*(H(\eta)). \leqno (8.6)
$$
\medskip

{\noindent\bf Theorem 8.3.} {\it Let $(\HH,j)$ be as above and  $0
\neq a \in QH^*(M)$. Suppose that $\eta_0 \in [0,1]$ is a Floer
point at which $a$ carries a tight Floer cycle for $H(\eta_0)$. Then
there exists $\delta_4 = \delta_4(\HH,j^{\eta_0}, a;\eta_0) > 0$
such that on each of the semi-intervals $[\eta_0, \eta_0 +
\delta_4)$ or $(\eta_0 - \delta_4,\eta_0]$ $a$ carries tight cycles
given by the transferred cycles
$h_{(\eta_0\eta;\rho_0)}(\alpha_\pm)$ of some tight cycles
$\alpha_\pm$ at $\eta_0$ respectively.}
\medskip

We would like to remark that when $\eta_0 \not \in
\CC^{nd}ross(\HH)$, we may take $\alpha_- = \alpha_+$. On the other
hand, if $\eta_0 \in \CC^{nd}ross(\HH)$, the two cycles $\alpha_-$
and $\alpha_+$ could be different.

The rest of the section will be occupied by the proof of this
theorem.

We note that since $\eta_0$ is a Floer point, there are only a
finite number of periodic orbits of $H(\eta_0)$, and can apply
Theorem 6.7 for {\it all} Floer chains of $(H(\eta_0),J(\eta_0))$.

By Theorem 4.6, the chain map (8.6) is defined at any point $\eta
\in I(\HH, j^{\eta_0};\eta_0)$, a dense subset of $[0,1]$. {\it We
emphasize that we need to choose $j^{\eta_0}$ depending on the
point $\eta_0$ to ensure the properties stated in Theorem 4.6. }

We now prove the following key proposition. We would like to
emphasize that this kind of continuity statement in the levels of
cycles, not in the levels of the corresponding homologies, does
not hold in general, and even if it holds so, proving such a
continuity statement is a highly non-trivial matter due to the
phenomenon of {\it cancellation of critical points}. This is the
reason why the structure theorem, Theorem 6.7, of Floer cycles
proven in section 6 is so crucial in our proof.
\medskip

{\noindent\bf Proposition 8.4.} {\it Let $(\HH,j)$, $\alpha_0$, $\eta_0$ and
$j^{\eta_0}$ be as above. Then there exists some $\delta_5 =
\delta_5(\HH,j^{\eta_0},\eta_0)
> 0$ such that the assignment
$$
\eta \mapsto \lambda_{H(\eta)}(h_{(\eta_0\eta;\rho_0)}(\alpha_0))
$$
is continuous on $(\eta_0 - \delta_5, \eta_0 + \delta_5) \cap
I(\HH, j^{\eta_0};\eta_0)$.
}

\begin{proof} Let $\delta > 0$ which is to be determined later and
consider the function
$$
\mu(\eta) = \lambda_{H(\eta)}(h_{(\eta_0\eta;\rho_0)}(\alpha_0)).
$$
We remind the readers that  by definition  all the generators of a
Floer cycle representing the class dual to a given $a \in QH^*(M)$
have the same Conley-Zehnder indices.

For the proof of the proposition, we will follow the scheme used in
the appendix of \cite{Oh2}. We first note that
$$
\mu(\eta)=\lambda_{H(\eta)}(h_{(\eta_0\eta;\rho_0)}(\alpha_0))
\geq \rho(H(\eta);a)
$$
by the definition of $\rho(H(\eta);a)$ since we have
$[h_{\eta_0\eta}(\alpha_0)] = a^\flat$. Hence $\mu(\eta)$ is
finite and well-defined.

Once the finiteness of $\mu(\eta)$ is shown, the proposition will
be an immediate consequence of the following lemma.
\medskip

{\noindent\bf Lemma 8.5.} {\it There exists $\delta_5 =
\delta_5(\HH,j^{\eta_0};\eta_0)
> 0$ for which we have the following inequality
$$
-E(\HH_{\eta\eta'}) \leq \mu(\eta') - \mu(\eta) \leq
E(\HH_{\eta\eta'}) \leqno (8.7)
$$
for any Floer points $\eta, \, \eta' \in (\eta_0 - \delta_5,
\eta_0 - \delta_5)$. In particular $\mu$ is continuous at
$\eta_0$.
}

\medskip

Therefore we will prove this lemma in the rest of the proof of
Proposition 8.4. The proof of this lemma is quite long and
intricate, and various arguments used in the proof touch the heart
of the chain level Floer theory.

We first consider the case $\eta_0 < \eta < \eta'$ and provide
complete details of the proof of (8.7) in this case. We will
briefly mention the proof of (8.7) for the other cases in the end
of the proof.

We compare the homotopy $\HH_{\eta_0\eta'}$ with the glued
homotopy $\HH_{\eta\eta'}\#\HH_{\eta_0\eta}$. We recall the
general gluing identity from (5.24)
$$
h_{\HH_{\eta_0\eta}\#_{(\rho;R)} \HH_{\eta\eta'}} =
h_{(\eta\eta';\rho_1)}\circ h_{(\eta_0\eta;\rho_0)} \leqno (8.8)
$$
for $\rho =(\rho_0,\rho_1)$ and for a sufficiently large $R>0$. We
also assume $\rho_1$ is monotone. Then, if
$[z'_{\eta'},w'_{\eta'}]$ is a peak of the cycle
$h_{\HH_{\eta_0\eta}\#_{(\rho;R)} \HH_{\eta\eta'}}(\alpha_0)$, by
the definition of the chain map $h_{(\eta\eta';\rho_1)}$, there
must exist some $[z',w'] \in h_{\eta_0\eta}(\alpha_0)$ for which
$$
\MM\Big((\HH_{\eta\eta'},j_{\eta\eta'};\rho_1);
[z',w'],[z_{\eta'}',w_{\eta'}']\Big) \neq \emptyset \leqno (8.9)
$$
holds. Now we state the following easy general lemma.
\medskip

{\noindent\bf Lemma 8.6.} {\it Let $F$ and $K$ be two nondegenerate
Hamiltonians and $\HH$ be a homotopy from $F$ to $K$ and $j$ be
given such that $(\HH,j)$ is Floer regular. Let $\rho$ be a given
monotone cut-off function and let $h_{\HH, j;\rho}$ be the
corresponding Floer homotopy map. Then we have
$$
\lambda_{K}(h_{(\HH,j;\rho)}(\alpha)) \leq \lambda_{F}(\alpha) +
E^-(\HH)
$$
for any Floer chain $\alpha$ of $F$.
}
\begin{proof} Let $[z',w']$ be a peak of the cycle
$h_{(\HH,j;\rho)}(\alpha)$. By the definition of
$h_{(\HH,j;\rho)}(\alpha)$, there is a generator $[z,w] \in
\alpha$ such that
$$
\MM\Big((\HH,j;\rho);[z,w],[z',w']\Big) \neq \emptyset.
$$
Then (5.10) implies
$$
\AA_K([z',w']) \leq \AA_F([z,w]) + E^-(\HH_{\eta\eta'}). \leqno (8.10)
$$
Since $[z',w']$ is a peak of the cycle $h_{(\HH;\rho)} (\alpha)$,
i.e., $\AA_K([z',w']) = \lambda_K(h_{(\HH,j;\rho)}(\alpha))$ and
$\AA_F([z,w]) \leq \lambda_{F}(\alpha)$, the lemma follows from
(8.10).
\end{proof}

Going back to the proof, Lemma 5.5 \& 5.6, (8.8) and Lemma 8.6
imply
$$
\begin{aligned} \lambda_{H(\eta')}(h_{\HH_{\eta_0\eta}\#_{(\rho;R)}
\HH_{\eta\eta'}}(\alpha_0)) & \leq
\lambda_{H(\eta)}(h_{(\eta_0\eta;\rho_0)}(\alpha_0) )+
E^-(\HH_{\eta\eta'}) \\
& = \mu(\eta) +  E^-(\HH_{\eta\eta'}).
\end{aligned}\leqno (8.11)
$$
We next deform the homotopy $\HH_{\eta_0\eta}\#_{(\rho;R)}
\HH_{\eta'\eta}$ to the $\rho_0$-elongated homotopy of
$\HH_{\eta_0\eta'}$ for the given fixed choice of $\rho_0$ that we
used in the definition of $\mu$ in the beginning. On the other hand,
since we assume that $\eta_0 < \eta < \eta'$, {\it
$\HH_{\eta_0\eta}\#_{(\rho;R)} \HH_{\eta\eta'}$ is also an
elongation of the sub-homotopy $\HH_{\eta_0\eta'}$} (by a monotone
cut-off function). More explicitly, if we denote by $\rho_{elng}$
the corresponding cut-off function, $\rho_{elng}$ is given by the
formula (5.22).

Now we consider the homotopy
$$
(\overline j_{\eta_0\eta'}, \overline
\HH_{\eta_0\eta'};\overline\rho)
$$
between the $\rho_0$-elongation of $(\HH,j^{\eta_0})$ and the
$\rho_{elng}$-elongation of $(\HH,j^{\eta_0})$ connected by the
homotopy $\overline\rho = \{\rho(\kappa)\}_{0 \leq \kappa \leq 1}$
of cut-off functions
$$
\overline\rho=\{\rho(\kappa)\}_{0 \leq \kappa \leq 1}; \quad
\rho(0) = \rho_0, \, \rho(1) = \rho_{elng}.
$$
Since $\alpha_0$ is a cycle, we derive, from Proposition 5.7,
$$
h_{(\eta_0\eta';\rho_0)}(\alpha_0) =
h_{\HH_{\eta_0\eta}\#_{(\rho;R)} \HH_{\eta\eta'}}(\alpha_0) +
\part_{(J(\eta'),H(\eta'))}\circ
H_{(\overline j_{\eta_0\eta'}, \overline
\HH_{\eta_0\eta'};\overline\rho)} (\alpha_0) \leqno (8.12)
$$
with the inequality
$$
\lambda_{H(\eta')}(H_{(\overline j_{\eta_0\eta'}, \overline
\HH_{\eta_0\eta'};\overline\rho)} (\alpha_0)) \leq
\lambda_{H(\eta_0)}(\alpha_0) + E(\HH_{\eta_0\eta'}). \leqno (8.13)
$$
Now using the non-Archimedean triangle inequality
$$
\lambda_{H(\eta')}(\alpha + \beta) \leq \max\{ \lambda_{H(\eta')}
(\alpha), \lambda_{H(\eta')}(\beta)\},
$$
we estimate the level
$$
\lambda_{H(\eta')}\Big(h_{\HH_{\eta_0\eta}\#_{(\rho;R)}
\HH_{\eta\eta'}}(\alpha_0) +
\part_{J(\eta'),H(\eta')}\circ
H_{(\overline j_{\eta_0\eta'}, \overline
\HH_{\eta_0\eta'};\overline\rho)} (\alpha_0)\Big).\leqno (8.14)
$$
For the first term of (8.14), we have (8.11). For the second term,
we recall Theorem 6.7 and Proposition 7.4 which imply
$$
\begin{aligned}
&\ \lambda_{H(\eta')}\Big(
\part_{J(\eta'),H(\eta')}\circ
H_{(\overline j_{\eta_0\eta'}, \overline
\HH_{\eta_0\eta'};\overline\rho)} (\alpha_0)\Big) \\
 \leq &\
\lambda_{H(\eta')}\Big(H_{(\overline j_{\eta_0\eta'}, \overline
\HH_{\eta_0\eta'};\overline\rho)} (\alpha_0)\Big) -
A^1_{(\HH,j^{\eta_0};\eta_0)}
\end{aligned}\leqno (8.15)
$$
for all Floer points $\eta'$ satisfying
$$
|\eta'-\eta_0| < \delta \quad \text{and} \quad d(\eta',\SS
ing(\HH)) > \delta
$$
for some $\delta = \delta(\HH,j,\eta_0)$. Therefore combining
(8.11)-(8.15), we derive
$$
\begin{aligned}\mu(\eta') & = \lambda_{H(\eta')}(h_{(\eta_0\eta';
\rho_0)}(\alpha_0))
\quad \\
& \leq \max\Big\{\mu(\eta) + E^-(\HH_{\eta\eta'}), \,
\lambda_{H(\eta_0)}(\alpha_0) + E(\HH_{\eta_0\eta'}) -
A^1_{(\HH,j^{\eta_0};\eta_0)} \Big \}.
\end{aligned}\leqno (8.16)
$$
We will study under what conditions, we have
$$
\mu(\eta) + E^-(\HH_{\eta\eta'}) \geq
\lambda_{H(\eta_0)}(\alpha_0) + E(\HH_{\eta_0\eta'}) -
A^1_{(\HH,j^{\eta_0};\eta_0)}
$$
and in turn $\mu(\eta') \leq \mu(\eta) + E^-(\HH_{\eta\eta'})$.

We apply the homotopy $(\HH_{\eta_0\eta'}^{-1})_{\rho_0}$ to (8.8)
and obtain
$$
h_{(\HH_{\eta_0\eta'}^{-1};\rho_0)}\circ
h_{\HH_{\eta_0\eta}\#_{(\rho;R)} \HH_{\eta\eta'}} =
h_{(\HH_{\eta_0\eta'}^{-1};\rho_0)} \circ
h_{(\eta\eta';\rho_1)}\circ h_{(\eta_0\eta;\rho_0)} \leqno (8.17)
$$
where $\rho = (\rho_0,\rho_1)$. It follows from (8.17) that the
cycle $h_{(\HH_{\eta_0\eta'}^{-1};\rho_0)} \circ
h_{(\eta\eta';\rho_1)}\circ h_{(\eta_0\eta;\rho_0)}(\alpha_0)$ is
homologous to $\alpha_0$ since
$h_{(\HH_{\eta_0\eta'}^{-1};\rho_0)}\circ
h_{\HH_{\eta_0\eta}\#_{(\rho;R)} \HH_{\eta\eta'}}$ is chain
homotopic to the identity. Therefore, by the tightness assumption
on $\alpha_0$, we derive
$$
\lambda_{H(\eta_0)}(h_{(\HH_{\eta_0\eta'}^{-1};\rho_0)} \circ
h_{(\eta\eta';\rho_1)}\circ h_{(\eta_0\eta;\rho_0)}(\alpha_0))
\geq \lambda_{H(\eta_0)}(\alpha_0). \leqno (8.18)
$$
Applying Lemma 8.6 to $h_{(\HH_{\eta_0\eta'}^{-1};\rho_0)}$ and
the cycle $h_{(\eta\eta';\rho_1)}\circ
h_{(\eta_0\eta;\rho_0)}(\alpha_0)$, we obtain
$$
\begin{aligned} \lambda_{H(\eta_0)}(h_{(\HH_{\eta_0\eta'}^{-1};\rho_0)} &
\circ h_{(\eta\eta';\rho_1)}\circ
h_{(\eta_0\eta;\rho_0)}(\alpha_0)) \\
& \leq  \lambda_{H(\eta')}( h_{(\eta\eta';\rho_1)}\circ
h_{(\eta_0\eta;\rho_0)}(\alpha_0)) + E^-(\HH_{\eta_0\eta'}^{-1})\\
& = \lambda_{H(\eta')}( h_{(\eta\eta';\rho_1)}\circ
h_{(\eta_0\eta;\rho_0)}(\alpha_0)) + E^+(\HH_{\eta_0\eta'}).
\end{aligned}\leqno (8.19)
$$
Then (8.17)-(8.19) give rise to
$$
\begin{aligned} \lambda_{H(\eta_0)}(h_{(\eta\eta';\rho_1)}\circ
h_{(\eta_0\eta;\rho_0)}(\alpha_0)) & \geq
\lambda_{H(\eta_0)}(\alpha_0) - E^+(\HH_{\eta_0\eta'})\end{aligned}  \leqno (8.20)
$$
$$\begin{aligned}
& \geq \lambda_{H(\eta)}(h_{(\eta_0\eta;\rho_0)}(\alpha_0)) -
E^-(\HH_{\eta_0\eta}) - E^+(\HH_{\eta_0\eta'})
\end{aligned}\leqno (8.21)
$$
where we applied Lemma 8.6 for the latter inequality. Now after
rewriting (8.12) as
$$
h_{\HH_{\eta_0\eta}\#_{(\rho;R)} \HH_{\eta\eta'}}(\alpha_0) =
h_{(\eta_0\eta';\rho_0)}(\alpha_0) -
\part_{(J(\eta'),H(\eta'))}\circ
H_{(\overline j_{\eta_0\eta'}, \overline
\HH_{\eta_0\eta'};\overline\rho)} (\alpha_0),
$$
we derive,  from this, (8.8), (8.14) and (8.15), that we have
$$
\begin{aligned} \lambda_{H(\eta')}(h_{(\eta_0\eta;\rho_0)} & \circ
h_{(\eta\eta';\rho_1)}(\alpha_0)) =
\lambda_{H(\eta')}(h_{\HH_{\eta_0\eta}\#_{(\rho;R)}
\HH_{\eta\eta'}}(\alpha_0)) \\
\leq \max & \Big\{\lambda_{H(\eta')}(h_{(\eta_0\eta';\rho_0)}(\alpha_0)),\\
& \quad \lambda_{H(\eta_0)}(\alpha_0) + E(\HH_{\eta_0\eta'}) -
A^1_{(\HH,j^{\eta_0};\eta_0)}\Big\}
\end{aligned}
\leqno (8.22)
$$
as in (8.16). On the other hand, we derive, from (8.20) and
(8.22),
$$
\begin{aligned}
\max\Big\{\lambda_{H(\eta')}(h_{(\eta_0\eta';\rho_0)}(\alpha_0)),
\lambda_{H(\eta_0)}(\alpha_0) & +
E(\HH_{\eta_0\eta'}) - A^1_{(\HH,j^{\eta_0};\eta_0)}\Big\} \\
& \geq \lambda_{H(\eta_0)}(\alpha_0) - E^+(\HH_{\eta_0\eta'}).
\end{aligned}
\leqno (8.23)
$$
Now we choose $\delta_5= \delta_5(\HH,j^{\eta_0},\eta_0) > 0$ so
that
$$
E(\HH_{\eta_0\eta'}) - A^1_{(\HH,j^{\eta_0};\eta_0)} < -
E^+(\HH_{\eta_0\eta'})
$$
i.e.,
$$
 E(\HH_{\eta_0\eta'})+ E^+(\HH_{\eta_0\eta'})
< A^1_{(\HH,j^{\eta_0};\eta_0)} \leqno (8.24)
$$
holds for any $\eta, \, \eta' \in (\eta_0 - \delta_5, \eta_0 -
\delta_5)$. We would like to emphasize that {\it we can choose
$\delta_5$ so that it satisfies $(8.24)$ and depends only on
$(\HH,j,\eta_0)$}. With this choice of $\delta_5$ made, we have
$$
\lambda_{H(\eta_0)}(\alpha_0) + E(\HH_{\eta_0\eta'}) -
A^1_{(\HH,j^{\eta_0};\eta_0)} < \lambda_{H(\eta_0)}(\alpha_0) -
E^+(\HH_{\eta_0\eta'}). \leqno (8.25)
$$
Then (8.23) and (8.25) imply
$$
\lambda_{H(\eta')}(h_{(\eta_0\eta';\rho_0)}(\alpha_0)) \geq
\lambda_{H(\eta_0)}(\alpha_0) - E^+(\HH_{\eta_0\eta'}) \leqno (8.26)
$$
Now combining (8.16), (8.25) and (8.26), we first obtain
$$
\begin{aligned} \max\Big\{\mu(\eta) + E^-(\HH_{\eta\eta'}), \, &
\lambda_{H(\eta_0)}(\alpha_0) +
E(\HH_{\eta_0\eta'}) - A^1_{(\HH,j^{\eta_0};\eta_0)} \Big \} \\
& = \mu(\eta) + E^-(\HH_{\eta\eta'})
\end{aligned}
$$
which in turn implies
$$
\mu(\eta') \leq \mu(\eta) + E^-(\HH_{\eta\eta'}). \leqno (8.27)
$$

Next we compare the homotopy $\HH_{\eta_0\eta}$ and
$\HH_{\eta\eta'}^{-1}\# \HH_{\eta_0\eta'}$. Recall $\eta' > \eta$
and $\HH_{\eta'\eta}= \HH^{-1}_{\eta\eta'}$ by Definition 3.8. Again
we will have the identity
$$
h_{(\HH^{-1})_{\eta\eta'}\#_{(\rho;R)}\HH_{\eta_0\eta'}} =
h_{\eta'\eta;\rho_1} \circ h_{\eta_0\eta';\rho_0}
$$
for a sufficiently large $R$.  Now (8.11) is replaced by
$$
\lambda_{H(\eta)}(h_{(\HH^{-1})_{\eta\eta'}\#_{(\rho;R)}
\HH_{\eta_0\eta'}}(\alpha_0)) \leq \mu(\eta') +
E^-(\HH^{-1}_{\eta\eta'}) = \mu(\eta') + E^+(\HH_{\eta\eta'}).
$$
Here we would like to remark that
$(\HH^{-1})_{\eta\eta'}\#_{(\rho;R)}\HH_{\eta_0\eta'}$ is still an
elongation of $\HH_{\eta_0\eta}$ but by a non-monotone cut off
function.

By repeating the arguments above with $\HH_{\eta\eta'}$ replaced
by $\HH_{\eta\eta'}^{-1}$, we obtain
$$
\mu(\eta) \leq \mu(\eta') + E^+(\HH_{\eta\eta'}). \leqno (8.28)
$$
Combining (8.27) and (8.28), we have obtained
$$
-E^+(\HH_{\eta\eta'}) \leq  \mu(\eta')-\mu(\eta) \leq
E^-(\HH_{\eta\eta'}).
$$
This is indeed an inequality stronger than (8.7) and in particular
proves (8.7). The cases other than that of $\eta_0 < \eta < \eta'$
can be handled by the same arguments if we replace the
sub-homotopy $\HH_{ss'}$ by $(\HH^{-1})_{s's}$ every time the
reverse inequality $s > s'$ appears in the proof. This explains
appearance of $E$ in general instead of $E^-$ or $E^+$ in (8.7).
We leave the details to the readers. This finishes the proof of
Lemma 8.5 and hence Proposition 8.4.
\end{proof}

\medskip

{\noindent \bf Remark 8.7.} In fact, an examination of the above
proof combined with the discussion above Proposition 3.8 shows that
the map $\eta \mapsto
\lambda_{H(\eta)}(h_{(\eta_0\eta;\rho_0)}(\alpha_0))$ is
differentiable such that its derivative depends only on the periodic
orbit but not on its liftings and is given by
$$
-\int_0^1 \frac{\partial H(\eta)}{\partial \eta}(t, z_\eta(t))\, dt
$$
where $z_\eta$ is a periodic orbit with $[z_\eta,w_\eta]$ is a peak
of the cycle $h_{(\eta_0\eta;\rho_0)}(\alpha_0)$.
\medskip

Next we note that there is a canonical one-one correspondence
$$
\operatorname{Per}(H(\eta_0)) \leftrightarrow \operatorname{Per}(H(\eta))
$$
for any $\eta \in (\eta_0 - \delta,\eta_0 + \delta)$ for some
$\delta = \delta(\HH,\eta_0) > 0$. We denote the corresponding
family by $z_k(\eta)$ for $k = 1, \cdots, \#(\operatorname{Per}(H(\eta_0))$.
This map in turn induces a one-one correspondence
$$
\operatorname{Crit}\AA_{H(\eta_0)} \leftrightarrow \operatorname{Crit}\AA_{H(\eta_0)}
$$
given by
$$
[z_k(\eta_0),w_k(\eta_0)] \leftrightarrow [z_k(\eta),w_k(\eta_0)
\# u_{z_k(\eta_0)z_k(\eta)}^{can}].
$$
For the simplicity of notation, we denote
$u_{\eta_0\eta,k}^{can} = u_{z_k(\eta_0)z_k(\eta)}^{can}$.
By construction of $h_{(\eta_0\eta;\rho_0)}$, we then have
$$
h_{(\eta_0\eta;\rho_0)}([z_k(\eta_0),w_k(\eta_0)]) =
[z_k(\eta),w_k(\eta_0) \# u_{\eta_0\eta,k}^{can}] \leqno (8.29)
$$
modulo any filtration order as want, by choosing $\delta$ smaller if
necessary. We recall the identity $h_{(\eta\eta;\rho_0)} \equiv id$
in the chain level. We will always assume this holds for the rest of
the proof.

We start with the following result which is of independent interest.
A similar statement has been an important ingredient of Usher's
algebraic proof of Theorem II \cite{usher}. It was also proved in
\cite{fooo06} in the context of Lagrangian intersection Floer
cohomology by a purely algebraic way, but with
$\Lambda_{0,nov}^{(0)}$ the (positive) universal Novikov ring. (See
Proposition 26.9 \cite{fooo06}.) Since we use the Novikov ring
$\Lambda_\omega$, which does not have this positivity property, we
cannot deduce this result from Proposition 26.9  \cite{fooo06}. Here
we instead give a simple geometric proof based on the basic fact,
Lemma 8.9 below, on the Floer boundary operator for a small Morse
function.
\medskip

{\noindent \bf Proposition 8.8.} {\it There exists a constant $c >
0$ depending only on $\HH$ but independent of $\eta_0$ or $\lambda$
such that
$$
\partial\left(CF(H(\eta_0)\right) \cap CF^\lambda(H(\eta_0))
\subset \partial\left(CF^{\lambda + c}(H(\eta_0))\right)
$$
for all $\lambda \in \R$. }
\begin{proof}
Denote $H = H(\eta_0)$ and $\partial_{(H,J)}=\partial_H$ in this
proof for the simplicity of notations. We first recall the following
well-known result for a small Morse function $\e f$. See \cite{Fl4},
\cite{Oh5} for its proof.
\medskip

{\noindent \bf Lemma 8.9.} {\it Let $f$ be a Morse function and $\e>
0$ be a sufficiently small so that so that $\partial_{\e f} =
\partial^{Morse}(-\e f) \otimes \Lambda_\omega$. Let $c_0 = \e
\max|f|$. In particular, we have
$$
\partial (CF(\e f)) \cap CF^\lambda(\e f)
\subset \partial (CF_k^{\lambda + c_0}(\e f)).
$$
}

We now connect the given Cerf homotopy $\HH$ by another {\it fixed}
homotopy from $\e f$ and $H^0$ and denote the join homotopy again by
$\HH$ which now connects $\e f$ to $H$. Let $\gamma \in
\partial (CF(H)) \cap CF^\lambda(H)$. We transfer
$\gamma$ to $CF(\e f)$ and consider $h_{\eta_0 0}(\gamma)$. Note
$$
\lambda_{\e f}(h_{\eta_0 0}(\gamma)) \leq \lambda_H(\gamma) +
E^+(\HH) \leq \lambda + E^+(\HH).
$$
Therefore it follows from Lemma A.4  that there exists a chain
$\beta_0 \in CF(\e f)$ such that
$$
h_{\eta_0 0}(\gamma) = \partial_{\e f}(\beta_0), \quad \lambda_{\e
f}(\beta_0) \leq \lambda + E^+(\HH) + \e \max|f|.
$$
We estimate
$$
\lambda_H(h_{0\eta_0}(\beta_0))\leq \lambda + E^+(\HH) + \e \max|f|
+ E^-(\HH) = \lambda + E(\HH) + \e \max|f|. \leqno (8.30)
$$
We compute $\partial_H(h_{0\eta_0}(\beta_0))$
$$
\partial_H(h_{0\eta_0}(\beta_0)) = h_{0\eta_0}(\partial_{\e f}(\beta_0))
= h_{0\eta_0}(h_{\eta_0 0}(\gamma)).
$$
By the chain homotopy formula, we obtain
$$
h_{0 \eta_0}\circ h_{\eta_0 0}(\gamma) = \gamma +
\partial_H(\overline H_{\HH_{\eta_00}\# \HH_{0\eta_0}}(\gamma))
$$
and hence
$$
\partial_H(h_{0\eta_0}(\beta_0) - \overline H_{\HH_{\eta_00}\# \HH_{0\eta_0}}(\gamma)) =
\gamma.
$$
We consider $\beta = h_{0\eta_0}(\beta_0) - \overline
H_{\HH_{\eta_00}}(\gamma)$ and estimate its level
$$
\lambda_H(\beta) \leq \max\{\lambda_H(h_{0\eta_0}(\beta_0)),
\lambda_H(\overline H_{\HH_{\eta_00}\# \HH_{0\eta_0}}(\gamma))\}.
$$
But we have the estimates
$$
\lambda_H(h_{0\eta_0}(\beta_0)) \leq \lambda + E(\HH) + \e \max|f|
$$
from (8.30) and
$$
\lambda_H(\overline H_{\HH_{\eta_00}\# \HH_{0\eta_0}}(\gamma)) \leq
\lambda_H(\gamma) + E(\HH) \leq \lambda + E(\HH)
$$
from (5.25). Hence by taking $c = E(\HH) + \e \max|f|$, we have
finished the proof.
\end{proof}

Now we are ready to give the proof of Theorem 8.3.

\begin{proof}[Proof of Theorem 8.3]
Fix a positive constant $\delta < \min \lambda_1 /2$. Proposition
8.8 implies that any element from the submodule
$\partial\left(CF(H(\eta_0)\right) \cap CF^\lambda(H(\eta_0))$ with
$\lambda =\rho(H(\eta_0);a)+\delta$ can be represented by
$\partial_{H(\eta_0)}(\beta)$ with
$$
\lambda_{H(\eta_0)}(\beta) \leq \rho(H(\eta_0);a) + \delta+ c.
$$
Now we choose a tight cycle $\alpha_0$ and its decomposition
$$
\alpha_0 = \operatorname{peak}(\alpha_0) + \widetilde \alpha_0.
$$
We denote the gap of the operator $\partial_{(H(\eta_0),J(\eta_0))}$
by $\lambda_1 = \lambda_1(\eta_0): = A_{(H(\eta_0),J(\eta_0))}$
given in Proposition 6.2.

We then consider the transferred cycle $h_{\eta_0\eta}(\alpha_0):=
\alpha(\eta)$ for $\eta$ sufficiently close to $\eta_0$. We will
show the following two alternatives : there exists $\delta_4 =
\delta_4(\HH,j,a;\eta_0)$ such that
\begin{enumerate}
\item[(1)] either
$h_{\eta_0\eta}(\alpha_0)$ is tight for $H(\eta)$ for all $\eta \in
(\eta_0 - \delta_4, \eta_0 + \delta_4)$,
\item[(2)] or $h_{\eta_0\eta}(\alpha_0)$ is not tight
and there exist tight cycles $\alpha_\pm$ for $H(\eta_0)$ such that
$h_{\eta_0\eta}(\alpha_\pm)$ are tight on $(\eta_0 -
\delta_4,\eta_0]$ and on $[\eta_0,\eta_0 + \delta_4)$ respectively.
\end{enumerate}
In the second case, $\alpha_-$ and $\alpha_+$ could be different
cycles.

Now suppose that the cycles $h_{\eta_0\eta}(\alpha_0)$ is not
tight for $H(\eta)$ at least in one direction, say, for $\eta <
\eta_0$. We will try to find another tight cycle $\alpha_-$ at
$\eta_0$ for which the second alternative above holds.

Since $h_{\eta_0\eta}(\alpha_0)$ is not tight, there exists a chain
$\beta(\eta)$ such that
$$
\lambda_{H(\eta)}\left(h_{\eta_0\eta}(\alpha_0) +
\partial_{H(\eta)} (\beta(\eta))\right) <
\lambda_{H(\eta)}(h_{\eta_0\eta}(\alpha_0)). \leqno (8.31)
$$
This in particular implies that $\partial_{H(\eta)} (\beta(\eta))$
kills $\operatorname{peak}(h_{\eta_0\eta}(\alpha_0))$.

Since we have
$$
\lambda_{H(\eta)}(\beta(\eta)) \leq \lambda_{H(\eta_0)}(h_{\eta\eta_0}(\beta(\eta))
+ C_{(\HH,j)}(|\eta-\eta_0|)
$$
we may assume by Proposition 8.8 that
$$
\lambda_{H(\eta)}(\beta(\eta)) \leq \rho(H(\eta_0);a) + 2\delta + c
$$
if we choose $\eta$ so that $C_{(\HH,j)}(|\eta-\eta_0|) < \delta$,
where the constant $C_{(\HH,j)}(r)$ is the one given in Proposition
6.6.

If $h_{\eta\eta_0}(\operatorname{peak}(\alpha_0))=
\partial_{H(\eta)} \gamma'$ for some $\gamma'$, applying the
homotopy formula between $h_{\eta_0\eta}\circ h_{\eta\eta_0}$ and
the identity, we derive
$$
\operatorname{peak}(\alpha_0) = \partial_{H(\eta_0)}(
h_{\eta_0\eta}(\gamma')) +  \partial_{H(\eta_0)}(\overline H
\operatorname{peak}(\alpha_0)) + \overline H (\partial_{H(\eta_0)}(
\operatorname{peak}(\alpha_0))).
$$
It is easy to show that the levels of the last two terms are less
than equal to
$$
\lambda_{H(\eta_0)}(\operatorname{peak}(\alpha_0)) -
\frac{2\lambda_1(\eta_0)}{3} + E(\HH_{\eta_0\eta}) <
\lambda_{H(\eta_0)}(\operatorname{peak}(\alpha_0)),
$$
if we choose $|\eta_0 - \eta|$ so small that $E(\HH_{\eta_0\eta})
\leq \frac{\lambda_1}{3}$. Therefore $\partial_{H(\eta_0)} \gamma$
with $\gamma = h_{\eta_0\eta}(\gamma')$ kills
$\operatorname{peak}(\alpha_0))$ modulo terms of level less than
$\lambda_{H(\eta_0)}(\operatorname{peak}(\alpha_0))$. Hence we have
$$
\lambda_{H(\eta_0)}(\alpha_0 - \partial_{H(\eta_0)} \gamma) <
\lambda_{H(\eta_0)}(\alpha_0)
$$
which violates tightness of $\alpha_0$. Therefore we may assume
that $h_{\eta_0\eta}(\operatorname{peak}(\alpha_0))$ is not a
boundary.

Furthermore if $\beta(\eta)$ has a generator $[z',w']$ not
connected to $h_{\eta_0\eta}(\operatorname{peak}(\alpha_0))$,
i.e., if $\partial_{H(\eta)}([z',w'])$ is not contributed by any
generator thereof, we can safely replace $\beta(\eta)$ by
$\beta(\eta) - a' [z',w']$ for some $a' \in \Q$ to get rid of the
generator $[z',w']$ without increasing the level of
$h_{\eta_0\eta}(\alpha_0) +
\partial_{H(\eta)}(\beta(\eta))$ :
Since $h_{\eta_0\eta}(\operatorname{peak}(\alpha_0))$ is not a
boundary and since the peak of $h_{\eta_0\eta}(\alpha_0)$ is
contained in $h_{\eta_0\eta}(\operatorname{peak}(\alpha_0))$, this
$\partial \beta(\eta)$ must have a non-zero remainder, after it
kills $\operatorname{peak}( h_{\eta_0\eta}(\alpha_0))$. The level of
the remainder must be lower than
$\lambda_{H(\eta)}(h_{\eta_0\eta}(\alpha_0))$ by (8.31). Therefore
we may assume that all the generators of $\beta(\eta)$ is connected
to some generator of
$h_{\eta_0\eta}(\operatorname{peak}(\alpha_0))$. This in turn
implies that all the generators of
$h_{\eta_0\eta}^{-1}(\beta(\eta))$ are also connected to
$\operatorname{peak}(\alpha_0)$ by (8.29).

We also derive from (8.31)
$$
\lambda_{H(\eta_0)}(h_{\eta_0\eta}^{-1}(\beta(\eta)) <
\rho(H(\eta_0);a) + 3\delta + c. \leqno (8.32)
$$
By the Gromov-Floer compactness, there are only finitely many,
say,
$$
N_0 = N_0(\HH,\eta_0, a, \rho(H(\eta_0);a),
\rho(H(\eta_0);a) + 3\delta + c)
$$
critical points with
$$
[z_1(\eta_0), w_1(\eta_0)], \cdots, [z_{N_0}(\eta_0),w_{N_0}(\eta_0)],
$$
connected to $(\AA_{H(\eta_0)}^{-1}(\rho(H(\eta_0);a)) \cap
\operatorname{Crit} \AA_{H(\eta_0)}$ (and so those connected to
$\operatorname{peak}(\alpha_0)$). Since $\alpha_0$ is tight, there
must exist at least one $[z_\ell,w_\ell]$ for which there is
another critical point $[z_\ell',w_\ell']$ of $\AA_{H(\eta_0)}$
such that
$$
\#(\MM([z_\ell,w_\ell],[z_\ell',w_\ell']) \neq 0
$$
and
$$
\AA_{H(\eta_0)}([z_\ell',w_\ell']) \geq
\lambda_{H(\eta_0)}(\alpha_0). \leqno(8.33)
$$

We divide our discussion into two cases : one for which the cycle
$\alpha_0': = \alpha_0 +
\partial_{H(\eta_0)}(h_{\eta_0\eta}^{-1}(\beta(\eta))$ is tight
and so
$$
\lambda_{H(\eta_0)}(\alpha_0) = \lambda_{H(\eta_0)}(\alpha_0')
\leqno (8.34)
$$
and the other for which $\alpha_0'$ is not tight and so
$$
\lambda_{H(\eta_0)}(\alpha_0') > \lambda_{H(\eta_0)}(\alpha_0).
\leqno (8.35)
$$

We start with the case (8.35). In this case, the inequality (8.33)
must be strict for at least one $\ell$ : Otherwise the cycle
$\alpha_0'$ itself will be tight and so belongs to the category
(8.34). We recall that there are only finitely many, say $N_0 =
N_0(H(\eta_0),a)$, critical points of $\AA_{H(\eta_0)}$ at the
level $y_0 = \rho(H(\eta_0);a)$. This number depends only on
$H(\eta_0)$ and $a$. And we note that there are at most finitely
many $[z_\ell'(\eta_0),z_\ell'(\eta_0)]$ connected to
$[z_\ell(\eta_0),z_\ell(\eta_0)]$, say,
$$
N_1(\ell) = N_1(\HH,\eta_0,a, \rho(H(\eta_0);a),\rho(H(\eta_0);a) +
3\delta + c;\ell)
$$
for each $\ell$. Therefore if we define
$$
N_2 = \sum_{\ell =1}^{N_0} N_1(\ell),
$$
there will be at most $N_2$ such $[z_\ell',w_\ell']$ in total.
Furthermore by the definition of the gap $\lambda_1(\eta_0)$ of the
boundary map $\partial_{H(\eta_0)}$, we have
$$
\rho(H(\eta_0);a) + \lambda_1(\eta_0) \leq
\lambda_{(H(\eta_0))}([z_\ell(\eta_0), w_\ell(\eta_0)]) \leq
\rho(H(\eta_0);a) + 3\delta + c
$$
Then from (8.35) and from the finiteness of such $[z_\ell',w_\ell']$
connected to $[z_\ell,w_\ell]$ which in turn connected to
$[z_0,w_0]$ a generator of $\operatorname{peak}(\alpha_0)$, we
derive that there must be a gap denoted by
$$
\lambda_2(\eta_0)=\lambda_2(\HH,\eta_0,a,\rho(H(\eta_0);a),\rho(H(\eta_0);a)
+ 3\delta +c)
$$
such that
$$
\AA_{H(\eta_0)}([z_\ell',w_\ell']) > \lambda_{H(\eta_0)}(\alpha_0) +
\lambda_2 \leqno (8.36)
$$
for all such $[z_\ell',w_\ell']$. (In Lemma 8.10 below, we will in
fact show that we can choose $\lambda_2$ depending only on
$(\HH,j)$ and $3\delta + c$.)

By nondegeneracy of $H(\eta_0)$ and finiteness of $\#
\operatorname{Per}(H(\eta_0))$, we can choose a sufficiently small
$r:=|\eta-\eta_0| > 0$ so that
$$
C_{(\HH,j)}(r) < \min\left\{\frac{\lambda_2}{4}, \lambda_1 - 2
C_{(\HH,j)}(r) \right\}. \leqno (8.37)
$$
Furthermore
we derive from (8.36) that for any such
$[z_\ell'(\eta_0),w_\ell'(\eta_0)]$ the action of transferred
critical points $[z_\ell'(\eta),w_\ell'(\eta)]$ satisfies
$$
\AA_{H(\eta)}([z_\ell'(\eta),w_\ell'(\eta)]) \geq
\AA_{H(\eta)}([z_0(\eta),w_0(\eta)]) + \frac{\lambda_2}{2} =
\lambda_{H(\eta_0)}(\alpha_0) + \frac{\lambda_2}{2}
$$
by the choice of $\eta$. This then implies that as long as
$|\eta-\eta_0| \leq r$ as above in (8.37), we estimate the level of
the cycle $h_{\eta_0\eta}(\alpha_0)$ as
$$
\lambda_{H(\eta)}(h_{\eta_0\eta}(\alpha_0) + \partial \beta(\eta))
\geq \lambda_{H(\eta_0)}(\alpha_0) + \frac{\lambda_2}{2} >
\lambda_{H(\eta)}(h_{\eta_0\eta}(\alpha(\eta)).
$$
This contradicts to the hypothesis (8.31). This finishes the proof
by setting $\delta_4 = r$ for the case (8.35).

For the case (8.34), the point $(\eta_0, y_0)$ with $y_0
=\lambda_{H(\eta_0)}(\alpha_0)$ is a nondegenerate crossing in the
bifurcation diagram of $\HH$. By the Cerf property, the two
continuous functions
$$
\begin{aligned}
\mu_1(\eta) & : = \lambda_{H(\eta_0)}(h_{\eta_0\eta}(\alpha_0))\\
\mu_2(\eta) & : = \lambda_{H(\eta_0)}(h_{\eta_0\eta}(\alpha_0'))
\end{aligned}
$$
provide the two possible branches at $(\eta_0,y_0)$. By applying the
above consideration to $\alpha_0'$ instead of $\alpha_0$, $\alpha_0$
must belong to the category (1) in the above alternatives and so its
transferred cycles $h_{\eta_0\eta}(\alpha_0')$ must be tight. This
finishes the proof.
\end{proof}

For the later purpose, we provide a precise description of the
above constant $\lambda_2$ in the following lemma and its lower
semi-continuity property.
\medskip

{\noindent \bf Lemma 8.10} {\it Let $K > 0$ and define the
constant $\lambda_2(\eta;K)$ as follows : for $\eta \not\in \SS
ing(\HH)$,
$$
\begin{aligned}
\lambda_2(\eta;K) & : = \inf_{z,C, C'} \{ E_{(H(\eta),J(\eta))}(u)
- E_{(H(\eta),J(\eta))}(u') \mid u\in \MM(J,H; z, (\cdot) ;C),\\
& \quad  u' \in \MM(J,H; z, (\cdot) ;C'),   0 <
E_{(H(\eta),J(\eta))}(u), \, E_{(H(\eta),J(\eta))}(u') \leq K, \\
& \quad E_{(H(\eta),J(\eta))}(u) - E_{(H(\eta),J(\eta))}(u') > 0, \,
\mu_{H(\eta)}(C) = \mu_{H(\eta)}(C') = 1\}
\end{aligned}
$$
and for $\eta \in \SS ing(\HH)$ by the same except that we add the
requirement that at least one of the $\omega$-limits $(\cdot)$ of
$u, \, u'$ is not of the form $z^-$. Then the function $\eta \in
[0,1] \mapsto \lambda_2(\eta;K)$ is lower semi-continuous for each
fixed $K$. In particular we have
$$
\lambda_{2,K} := \min_{\eta \in [0,1]} \lambda_2(\eta;K) > 0.
$$}
\begin{proof}
This is a consequence of Gromov-Floer compactness which we can
apply as in the proofs of Proposition 6.1 and 6.2 using the
following ingredients :
\begin{itemize}
\item[(1)] finiteness of $\#(\operatorname{Per}H(\eta))$ for each
$\eta$, \item[(2)] finiteness of homotopy classes of $u, \, u'$
with the same $\alpha$-limit under the given energy bounds,
\item[(3)] and finally by isolatedness of $\MM(z,(\cdot);C)$ when
$\mu_H(C) = 1$.
\end{itemize}

\end{proof}

\section*{\bf \S 9. Proof of the nondegenerate spectrality}

Finally we are ready to prove Theorem II, which we reiterate here.
Let $H$ be a nondegenerate Hamiltonian and
$$
a = \sum a_A q^{-A}, \quad a_A \in H^*(M) \leqno (9.1)
$$
be a non-zero quantum cohomology class.

\medskip{\noindent \bf  Theorem 9.1.} {\it Let $H$ be any  nondegenerate
one-periodic Hamiltonian and $J$ be such that $(H,J)$ is Floer
regular. Then for any nonzero quantum cohomology class $a$, the
mini-max value $\rho(H;a)$ is realized by the level of a tight Floer
cycle of $(H,J)$. }
\medskip

Besides the preparatory materials proven in section 3-7, our proof
of Theorem II also exploits the following two ingredients in an
essential way:

\begin{itemize} \item[(1)] $\rho(\e f;a)$ is homologically essential for any
Morse function if $\e>0$ is sufficiently small.

\item[(2)] The mini-max value $\rho(H;a)$ is tied to a nontrivial
topological property, not an accidental critical value.

\end{itemize}

Now we fix a Morse function $f$ and a Cerf homotopy $\HH \in
\PP^{Cerf}(\HH_m(M);\e f, H)$ satisfying the properties described
in section 7,
$$
\HH=\{H(\eta)\}_{0 \leq \eta \leq 1},
$$
such that
$$
H(0) = \e f, \quad H(1) = H.
$$
For simplicity, without loss of any generality, we will assume
that {\it the values of the Morse function $f$ at critical points
are all distinct}.

We also consider the set of $j = \{J(\eta)\}_{0\leq \eta \leq 1} $
lying in $\PP(j_\omega;\HH)$, In particular, by the choice of $j$,
the Floer homology $HF_*(H(\eta),J(\eta))$ is defined for any
$\eta \in I(\HH,j)$. We
\medskip

The proof of Theorem 9.1 will be done by a continuation argument. We
define a subset of [0,1]
$$
S(\HH)  =\{\eta \in [0,1] \mid \text{ $\eta$ satisfies one of the
following two conditions} \}: \leqno (9.2)
$$
\begin{itemize} \item[(1)]  when $\eta \in [0,1] \setminus \SS ing(\HH)$,
$\rho(H(\eta);a)$ is homologically essential, i.e., there is a $J$
for which $(H,J)$ is Floer-regular and carries a Floer cycle
$\alpha$ with $\rho(H(\eta);a) = \lambda_H(\alpha)$.

\item[(2)] when $\eta \in \SS ing(\HH)$, there is $j \in
\PP(j_\omega;\HH)$ such that it carries a sequence of Floer points
$\eta_i \to \eta$ for which the mini-max value $\rho(H(\eta_i);a)$
is homologically essential for each $i$.
\end{itemize}
\medskip

We will prove $S(\HH) = [0,1]$ by a continuation argument starting
from $H(0) = \e f$, which will in particular prove Theorem II.
\medskip

\subsection*{\it Step 1: $S(\HH)$ is nonempty}

\medskip

We will show that $\rho(\e f;a)$ is a homologically essential
critical value.

Consider an almost complex structure $J_0$ such that $-\e f$ is
Morse-Smale for the metric $g = \omega(\cdot, J_0\cdot)$. We
denote by $CM_*(-\e f)$ the graded Morse complex associated to
$(-\e f,J_0)$. An element of $CM_\ell(-\e f)$ has the form
$$
\sum_{k} a_k [p_k], \quad a_k \in\Q, \, p_k \in \text{Crit}_\ell(-\e
f).
$$
For given $p \in \text{Crit}_*(-\e f)$ and $A \in \Gamma$, we denote
$p \cdot q^A := [p, \widehat p \# A]$ where $\widehat p$ is the
constant disc $p$.

We represent $a^\flat$ by a Novikov Morse cycle
$$
\gamma = \sum_A\gamma_{A} q^{A}, \quad \gamma_j \in CM_*(-\e f)
$$
of $- \e f$. As we argued in [section 5, 15], we may assume that the
representative $\gamma$ of $a^\flat$ has the form
$$
\gamma = \sum_{A \in \Gamma(a)}\gamma_A q^A \leqno (9.3)
$$
where $\gamma_A$ is a Morse cycle of $- \e f$. By the definition of
Novikov Floer chains, we can enumerate $\Gamma(a)$ so that
$\lambda_j = \omega(A_j)$ satisfy
$$
\lambda_1 > \lambda_2 > \cdots.
$$
Then it is easy to see
$$
\lambda_{\e f}(\gamma) = \lambda_{\e f}(\gamma_1 q^{A_1}), \leqno
(9.4)
$$
provided $\e > 0$ satisfies
$$
\e (\max f - \min f) \leq c(a)= \lambda_1 - \lambda_2. \leqno (9.5)
$$
(See the proof of [Lemma 7.4, 15].) On the other hand, we compute
$$
\lambda_{\e f}(\gamma_1 q^{A_1}) = \AA_{\e f}(p_{\gamma_1} q^{A_1})
= - \omega(A_1) - \e f(p_{\gamma_1})
$$
where $p_{\gamma_1}$ is the peak of the Morse cycle $\gamma_1$
measured by the values of $- \e f$. Recall that {\it as long as we
fixed the quantum cohomology class $a$}, the collection $\Gamma(a)$
is fixed and so $A_1$ is fixed and the coefficient Morse cycles
$\gamma_1$ varies inside the homology classes $PD(a_1) \in
H_*(M,\Q)$.

Therefore if we choose $\e$ so that (9.5) is satisfied, the mini-max
value $\rho(\e f;a)$ becomes
$$
\rho(\e f;a) = -\omega(A_1) + \inf_{\gamma_1 \in PD(a_1)}\max_{p}\{-
\e f(p) \mid p \in \gamma_1 \}.
$$
However it is an easy consequence of {\it compactness} of $M$ and is
well-known in the finite dimensional critical point theory that on a
compact manifold $M$, the mini-max value
$$
\inf_{\gamma_1 \in PD(a_j)}\max_{p}\{- \e f(p) \mid p \in \gamma_1
\}
$$
can be realized by the level of a tight Morse cycle $\gamma_{tgt}$.
Then if we re-choose $\gamma_1 = \gamma_{tgt}$ and fix other
$\gamma_j$'s for $j \geq 2$,  we have
$$
\rho(\e f;a) = \AA_{\e f}([p_{(f;a_1)}, p_{(f;a_1)}q^{A_1}]) =
\lambda_{\e f}(\gamma). \leqno (9.6)
$$
Here $p_{(f;a_1)}$ is the unique critical point  of $f$ that is the
peak of the tight Morse cycle $\gamma_1$ associated to the
cohomology class $a_1 \in H^*(M)$. This proves that $\rho(\e f;a)$
is a homologically essential critical value.

\medskip

\subsection*{\it Step 2: $S(\HH)$ is open in $[0,1]$}

\medskip

Let $\eta_0 \in S(\HH)$. We would like to show that there exists
$\delta > 0$ such that $(\eta_0-\delta,\eta_0+\delta) \subset
S(\HH)$. We consider two cases separately: one is the case where
$\eta_0 \in [0,1] \setminus \SS ing(\HH)$ and the other the case
with $\eta_0 \in \SS ing(\HH)$.

For the case where $\eta_0\in [0,1] \setminus (\SS ing(\HH)$, we
choose $j \in \PP^{sub}(j_\omega;\HH;\eta_0)$ i.e., $\eta_0$ is a
Floer point for $(\HH,j)$. See section 4 for the definition of
$\PP^{sub}(j_\omega;\HH;\eta_0)$. Then the stability theorem,
Theorem 8.3 implies existence of such $\delta > 0$.

Next, we consider the case when $\eta_0$ is in $\SS ing(\HH)\cap
S(\HH)$. In this case, by the definition of $S(\HH)$, there is $j
\in \PP(j_\omega;\HH)$ for which we have a sequence of the Floer
points $\eta_k \to \eta_0$ and tight Floer cycles $\alpha_k \in
CF_*(H(\eta_k))$. Under this assumption, we would like to prove that
there exists $\delta_6 = \delta_6(\HH,j;\eta_0)
> 0$ such that all $\eta \in (\eta_0 - \delta_6, \eta_0 + \delta_6)\setminus \{\eta_0\}$
allows tight cycles with its level $\rho(H(\eta);a)$. As in
section 8, we will assume, without loss of any generality,
$$
\eta_k \nearrow \eta_0.
$$

We denote $[z_k,w_k]$ be a peak of $\alpha_k$. Then we have
$$
\lim_{k \to \infty} \AA_{H(\eta_k)}([z_k,w_k]) = \lim_{k \to
\infty}\rho(H(\eta_k);a) = \rho(H(\eta_0);a) \leqno (9.7)
$$
by the continuity of $\rho_a$. After choosing a subsequence, we may
assume that $z_k$ converges to a periodic orbit $z_\infty \in
\text{Per}(H(\eta_0))$.

Using the Cerf property of $\HH$, there are two cases to consider:
one is the case where
$$
z_k \neq z^\pm(\eta_k) \leqno (9.8)
$$
and the other
$$
z_k = z^+(\eta_k) \quad \text{or} \quad z_k = z^-(\eta_k) \leqno
(9.9)
$$
after choosing a subsequence of $\eta_k$ if necessary.

Note that for the case (9.9), the limit orbit is nothing but
$$
z_\infty = z_0
$$
where $z_0$ is the unique degenerate periodic orbit of $H(\eta_0)$
and for the case (9.8), the limit $z_\infty$ is far away from
$z_0$. The following proposition reduces the proof to the case of
(9.8).
\medskip

{\noindent \bf Proposition 9.2} {\it Let $a\neq 0$ be a given
quantum cohomology class and denote by $\alpha$ a Floer cycle with
$[\alpha] = a^\flat$. Then there exists $0 < \delta_9 \leq
\delta_1$ with $\delta_9 = \delta_9(\HH,j)$ such that at any Floer
point $\eta$ satisfying $d(\eta,\SS ing(\HH)) < \delta_9$ there
exists a tight Floer cycle $\alpha$ of $(H(\eta),J(\eta))$, no
peak of which is of the form
$$
[z^\pm(\eta), w]
$$
for any bounding disc $w$.}
\begin{proof}
We have already shown in Proposition 6.9 that no peak of $\alpha$ is
of the form $[z^+(\eta),w]$. We now prove applying Theorem 6.7 (2)
and tightness of $\alpha$ that there exists $\delta_9$ such that
whenever $|\eta-\eta_0| < \delta_9$, if there is a peak of the tight
cycle $\alpha$ having the form $[z^-(\eta),w^-(\eta)]$, we can
cancel the peak by adding $\partial_{(H(\eta),J(\eta))}(c\cdot
[z^+(\eta),w^+(\eta)])$ with $w^+(\eta) = w^-(\eta) \#
u^{can}_{z^-(\eta)z^+(\eta)}$ for a suitable $c \in \Q$ keeping the
level unchanged and hence keeping tightness of $\alpha$ as well.
Applying this cancelling repeatedly, we can cancel all such peaks.
This finishes the proof.
\end{proof}

We therefore safely assume that we are in the case of (9.8). Then
an examination of the proof of Theorem 8.3 proves the following
\medskip

{\noindent \bf Lemma 9.3.} {\it There exists $\delta'_4 =
\delta'_4(\HH,j;a,\eta_0)> 0$ and a sufficiently large $N \in \N$
such that the transferred cycles
$$
h_{(\eta_N\eta;\rho_0)}(\alpha_{N,\pm})
$$
are tight on each of the semi-intervals of $(\eta_0 -
\delta_4',\eta_0 + \delta_4')\setminus \{\eta_0\}$.}
\begin{proof}
The proof is a variation of that of Proposition 8.4 except that in
the current case, $\eta_0 \in \SS ing(\HH)$ and so the Floer
homology itself at $\eta_0$ is not defined. On the other hand, by
the Cerf property of $\HH$, there is no critical point $[z,w]$ of
$\AA_{H(\eta_0)}$ at the same level of the form
$\AA_{H(\eta_0)}([z_0,w_0])$ where $z_0$ is a degenerate periodic
orbit.

Consider the constants $\lambda_1$ and $\lambda_2$ that appears in
the proof of Theorem 8.3. An examination of the proof shows that we
can choose $\lambda_1$ depending only on $(\HH,j)$ but independent
of $\eta_0$ if we replace $A_{(H(\eta_0),J(\eta_0))}$ by the
constant $A_{(\HH,j)}$ provided in Theorem 6.7.

For the constant $\lambda_2=\lambda_{2,K}$, we take $\lambda_{2,K}$
given in Lemma 8.10 for
$$
K = 3\delta + c
$$
where $\delta, \, c$ are as in the proof of Theorem 8.3. Once we
have this, the same argument as that of Theorem 8.3 with
$\lambda_2(\eta_0)$ replaced by $\lambda_{2,K}$ provides the
constant $\delta_4' =\delta'_4(\HH,j;a,\eta_0)> 0$ which finishes
the proof.
\end{proof}

If we choose $N$ sufficiently large, then $\eta_0 \in (\delta_4' -
\eta_N,\delta_4' + \eta_N)$ and hence follows openness of $S(\HH)$
at $\eta_0$. Combining the above altogether, we have finished the
proof of openness of $S(\HH)$.
\medskip

We remark that the corresponding level functions
$$
\begin{aligned}
\mu_-(\eta) & = \lambda_{H(\eta)}(h_{\eta_0\eta}(\alpha_{N,-}))\\
\mu_+(\eta) & = \lambda_{H(\eta)}(h_{\eta_0\eta}(\alpha_{N,+}))
\end{aligned}
$$
together define a continuous function on $(\eta_N -
\delta_4',\eta_N + \delta_4')$ which extends continuously across
$\eta_0$.

\subsection*{\it Step 3: $S(\HH)$ is closed in $[0,1]$}

\smallskip

Let $\eta_k \to \eta_\infty$ be a sequence of Floer points such
that each $\rho(H(\eta_i);a)$ is a homologically essential
critical value. First consider the case where $\eta_\infty$ lies
in $[0,1] \setminus \SS ing(\HH)$. We choose $j \in
\PP^{sub}(j_\omega;\HH;\eta_\infty)$ and consider the transferred
cycles
$$
\alpha_k(\eta): = h_{(\eta_k\eta;\rho_0)}(\alpha_k).
$$
Proposition 8.2 shows that the existence of tight cycles at
$\eta_k$ does not depend on the choice of $j$ and so we may assume
that all $\eta_k$ including $\eta_\infty$ are Floer points of
$(\HH,j)$. Then by the same arguments used for the proof of
Theorem 8.3 and using Lemma 8.10, there is a constant $\delta_4 =
\delta_4(\eta_\infty)
> 0$ such that on each of the semi-intervals of $(\eta_k -
\delta_4,\eta_k]$ or $[\eta_k,\eta_k + \delta_4)$, we can find a
tight cycle $\alpha_k$ for which $\alpha_k(\eta)$ is also tight on
the corresponding intervals. Obviously if $k$ is sufficiently
large, then $\eta_k \in (\eta_\infty - \delta,\eta_\infty +
\delta)$. In particular,
$h_{(\eta_k\eta_\infty;\rho_0)}(\alpha_k)$ is a required tight
Floer cycle of $H(\eta_\infty)$ at the level
$\rho(H(\eta_\infty);a)$. This takes care of the case when
$\eta_\infty$ lies in $[0,1] \setminus \SS ing(\HH)$.

On the other hand when $\eta_\infty \in S(\HH) \cap \SS ing(\HH)$,
there is nothing to prove by the definition of $S(\HH)$. This proves
that $S(\HH)$ is closed.

\medskip

\subsection*{\it Step 4: Finish-up of the proof}

\smallskip

Combining Step 1-3 and the fact that $[0,1]$ is connected, we have
proved $S(\HH) = [0,1]$ and so the proof of Theorem 9.1. This in
turn finishes the proof of Theorem II at last.

\medskip{\noindent \bf Remark 9.4.}
In fact, an examination of the above proof seems to show that the
spectrality axiom holds for the Hamiltonians either of the
Bott-Morse types or {\it of a finite type}: We call a Hamiltonian
$H$ a finite type, if all of its periodic orbits are isolated and
the degeneracy of the linearization is of finite order. It would be
interesting to see if the spectrality axiom holds for arbitrary
Hamiltonians or not.
\medskip

We are now ready to prove Theorem V stated in the introduction as a
by-product of the arguments used in the above proof. We re-state the
theorem here.
\medskip

{\noindent \bf Theorem 9.5.} {\it Let $\HH$ be a Cerf-homotopy. Then
the spectral function $s \mapsto \rho(H^s;a)$ is smooth away from a
countable subset of $\CC^{nd} ross(\HH)$. }
\begin{proof}
Let $\eta_0 \in (0,1) \setminus \SS ing(\HH)$. We choose a generic
$J^{\eta_0}$ such that the pair $(H(\eta_0),J^{\eta_0})$ is
Floer-regular and then extend $J^{\eta_0}$ to a family $j =
\{J^\eta\}_{0 \leq \eta\leq 1}$ so that the pair $(\HH,j)$ is a
Floer homotopy, i.e., satisfies the properties of Theorem 4.6. Then
the above proof shows that  if $\eta_0$ is not in $\CC^{nd}
ross(\HH)$, there exists a tight Floer cycle $\alpha_0$ for
$(H(\eta_0),J^{\eta_0})$ such that the function $\mu$
$$
\mu(\eta) = \lambda_{H(\eta)}(h_{\eta_0\eta}(\alpha_0))
$$
is well-defined on $(\eta_0-\delta,\eta_0+\delta)$ for some $\delta
> 0$ which is smooth thereon. Furthermore by making $\delta > 0$ smaller if necessary,
the proof of Theorem 8.3 shows that the cycles
$h_{\eta_0\eta}(\alpha_0)$ are all tight and hence $ \rho(H(\eta);a)
= \mu(\eta)$. On the other hand if $\eta_0 \in \CC^{nd} ross(\HH)$,
the last paragraph of the proof of Theorem 8.3 shows that
$\rho(H(\eta);a)$ may be realized by two different branches
$\mu_-:(\eta_0 - \delta,\eta_0] \to \R$ and $\mu_+:[\eta_0,\eta_0 +
\delta) \to \R$ with $\mu_-(\eta_0) = \mu_+(\eta_0)$ given by
$$
\begin{aligned}
\mu_-(\eta) & = \lambda_{H(\eta)}(h_{\eta_0\eta}(\alpha_-))\\
\mu_+(\eta) & = \lambda_{H(\eta)}(h_{\eta_0\eta}(\alpha_+))
\end{aligned}
$$
where $\alpha_\pm$ are two different tight Floer cycles of
$H(\eta_0)$. This proves that the function $\eta \mapsto
\rho(H(\eta);a)$ is differentiable possibly except at such points
from $\CC^{nd} ross(\HH)$. This finishes the proof.
\end{proof}

\section*{\bf \S10. Spectral invariants of Hamiltonian diffeomorphisms}

We recall that the invariants $\rho(H;a)$ were constructed for
arbitrary Hamiltonian functions $H: S^1 \times M \to \R$ in
\cite{Oh5}. We first summarize the basic properties of the
invariants $\rho=\rho(H;a)$. Except the nondegenerate spectrality
axiom proven in the present paper, all other axioms are proved in
\cite{Oh5}.
\medskip

{\noindent\bf Theorem 10.1.} {\it Let $(M,\omega)$ be arbitrary closed
symplectic manifold. For any given quantum cohomology class $0
\neq a \in QH^*(M)$, we have a continuous function denoted by
$$
\rho =\rho(H; a): C_m^\infty(S^1 \times M) \times QH^*(M) \to \R
$$
such that they satisfy the following axioms: Let $H, \, F  \in
C_m^\infty(S^1 \times M)$ be smooth Hamiltonian functions and $a
\neq 0 \in QH^*(M)$. Then $\rho$ satisfies the following axioms:

\begin{itemize} \item[(1)] {\bf (Nondegenerate spectrality)} For each $a\in
QH^*(M)$, $\rho(H;a) \in \text{Spec}(H)$ if $H$ is nondegenerate.

\item[(2)] {\bf (Projective invariance)} $\rho(H;\lambda a) =
\rho(H;a)$ for any $0 \neq \lambda \in \Q$.
\item[(3)] {\bf
(Normalization)} For $a = \sum_{A \in \Gamma} a_A q^{-A} $, we
have $\rho(\underline 0;a) = v(a)$ where $\underline 0$ is the
zero function and
$$
v(a): = \min \{\omega(-A) ~|~  a_A \neq 0 \} = - \max \{\omega(A)
\mid a_A \neq 0 \}.
$$
is the $($upward$)$ valuation of $a$.
\item[(4)] {\bf (Symplectic
invariance)} $\rho(\eta^*H ;\eta^*a) = \rho(H ;a)$ for any
symplectic diffeomorphism $\eta$
\item[(5)] {\bf (Triangle inequality)}
$\rho(H \# F; a\cdot b) \leq \rho(H;a) + \rho(F;b) $
\item[(6)] {\bf
($C^0$-continuity)} $|\rho(H;a) - \rho(F;a)| \leq \|H \# \overline
F\| = \|H - F\|$ where $\| \cdot \|$ is the Hofer's pseudo-norm on
$C_m^\infty(S^1 \times M)$. In particular, the function $\rho_a: H
\mapsto \rho(H;a)$ is $C^0$-continuous.
\end{itemize}
}
\medskip

(In the symplectic invariance axiom, $\eta^*a$ denotes the canonical
pull-back of $a$ under the symplectic diffeomorphism $\phi$. In
general $\eta^*a \neq a$ {\it unless $\phi$ is symplectically
isotopic to the identity}. In \cite{Oh5}, the symplectic invariance
was stated incorrectly as $\rho(\eta^*H ;a) = \rho(H ;a)$. We thank
Polterovich for pointing out this error.)

By the one-one correspondence between (normalized) $H$ and its
associated Hamiltonian path $\phi_H: t \mapsto \phi_H^t$, one can
regard the spectral function
$$
\rho_a: \HH_m(M) = C^\infty_m([0,1]\times M) \to \R
$$
as a function defined on $\PP(Ham(M,\omega);id)$, i.e.,
$$
\rho_a: \PP(Ham(M,\omega), id) \to \R \leqno (10.1)
$$
as described in \cite{Oh5}. Here we denote by $\PP(Ham(M,\omega),
id)$ the set of smooth Hamiltonian paths in $Ham(M,\omega)$ and by
$\widetilde{Ham}(M,\omega)$ the set of path homotopy classes on
$\PP(Ham(M,\omega), id)$, i.e., the (\'etale) universal covering
space of $Ham(M,\omega)$. We equip $\widetilde{Ham}(M,\omega)$ with
the quotient topology. The following corollary shows that the
function $\rho_a$ pushes down to $\widetilde{Ham}(M,\omega)$ as a
continuous function. The proof of the following fact is entirely the
same as in the rational case \cite{Oh5} summarized in the
introduction of the present paper, and so omitted.
\medskip

{\noindent\bf Corollary 10.2.} {\it Let $(M,\omega)$ be an arbitrary closed
symplectic manifold. For any nondegenerate $\widetilde \phi \in
\widetilde{Ham}(M,\omega)$, we have
$$
\rho(H;a) = \rho(K;a) \leqno (10.2)
$$
for any smooth functions $H \sim K$ satisfying $[H]=[K] =
\widetilde \phi$.
}
\medskip

We like to emphasize that at this moment, because we do not know
validity of the spectrality axiom for degenerate Hamiltonians, we
do not have the equality
$$
\rho(H;a) = \rho(K;a)
$$
yet for two $H$ and $K$ representing the same $\widetilde\phi$, if
the latter is degenerate. The scheme of the proof used to prove
(1.3) for the rational $(M,\omega)$ cannot be applied without the
spectrality axiom. In this regard, the following is still a
non-trivial theorem to prove. The argument used in the proof is
similar to the proof of \cite[Lemma 5.1]{Oh2} or \cite[Theorem
5.1]{Oh5}.
\medskip

{\noindent\bf Theorem 10.3.} {\it For any pair $(H,K)$ satisfying $H \sim
K$, we have
$$
\rho(H;a) = \rho(K;a). \leqno (10.3)
$$
}
\begin{proof} For the nondegenerate case, Corollary 10.2 proves
(10.3). It remains to prove (10.3) for the degenerate cases.

Suppose $H \sim K$. We approximate $H$ and $K$ by sequences of
nondegenerate Hamiltonians $H_i$ and $K_i$ in the $C^\infty$
topology respectively. We note that the Hamiltonian
$$
K \# H_i \# \overline K
$$
generates the flow $\phi_K^t\circ \phi_{H_i}^t \circ
(\phi_K^t)^{-1}$, which is conjugate to the flow $\phi_{H_i}^t$
and is nondegenerate. Therefore we have
$$
\rho(H_i;a) = \rho(K \# H_i \# \overline K;a) \leqno (10.4)
$$
by the symplectic invariance of $\rho$. On the other hand, since
$H\sim K$, we have
$$
K \# H_i \# \overline K \sim K \# H_i \# \overline H.
$$
Since both are nondegenerate, Corollary 10.2 implies
$$
\rho(K \# H_i \# \overline K;a) = \rho(K \# H_i \# \overline H;a).
\leqno (10.5)
$$
Here we remind the readers that the definition of the {\it
nondegeneracy} of a Hamiltonian depends only on its time-one map.
By taking the limits of (10.4) and (10.5) and using the continuity
of $\rho(\cdot;a)$, we get
$$
\rho(H;a) = \rho(K \# H \# \overline K;a) = \rho(K\# H \#
\overline H;a) = \rho(K;a)
$$
where the last equality comes since $H \# \overline H = 0$. Hence
the proof.
\end{proof}

Therefore, for any nondegenerate $\widetilde\phi$, we can define
the value $\rho_a(\widetilde \phi)$ by
$$
\rho_a(\widetilde \phi) : = \rho(H;a) \leqno (10.6)
$$
for $H$ satisfying $[H] = \widetilde \phi$. This defines a
well-defined continuous function
$$
\rho_a: \widetilde{Ham}^{nd}(M,\omega) \to \R \leqno (10.7)
$$
where $\widetilde{Ham}^{nd}(M,\omega)$ is the subset of
$\widetilde{Ham}(M,\omega)$ consisting of nondegenerate
$\widetilde \phi$'s.
\medskip

{\noindent\bf Theorem 10.4.} {\it The function $\rho_a$ defined by {\rm(10.7)}
extends to continuously $\widetilde{Ham}(M,\omega)$ in the
quotient topology of $\widetilde{Ham}(M,\omega)$ induced from
$$\PP(Ham(M,\omega), id).$$
}
\begin{proof} Recall the definition of the quotient topology under
the projection
$$
\pi:\PP(Ham(M,\omega), id) \to \widetilde{Ham}(M,\omega).
$$
We proved that the assignment
$$
H \mapsto \rho(H;a) \leqno (10.8)
$$
is continuous on $C^\infty([0,1]\times M)$ in \cite{Oh5}. By the
definition of the quotient topology,
$$
\rho_a : \widetilde{Ham}(M,\omega) \to \R
$$
is continuous, because the composition
$$
\rho_a \circ \pi: \PP(Ham(M,\omega), id) \to \R,
$$
which is nothing but (10.8), is continuous.
\end{proof}
\medskip

{\noindent\it Remark 10.5.}
One cheap way of defining a function on the covering space
$\widetilde{Ham}(M,\omega)$ out of the spectral invariants
$\rho(H;a)$ is to take the infimum of $\rho(H;a)$ among all $H$'s
connecting the identity to the given time one map $\phi=\phi_H^1$
in the same homotopy class of paths: i.e.,
$$
\rho'(\widetilde \phi;a) : = \inf_{\widetilde \phi = [H]}\rho(H;a)
$$
similarly as in the cases of Hofer pseudo-norm
$$
\|\widetilde \phi\| = \inf_{\widetilde \phi = [H]} \|H\|.
$$
However, {\it unless we had the spectrality axiom}, more
specifically without its consequence Theorem 10.3, there would be
no reason why the function $\rho': \widetilde{Ham}(M,\omega) \to
\R$ is continuous just as the function $\widetilde \phi \mapsto
\|\widetilde \phi\|$ is not a priori continuous in the (\'etale)
topology of $\widetilde{Ham}(M,\omega)$.

Finally, for the reader's convenience, we summarize the basic
axioms of the invariant $\rho: \widetilde{Ham}(M, \omega) \times
QH^*(M) \to \R$ in the following theorem, whose proofs immediately
follow from Theorem 10.1 and 9.4.
\medskip

{\noindent\bf Theorem 10.6.} {\it Let $(M,\omega)$ be any closed symplectic
manifold. Let $\widetilde \phi, \, \widetilde \psi \in
\widetilde{Ham}(M,\omega)$ and $0 \neq a \in QH^*(M)$. Then for
each $0 \neq a \in QH^*(M)$, the function
$$
\rho_a: \widetilde{Ham}(M,\omega) \to \R
$$
is continuous, and the function
$$
\rho: \widetilde{Ham}(M,\omega) \times QH^*(M) \to \R
$$
satisfies the following axioms:
\begin{itemize}
\item[(1)] {\bf (Nondegenerate
spectrality)} For each $a\in QH^*(M)$, $\rho(\widetilde \phi;a)
\in \text{Spec}(\widetilde \phi)$, if $\widetilde \phi$ is
nondegenerate.
\item[(2)] {\bf (Projective invariance)}
$\rho(\widetilde\phi;\lambda a) = \rho(\widetilde\phi;a)$ for any
$0 \neq \lambda \in \Q$.
\item[(3)] {\bf (Normalization)} For $a =
\sum_{A \in \Gamma} a_A q^{-A} $, we have $\rho(\underline 0;a) =
v(a)$ where $\underline 0$ is the identity in
$\widetilde{Ham}(M,\omega)$ and
$$
v(a): = \min_A \{\omega(-A) ~|~  a_A \neq 0 \} = - \max
\{\omega(A) \mid a_A \neq 0 \}.
$$
is the (upward) valuation of $a$.
\item[(4)] {\bf (Symplectic
invariance)} $\rho(\eta^{-1} \widetilde \phi \eta;\eta^*a) =
\rho(\widetilde \phi;a)$ for any symplectic diffeomorphism $\eta$
\item[(5)] {\bf (Triangle inequality)} $\rho(\widetilde \phi \cdot
\widetilde \psi; a\cdot b) \leq \rho(\widetilde \phi;a) +
\rho(\widetilde \psi;b) $
\item[(6)] {\bf ($C^0$-continuity)}
$|\rho(\widetilde \phi;a) - \rho(\widetilde \psi;a)| \leq
\|\widetilde \phi \circ \widetilde\psi^{-1} \| $ where $\| \cdot
\|$ is the Hofer's pseudo-norm on $\widetilde{Ham}(M,\omega)$. In
particular, the function $\rho_a: \widetilde \phi \mapsto
\rho(\widetilde \phi;a)$ is $C^0$-continuous.
\end{itemize}
}

We like to remind the readers the spectrality axiom holds for any
$\widetilde \phi$ for the rational symplectic manifolds \cite{Oh5}.
It is an interesting problem to investigate to understand what kind
of Hamiltonians, besides nondegenerate ones, satisfy the spectrality
axiom, which will be a subject of future study.

\section*{\bf \S 11. Applications to Hofer's geometry}

In this section, we provide two immediate applications of the
spectrality axiom in the study of Hofer's geometry of the
Hamiltonian diffeomorphism group. We first recall the following
definitions
$$
\begin{aligned} E^-(\widetilde \phi) & = \inf_{[\phi,H] = \widetilde
\phi}
\int_0^1 - \min H_t \, dt \\
E^+(\widetilde \phi) & = \inf_{[\phi,H] = \widetilde \phi}
\int_0^1 \max H_t \, dt
\end{aligned}
\leqno (11.1)
$$
(See \cite{Po}, and \cite{Mc} for example). Note that we have
$$
E^-(\overline H) = E^+(H)
$$
for the Hamiltonian $\overline H \mapsto \phi^{-1}$. This in turn
implies
$$
E^-(\widetilde \phi^{-1}) = E^+(\widetilde\phi)
$$
and hence
$$
0 \leq E^+(\widetilde \phi) + E^-(\widetilde \phi) \leq
\inf_{[\phi,H] = \widetilde \phi} \int(\max H_t - \min H_t) \, dt.
\leqno (11.2)
$$
In particular we have the inequality for the Hofer pseudo-norm
$\|\phi\|$ and $\|\phi\|_{med}$
$$
\|\phi\|_{med}:= \inf_{\pi(\widetilde \phi) = \phi}(E^+(\widetilde
\phi) + E^-(\widetilde \phi)) \leq \|\phi\|. \leqno (11.3)
$$
Now we consider the invariant $\rho(H;1)$. We have shown, by
definition,
$$
\rho(\widetilde \phi;1) = \rho(H;1) \leqno (11.4)
$$
for any $H$ with $\widetilde \phi = [H]$. Since we have
$$
\rho(H;1) \leq E^-(H)
$$
(11.4) implies
$$
\rho(\widetilde \phi;1) \leq E^-(H) \leqno (11.5)
$$
for all $H$ with $[H] = \widetilde \phi$. By taking the infimum of
(11.5) over all such $H$, we have proved the following inequality
$$
\rho(\widetilde \phi;1) \leq E^-(\widetilde \phi). \leqno (11.6)
$$
Applying the same argument to $\overline H$ and using
$E^-(\overline H) = E^+(H)$, we obtain
$$
\rho(\widetilde \phi^{-1};1) \leq E^+(H). \leqno (11.7)
$$
We now define
$$
\widetilde \gamma(\widetilde \phi) = \rho(\widetilde \phi;1) +
\rho(\widetilde \phi^{-1};1).
$$
Furthermore we have
$$
\widetilde \gamma(\widetilde \phi) = \rho(\widetilde \phi;1) +
\rho(\widetilde \phi^{-1};1) \geq \rho(id;1) = 0.
$$
We recall the definition of the spectral norm $\gamma$ from
\cite{Oh6}
$$ \gamma(\phi) = \inf_{H \mapsto \phi}\Big(\rho(H;1) +
\rho(\overline H;1)\Big) = \inf_{\pi(\widetilde\phi) = \phi}
(\widetilde\gamma(\widetilde\phi)).
$$

Combining (11.5) and (11.7), we have proved
\medskip

{\noindent\bf Theorem 11.1.} {\it For any $\widetilde \phi \in
\widetilde{Ham}(M,\omega)$, we have
$$
\widetilde \gamma(\widetilde \phi) \leq \|\widetilde \phi\|. \leqno
(11.8)
$$
In particular, we have
$$
\gamma(\phi) \leq \|\phi\|_{med}. \leqno (11.9)
$$
}

(11.9) is an improvement of the inequality $\gamma(\phi) \leq
\|\phi\|$ proven in \cite{Oh6} and give a different proof of
nondegeneracy, with a lower bound, of the {\it medium Hofer norm}
$\|\cdot \|_{med}$. Previously McDuff \cite{Mc} proved that this is
nondegenerate by a different method.

Next we define a partial order on $\widetilde{Ham}(M,\omega)$ out
of the invariant $\rho(\widetilde \phi;1)$. We note that
$\rho(\widetilde \phi;1)$ can be strictly negative although the
sum $\rho(\widetilde \phi;1) + \rho(\widetilde \phi^{-1};1)$
cannot. We recall our convention of the action functional is
$$
\AA_H([z,w]) = - \int w^*\omega - \int H(t,z(t))\, dt
$$
emphasizing the `$-$' in front of the integral $\int H(t,z(t))\,
dt$, and that $\rho(H;a)$ is defined in terms of the action
functional, not in terms of $H$ itself. This explains somewhat
contradictory usage of positiveness in the definition.
\medskip

{\noindent\bf Definition 11.2.} We say that a Hamiltonian $H$
is {\it homologically positive} if
$$
\phi(H;1) \leq 0,
$$
and $\widetilde \phi \in \widetilde{Ham}(M,\omega)$ is {\it
homologically positive} if its representing normalized Hamiltonian
$H$ is homologically positive. We also call the corresponding
Hamiltonian path $\{\phi_H^t\}_{0 \leq t\leq 1}$ {\it
homologically positive} if $[\phi,H]=\widetilde \phi$ is positive.
We define
$$
\begin{aligned} \widetilde{Ham}_+(M,\omega) & = \{\widetilde \phi \mid
\widetilde \phi \,\, \text{homologically positive} \} \\
C^+_m([0,1]\times M) & = \{H \in C^\infty_m(S^1 \times M) \mid H
\,\, \text{homologically positive} \}
\end{aligned}
$$
and denote
$$
\begin{aligned} \PP^+(Ham(M,\omega),id) & = \{f: [0,1] \to Ham(M,\omega)
\mid f(0) = id, \\
& f(t) = \phi_H^t, \, H \in C_m^+(S^1 \times M) \}
\end{aligned}
$$
for the set of homologically positive Hamiltonian paths issued at
the identity.
\medskip

We refer readers to \cite{ElP} for a general discussion on partially
ordered groups and the definition of the {\it normal cone} used
below in the context of contact transformations.
\medskip

{\noindent\bf Theorem 11.3.} {\it The subset
$C:=\widetilde{Ham}_+(M,\omega)$ forms a normal cone in
$\DD:=\widetilde{Ham}(M,\omega)$, i.e., $C$ satisfies \begin{itemize}
\item[(1)] If $f, \, g \in C$, $fg \in C$
\item[(2)] If $f \in C$ and $h \in
\DD$, $hfh^{-1} \in C$
\item[(3)] $id \in C$
\end{itemize}
}
\begin{proof} (1) follows from the triangle inequality
$$
\rho(f;1) + \rho(g;1) \geq \rho(fg;1)
$$
and (2) from the symplectic invariance of $\rho$ and (3) from the
identity $\rho(\underline 0;1) = 0$.
\end{proof}

Following \cite{ElP}, we define the partial order associated to this
normal cone on $\widetilde{Ham}(M,\omega)$ by
$$
f \geq g \, \text{ on } \, \DD \,\text{ if and only if }\, fg^{-1}
\in C.
$$
The question whether this is non-trivial, i.e., satisfies the
axiom
$$
f \leq g \, \, \& \, \, g\leq f \, \text{ if and only if }\, f = g
$$
is an interesting problem to study and is related to the study of
Hamiltonian loops $h$ and the corresponding spectral invariants
$\rho(h;1)$. This is a subject of future study. Viterbo \cite{V} had
earlier introduced the notion of positive Hamiltonians and a similar
partial order for the set of compactly supported Hamiltonians on
$\R^{2n}$ and proved nondegeneracy of the partial order.

\bigskip

\section*{\bf Appendix}

\medskip

\subsection*{\it A.1. Proof of Proposition 2.7}

\medskip

We first recall the definition from \cite{CZ}, \cite{SZ} of the
Conley-Zehnder index for a path $\alpha$ lying in $\SS P^*(1)$ where
we denote
$$
\SS P^*(1) = \{ \alpha : [0,1] \to Sp(2n,\R) \mid \alpha(0) = id,
\, \det(\alpha(1) - id) \neq 0 \} \leqno (A.1)
$$
following the notation from \cite{SZ}. We denote by
$\mu_{CZ}(\alpha)$ the Conley-Zehnder index of $\alpha$ given in
\cite{SZ}.

Next we note that a given pair $[\gamma,w] \in \widetilde
\Omega_0(M)$ determines a preferred homotopy class of
trivialization of the symplectic vector bundle $\gamma^*TM$ on
$S^1 = \part D^2$ that extends to a trivialization
$$
\Phi_w: w^*TM \to D^2 \times (\R^{2n},\omega_0)
$$
over $D^2$ of where $D^2 \subset \C$ is the unit disc with the
standard orientation. Any one-periodic solution $z: \R/\Z \to M$
of $\dot x = X_H(x)$ has the form $z(t) = \phi_H^t(p)$ for a fixed
point $p = z(0) \in \text{Fix}(\phi_H^1)$. For the given
one-periodic solution $z$ and its bounding disc $w:D^2 \to M$, we
consider the one-parameter family of the symplectic maps
$$
d\phi_H^t(z(0)): T_{z(0)}M \to T_{z(t)}M
$$
and define a map $\alpha_{[z,w]}: [0,1] \to Sp(2n,\R)$ by
$$
\alpha_{[z,w]}(t) = \Phi_w(z(t))\circ d\phi_H^t(z(0)) \circ
\Phi_{w}(z(0))^{-1}. \leqno (A.2)
$$
Obviously we have $\alpha_{[z,w]}(0) = id$, and nondegeneracy of
$H$ implies that
$$
\det(\alpha_{[z,w]}(1) -id) \neq 0
$$
and hence
$$
\alpha_{[z,w]} \in \SS P^*(1). \leqno (A.3)
$$
In general, according to the definition from \cite{CZ} of the
Conley-Zehnder index for a paths $\alpha$ lying in $\SS P^*(1)$, the
Conley-Zehnder index of $[z,w]$ is defined by
$$
\mu_H([z,w]) : = \mu_{CZ}(\alpha_{[z,w]}). \leqno (A.4)
$$

Now we start with the proof of Proposition 2.6. We can write
$$
z'(t) = \exp_{z(t)} \xi(t), \quad \xi(t) \in T_{z(t)}M
$$
since $z, \, z'$ is an associated pair of $H, \, H'$ that are
assumed to be sufficiently $C^2$ close. We denote by
$$
\Pi_{0,t}^{s,t}: T_{z(t)}M \to T_{z^s(t)}M
$$
the parallel translation along the short geodesics $r \mapsto
\exp_{z(t)}(r \xi(t))$ for $0 \leq r \leq s$ where $z^s: S^1 \to
M$ are the loops defined by
$$
z^s(t) = \exp_{z(t)}(s \xi(t)), \quad s \in [0,1].
$$
Now we make an identification
$$
\Pi: u^*TM \to [0,1] \times z^*TM
$$
for $u = u^{can}_{zz'}$ by the map
$$
v \in (u^*TM)_{(s,t)} \mapsto (s, (\Pi_{0,t}^{s,t})^{-1}(v)) \in
\{s\} \times (z^*TM)_t. \leqno (A.5)
$$
Then noting that we can write
$$
(w')^*TM = w^*TM \# u^*TM
$$
for $w' = w\# u$, $\Phi_w$ and $\Pi$ together induce a natural
trivialization
$$
\Phi_{w'}: (w')^*TM \to (D^2 \cup [0,1] \times S^1) \times \R^{2n}
$$
by the formula
$$
\Phi_{w'}(v) = \begin{cases} \Phi_w(\xi)  & \quad \xi \in w^*TM\\
(\Phi_w|_{\part D^2})(\Phi_u(\xi)) & \quad \xi \in u^*TM
\end{cases}
\leqno (A.6)
$$
By the nondegeneracy hypothesis of $z$, we know that the map
$\alpha_{[z,w]}: [0,1] \to Sp(2n,\R)$ defines a path in $\SS
P^*(1)$. We compare the two linearized vector fields of $X_H$ and
$X_{H'}$ along the corresponding periodic orbits $z$ and $z'$
$$
(\Phi_w|_{\part D})_*(DX_H(z)) \leqno (A.7)
$$
and
$$
(\Phi_{w'}|_{\part D})_*(DX_{H'}(z'))
$$
respectively. We can express
$$
\begin{aligned} (\Phi_{w'}|_{\part D})_*(DX_{H'}(z')) & = (\Phi_w|_{\part
D^2}) \circ (\Pi_0^1)^{-1} \circ (DX_{H'}(z')) \circ \Pi_0^1 \circ
(\Phi_w|_{\part D^2})^{-1} \\
& = (\Phi_w|_{\part D^2})\circ \Big( (\Pi_0^1)^{-1} \circ
DX_{H'}(z') \circ \Pi_0^1 \Big) \circ (\Phi_w|_{\part D^2})^{-1}.
\end{aligned}
$$
Since we assume that $H'$ is $C^2$ close to $H$ and $(z,  z')$ is
an associated pair, it can be easily seen
$$
(\Pi_0^1)^{-1} \circ DX_{H'}(z') \circ \Pi_0^1: z^*TM \to z^*TM
$$
is $C^1$ close to $DX_H(z): z^*TM \to z^*TM$. Therefore if we
write the flow of the linearization
$$
\Phi_{w'}(z'(t))\circ d\phi_{H'}^t(z'(0)) \circ
\Phi_{w'}(z'(0))^{-1}
$$
of $(\Phi_{w'})_*(DX_{H'}(z'))(1,t,\xi)$ as
$$
(1, t, \alpha'_{[z,w]}(t) \, \xi)
$$
the map $\alpha'_{[z,w']}: [0,1] \to Sp(2n,\R)$ is $C^1$ close to
$\alpha_{[z,w]}$. Therefore we can homotope $\alpha'_{[z',w']}$ to
$\alpha_{[z,w]}$ inside $\SS P^*(1)$. Then by an invariance property
of the Conley-Zehnder index \cite{SZ} under such a homotopy, we
obtain
$$
\mu(\alpha_{[z,w]}) = \mu(\alpha'_{[z',w']}).
$$
By the definition of the Conley-Zehnder index $\mu_H([z,w])$, this
implies $\mu_H{[z,w]}$ = $\mu_{H'}([z',w'])$ which finishes the
proof.
\medskip

\subsection*{\it A.2. Proof of Lemma 6.8}

\medskip

We first prove the following general lemma. The lemma is stated in
a more general form than needed for the purpose of using it also
in the proof of Proposition 9.2 in section 9. The present case
corresponds to the case of
$$
(H_\alpha,J_\alpha) = (H_\beta,J_\beta) = (H,J)\equiv (\HH,j).
$$

\smallskip

{\noindent\bf Lemma A.1.} {\it Let $(H,J)$ be an pair consisting of
smooth one-periodic families of Hamiltonians and almost complex
structures. Denote $(J_\alpha, H_\alpha)$ and $(J_\beta, H_\beta)$
be any regular pairs in the Floer theoretic sense for which the
Floer homologies are defined. Let $(\HH,j)$ be a homotopy pair with
$\HH$ connecting $H_\alpha, \, H_\beta$ for which the Floer chain
map $h_{(j,\HH';\rho)}$ is defined for a cut-off function. Let $U$
be any fixed tubular neighborhood of $z$ with smooth boundary $\part
U$ which is homologous to the one-cycle $[z]$. Let $z_\alpha, \,
z_\beta$ be one-periodic orbits of Hamilton's equation for $H_1, \,
H_2$ respectively such that the images of $z_\alpha$ and $z_\beta$
are contained in an open subset
$$
V \subset \overline V \subset U
$$
and denote by $u_{\alpha\beta}^{can}$ be the canonical thin
cylinder. Suppose that $u$ is a solution of $(5.3)$ satisfying $u
\sim u_{\alpha\beta}^{can}$
$$
u(-\infty) = z_\alpha, \quad u(\infty) = z_\beta.
$$
Then there exists a positive constant $\delta_0
> 0$ and $e=e(z; \delta_0,V \subset U) > 0$, which does not depend on $u$,
such that for any finite energy solution $u$ of $(5.3)$-$(5.4)$
whose image is not contained in $U$,  we have
$$
E_{(j, \HH, \rho)}(u) \geq e(z; \delta_0, V\subset U) > 0
$$
for any $H$ and $\HH$ as long as $H_m, \, m = 1, \, 2$ satisfy
$$
\|H_m - H\|_{C^2} < \delta_0, \quad \|\HH \|_{C^2} < \delta_0, \quad
\|j' - j\|_{C^1} < \delta_0
$$
}
\begin{proof} We will prove the lemma by contradiction. Suppose to
the contrary that there exists a sequence $\delta_k \to 0$ such that
we can find $H_{\alpha,k}, \, H_{\beta,k}$ and $(j_k,\HH_k)$
satisfying
$$
\|H - H_{\alpha,k}\|_{C^2}, \quad \|H - H_{\beta,k}\|_{C^2}, \quad
\|\HH_k - \HH \|, \quad \|j_k - j\| < \delta_k
$$
for which there exist periodic orbits $z_{\alpha,k} \in
\text{Per}(H_{\alpha,k})$ and $z_{\beta,k} \in
\text{Per}(H_{\beta,k})$ satisfying
$$
\text{Im }z_{\alpha,k}, \quad \text{Im }z_{\beta,k} \subset V,
$$
and a solution $u_k$ of (5.3)-(5.4) satisfying
$$
u_k(-\infty) = z_{\alpha,k}, \quad u(\infty) = z_{\beta,k} \leqno
(A.8)
$$
and
$$
u(\tau_k,t_k) \not \in U. \leqno (A.9)
$$
We recall the energy bound (5.21)
$$
E_{(j_k,\HH_k;\rho)}(u_k) \leq \AA_{H_{\beta,k}}([z_{\beta,k}, w_k
\# u]) - \AA_{H_{\alpha,k}}([z_{\alpha,k}, w_k]) + E^-(H_{\beta,k} -
H_{\alpha,k})
$$
for any bounding disc $w_k$ of $z_{\alpha,k}$. Since $u$ is assumed
to be homotopic to the canonical thin cylinder
$u_{\alpha\beta}^{can}$, we have
$$
\AA_{H_{\beta,k}}([z_{\beta,k}, w_k \# u]) =
\AA_{H_{\beta,k}}([z_{\beta,k}, w_k \# u^{can}_{\alpha\beta}])
$$
and hence this bound can be re-written as
$$
E_{(j_k,\HH_k;\rho)}(u_k) \leq \AA_{H_{\beta,k}}([z_{\beta,k}, w_k
\# u_{\alpha\beta}^{can}]) - \AA_{H_{\alpha,k}}([z_{\alpha,k}, w_k])
+ E^-(H_{\beta,k} - H_{\alpha,k}). \leqno (A.10)
$$
Using (A.8)-(A.10) and Lemma 2.4, we derive
$$
E_{(j_k,\HH_k;\rho)}(u_k) \to 0 \quad \text{as } \quad k \to \infty.
$$
In particular, $u_k$ can not bubble off and uniformly converges to a
stationary solution. We denote the stationary solution by
$z_\infty$, which will be a periodic orbit of $H$. Due to (A.9) and
(A.10), the image of $z_\infty$ must be contained in $\overline V$.
On the other hand, because of (A.11), the image of the limit of
$u_k$ cannot be contained in $U$ and gives rise to a contradiction.
This finishes the proof.
\end{proof}

Now we prove Lemma 6.8 by contradiction. Suppose the contrary that
there exists some sequence of Floer points $\eta_i$ with
$\text{dist}(\eta_i,\SS ing(\HH)) \to 0$ and a sequence of elements
$$
u_i \in \MM\Big(J(\eta_i),H(\eta_i);[z^+(\eta_i),w^+(\eta_i)],
[z_i,w_i]\Big)
$$
with $[z_i,w_i] \neq [z^-(\eta_i),w^+(\eta_i) \# u^{can}_{z^+z^-}]$
and
$$
E_{(J(\eta_i),H(\eta_i))}(u_i) \to 0. \leqno (A.11)
$$
After choosing a subsequence, we may assume
$$
\eta_i \to \eta_\infty \in \SS ing(\HH)
$$
and
$$
z_i \to z_\infty \in \text{Per}(H(\eta_\infty))
$$
as $i \to \infty$.

First if $z_i \neq z^\pm(\eta_i)$, then we have
$$
d_{C^1}(z_i,z^\pm(\eta_i)) \geq C
$$
for some $C>0$ independent of $i \to \infty$ because we are assuming
$\HH$ is a Cerf homotopy and so each singular point $\eta_\infty$
can contain only one bifurcation orbit. This contradicts to (A.11)
via Lemma A.3.

On the other hand, if $z_i = z^+(\eta_i)$, (A.11) implies $\text{Im
}u_i$ is contained in a small tubular neighborhood of $z_\infty$ and
so the compactified cycle $u_i$ are all homologous to the one
dimensional cycle $z_i$. But then Lemma 6.3 implies
$$
E_{(J(\eta_i),H(\eta_i))}(u_i) = 0
$$
for all sufficiently large $i$, after choosing a subsequence if
necessary, which contradicts that $u_i$ cannot be stationary.

Therefore we have proved
$$
z_i = z^-(\eta_i)
$$
for all $i$, after choosing a subsequence if necessary. Then we pick
a pair of sufficiently small tubular neighborhood $U$ of $\text{Im}
z_\infty$ such that it deformation retracts to the one dimensional
cycle $z_\infty$. Since $z_i \to z_\infty$, there exists another
smaller open neighborhood $V$ with $\overline V \subset U$, which
depends only on $z_\infty$ and contains the image of $z_i$ for all
sufficiently large $i$. We can also make the energy in (A.11)
satisfy
$$
\begin{aligned} E_{(J(\eta_i),H(\eta_i))}(u_i) < e(z_\infty; \delta_0, V
\subset U)
\end{aligned}
$$
by choosing $i$ sufficiently large. Then Lemma A.3 implies that the
image of $u_i$ must be contained in $U$. In particular, $z^+(\eta_i)
= u_i(-\infty)$ itself must be contained in $U$. In particular,
after compactifying $\R \cup\{-\infty, \infty \}$ to $[0,1]$, we
prove that $u_i$ must be homotopic to $u^{can}_{z^+z^-}$. This
proves $$[z_i,w_i] = [z^-(\eta_i),w^+(\eta_i) \# u^{can}_{z^+z^-}]$$
which is exactly what we wanted to prove. This proves Lemma 6.8.

\end{document}